\documentclass{article}
\usepackage[T2C]{fontenc}
\usepackage[russian,english]{babel}
\usepackage[tbtags]{amsmath}
\usepackage{amsfonts,amssymb,mathrsfs,amscd,longtable}
\usepackage{mathrsfs}
\usepackage[im]{sks2e-tr}
\usepackage[all]{xy}
\usepackage{hypbmsec}
\usepackage{mathtext}
\usepackage{hyperref}

\def\bad{\spaceskip=0.33emplus0.6emminus0.15em\immediate\write5{\string\bad}}

\issvol{86}
\rusissnum{4} 
\rusisspages{3--50}
\isspages{?--?}
\issyear{2022}
\numberwithin{equation}{section}

\theoremstyle{plain}
\newtheorem{theorem}{Theorem}
\newtheorem{lemma}{Lemma}
\newtheorem{corollary}{Corollary}
\newtheorem{proposition}{Proposition}
\theoremstyle{definition}

\newtheorem{proof}{Proof}
\newtheorem{remark}{Remark}

\newtheorem{convention}{Convention}

\begin{document}

\title({Canonical form of $C^*$-algebra of eikonals related to the metric graph}){{\bad Canonical form of $C^*$-algebra of eikonals related to the metric graph}}

\doi{10.1070/im9179} 

\author{M.~I.~Belishev}[Mikhail~I.~Belishev]%
\address{St. Petersburg Department\\
of Steklov Mathematical Institute\\
of Russian Academy of Sciences }

\email{belishev@pdmi.ras.ru}

\author{A.~V.~Kaplun}[Aleksandr~V.~Kaplun]%
\address{St. Petersburg Department\\
of Steklov Mathematical Institute\\
of Russian Academy of Sciences}

\email{alex.v.kaplun@gmail.com}

\date{22.04.2021\\09.10.2021}
\udk{517.538}

\maketitle

\markright{Canonical form of $C^*$-algebra of eikonals related to the metric graph}

\begin{fulltext}

\footnotetext[0] {This work was supported by the Russian Foundation for Basic Research (grant \No~20-01-00627A), the Volkswagen Stiftung Foundation, the Leonhard Euler Mathematical Institute (agreement \No~075-15-2019-1620), and in part by the Russian Young Mathematicians Contest.}

\begin{abstract}
The eikonal algebra $\mathfrak E$ of the metric graph $\Omega$ is an operator $C^*$--algebra defined by the dynamical system which describes the propagation of waves generated by sources supported in the boundary vertices of $\Omega$. This paper describes the canonical block form of the algebra $\mathfrak E$ of an arbitrary compact connected metric graph. Passing to this form is equivalent to constructing a functional model which realizes $\mathfrak E$ as an algebra of continuous matrix-valued functions on its spectrum $\widehat{\mathfrak{E}}$. The results are intended to be used in the inverse problem of reconstruction of the graph by spectral and dynamical boundary data.

Bibliography: 28 items.

\end{abstract}

\textbf{Published in Izvestiya RAN : Ser. Mat., Vol. 86, No 4, 3--50,
2022\,\,\,\\
(in Russian).  DOI: \href{https://doi.org/10.4213/im9179}{10.4213/im9179}}

\section{Introduction}

\label{s1}

\subsection{About work}
\label{ss1.1}

There is an approach to the inverse problems of mathematical physics  --- the \textit{boundary control method} (BC-method) \cite{1}. The approach has a pronounced interdisciplinary character: it is based on connections of inverse problems with system theory and control theory and uses asymptotic methods, functional analysis, operator theory, etc. An algebraic version of the BC-method, based on connections with Banach algebras, provided a new solution to the problem of reconstruction of a Riemannian manifold by boundary data \cite{1}--\cite{3}. Our perspective goal is to apply this version to inverse problems on graphs. This work is a step in that direction.

The algebraic version is based on a fundamental fact: a topological space can be characterized by an appropriate algebra. As an example, a compact Hausdorff space $\Omega$ is determined up to homeomorphism by the continuous functions algebra $\mathfrak{A}=C(\Omega)$ (I.\, M. Gelfand, 1943). Thus the spectrum of an algebra, i.e., the set $\widehat{\mathfrak{A}}$ of its irreducible representations endowed with a suitable topology, is homeomorphic to the space: $\widehat{\mathfrak{A}}\cong\Omega$. As a consequence, by taking any copy of $\mathfrak{A}'$ of the algebra $\mathfrak{A}$ and finding its spectrum $\widehat{\mathfrak{A}'}\cong\widehat{\mathfrak{A}}\cong\Omega$, we obtain a homeomorphic copy of the space $\Omega$. According to this scheme, the reconstruction problem is solved: a copy $\mathfrak{A}'$ is extracted from the inverse problem data, and its spectrum $\widehat{\mathfrak{A}'}$ is found, which provides the solution to the reconstruction problem --- a homeomorphic copy of the manifold $\Omega$ is recovered.  The substantive part of the approach consists in finding algebra $\mathfrak{A}'$ from the known data. As the latter, the \textit{eikonal algebra} $\mathfrak{E}$, determined by the dynamical system which describes the wave propagation in $\Omega$, is used.

Variants of the BC-method for inverse problems on graphs are proposed in \cite{4}--\cite{6}. The version using eikonal algebra originated in \cite{7} and extended in \cite{8}. Our work develops this approach. Its general direction is to study the relations between the properties of the algebra $\mathfrak{E}$ (block structure, algebraic invariants, representations) and the graph geometry. A promising goal is the reconstruction of the graph from its boundary data.

\subsection{Eikonal algebra}
\label{ss1.2}

For the sake of clarity and without loss of generality, the graph $\Omega$ can be imagined as a connected compact graph in $\mathbb{R}^3$ consisting of smooth curves (edges) $\{e_1,\dots,e_l\}=E$, connected in interior vertices $\{v_1,\dots,v_m\}=V$ \footnote{Any metric graph allows such realization.}. There are boundary vertices $\{\gamma_1,\dots,\gamma_n\}=\Gamma$, with only one edge coming out of them. The metric (inner distance) in $\Omega$ is induced by the Euclidean metric of $\mathbb{R}^3$.

The edges of the graph are ``material'': oscillations (waves) propagate along them, being initiated by point sources (controls), which are placed at the boundary vertices. The waves move from the boundary with unit velocity, gradually filling the graph. The process is described by the dynamical system
\begin{alignat*}{2}
&u_{tt}-\Delta u=0 &\quad &\text{in }\mathscr{H},\quad 0<t<T,
\\
&u|_{t=0}=u_t|_{t=0}=0 &\quad &\text{in }\Omega,
\\
&u=f &\quad &\text{on }\Gamma \times [0,T],
\end{alignat*}
where $\mathscr{H}=L_2(\Omega)$, $\Delta$ is the Laplacian defined on smooth functions that satisfy the matching (Kirchhoff) conditions in the inner vertices; $f=f(\gamma,t)$ is a \textit{boundary control} of the class $L_2(\Gamma\times [0,T])=: \mathscr{F}^T$; $u=u^f(x,t)$ is a solution (\textit{wave}), $u^f(\,{\cdot}\,,t)\in\mathscr{H}$ at $0\le t\le T$.

It is possible to control waves, not from the entire boundary but from its part $\Sigma\subset\Gamma$: in this case, the controls of the class
$$
\mathscr{F}^T_\Sigma:=\{f\in\mathscr{F}^T\mid \operatorname{supp} f\subset \Sigma\times[0,T]\}=\bigoplus\sum_{\gamma\in\Sigma}\mathscr{F}^T_\gamma
$$
are used.

Each boundary vertex is associated with a family of \textit{reachable sets} $\mathscr{U}^t_\gamma: =\{u^f(\,{\cdot}\,,t)\mid f\in\mathscr{F}^T_\gamma\}$, $0\le t\le T$, and the corresponding projectors $P^t_\gamma$ in $\mathscr{H}$ on $\mathscr{U}^t$. The operator $E^T_\gamma:=\int_0^Tt\,dP^t_\gamma$ is called \textit{eikonal} corresponding to the vertex of $\gamma$. The eikonals are self-adjoint operators -- elements of the bounded operator algebra $\mathfrak{B}(\mathscr{H})$.

For $C^*$-algebra $\mathfrak{A}$ and set $S\subset\mathfrak{A}$, we denote by $\vee S$ the $C^*$-algebra generated by this set, i.e. the minimal $C^*$--subalgebra in $\mathfrak{A}$ containing $S$. \textit{Eikonal algebra} corresponding to the selected family of boundary vertices $\Sigma\subset\Gamma$, is an operator $C^*$--algebra

\begin{equation}
\label{eq1.1}
\mathfrak{E}^T_\Sigma:=\vee\{E^T_\gamma\mid\gamma\in\Sigma\}\subset\mathfrak{B}(\mathscr{H}).
\end{equation}

\subsection{Results and comments}
\label{ss1.3}

As it was established in  \cite{7}, the algebra $\mathfrak{E}^T_\Sigma$ has a block structure: it is isomorphic to some subalgebra of the algebra $\bigoplus\limits_{j=1}^J C([0,\epsilon_j];\mathbb M^{m_j})$ and differs from the latter by the existence of relations between blocks. For the simplest graphs (three-beam stars) the nature of these connections and their evolution with changing $T$ were considered in \cite{8}.

The main result of this work is a \textit{canonical block form} of the eikonal algebra. This form is, firstly, distinguished by the absence of relations between blocks and, secondly, is invariant: up to trivial transformations (block permutations, parameterization replacement, etc.), it is determined by \textit{any} copy of the algebra \eqref{eq1.1}. The latter reinforces the hope for the usefulness of $\mathfrak{E}^T_\Sigma$ in inverse problems whose data determine it up to isometry.

  Of course, the effectiveness of the approach can be fully judged by specific applications to the inverse problems. There are no such applications in this paper but the above results seem to be an important step in this direction.

  The eikonal algebra belongs to the class of $C^*$-algebras with finite-dimensional representations of \textit{different} dimensions \cite{9}, \cite{10}. The reduction to a  canonical form is equivalent to constructing a functional model which realizes $\mathfrak{E}^T_\Sigma$ as an algebra of continuous matrix-valued functions on its spectrum. It belongs to the type of models described in the work of N.B. Vasilyev \cite{9}.

The definition of \eqref{eq1.1} reproduces the definition of the corresponding algebras used in \cite{1}--\cite{3} for the reconstruction of manifolds. The success of such an application motivated the attempt to transfer the approach to the problems on graphs. The obstacle is the non-commutativity of the algebra $\mathfrak{E}^T_\Sigma$. This problem has been encountered before in the inverse problem of electrodynamics \cite{3}, but there it was solved by factorization by the ideal of compact operators, which reduced the case to the  commutative $C(\Omega; \mathbb R)$. The noncommutativity $\mathfrak{E}^T_\Sigma$ is irremovable, which makes the study of it much more difficult.

Inverse problems on graphs are quite a relevant topic. Different formulations and approaches are contained in the works of S.A. Avdonin, P.B. Kurasov, M. Novachik, A.S.\ and V.S. Mikhailov, P.A. Kuchment, V.A. Yurko. In  \cite{11}--\cite{17} BC--method as well as other approaches are used to solve dynamic and spectral inverse problems for various classes of graphs. The works of V.\,A. Yurko and his followers use the spectral approach for inverse problems for differential operators on graphs \cite{18}--\cite{20}. Let us mention an informative review by P. Kuchment and G. Berkolayko on the whole subject of quantum graphs, including inverse problems on them.

The rather volumable introductory part of the paper essentially repeats the corresponding sections from \cite{7} and \cite{8}. This is unavoidable since the presentation of facts and results related to $\mathfrak{E}^T_\Sigma$ requires solid preparation. The technical part is rather complicated since the work deals with a maximally general object -- an arbitrary compact connected metric graph having a boundary.  The simple examples in  \cite{7}, \cite{8} with illustrations can be a useful tool for understanding. We also recommend \cite{22}, which describes in detail the procedure of reducing  $\mathfrak{E}^T_\Sigma$ to a canonical form for a simple graph.

This paper is addressed to specialists in $C^*$-algebras with a taste for applications and/or specialists in the field of mathematical physics who share the idea of the usefulness of abstractions. The eikonal algebra is a complex and property-rich object worthy of comprehensive research. For inverse problems, the possible relations between its algebraic invariants and graph geometry are particularly interesting. As an example, let us mention the hypothesis about the correspondence of clusters in the spectrum $\mathfrak{E}^T_\Sigma$ to inner vertices of $\Omega$. It is very interesting how the presence of cycles in the graph affects the structure of the spectrum.

\section{Waves on graph}
\label{s2}

\subsection{Graph}
\label{ss2.1}

Let $\Omega=E\cup W$ be a connected compact graph in $\mathbb{R}^3$ with  edges $\{e_1,\dots,e_L\}=E$ and vertices $\{w_1,\dots,w_M\}=W$. Edges are smooth\footnote{Everywhere in the paper, \textit{smooth} means $C^\infty$-smooth.} curves with ends serving as vertices. It is convenient to think of edges as open without including their ends. A vertex $w$ and an edge $e$ are \textit{incidental} (we write $w\prec e$) if $w$ is the end of $e$. The vertices $\{\gamma_1,\dots,\gamma_N\}=\Gamma$, which are all incident to a single edge, are called boundary vertices; the vertices $\{v_1,\dots,v_{M-N}\}=V=W\setminus\Gamma$ are the interior ones.

The number $\mu(w)$ of edges that are incident to vertex $w$ is called its valence; for $\gamma\in \Gamma$ we have $\mu(\gamma)=1$. Additionally, we assume that there are no vertices with $\mu(w)=2$, so that $\mu(v)\ge 3$ is satisfied in all interior vertices.

The graph is equipped with a metric (internal distance) $\tau$ induced by the Eucli\-dean metric of $\mathbb{R}^3$. Thus $\tau(a,b)$ is the minimum of lengths of piecewise smooth curves lying in $\Omega$ and connecting points $a$ and $b$.  For the set $A\subset\Omega$ its metric neighborhood of radius $r$ is denoted by
$$
\Omega^r_A:=\{x\in\Omega\mid \tau(x,A)<r\},\qquad r>0.
$$

Each edge $e$ is parameterized by the length $\tau$ counted from one of its ends. For the function $y$ on the graph, the sign of the derivative w.r.t. the length $dy/d\tau$ depends on the choice of the end, but the second derivative $d^2y/d\tau^2$ does not depend on this choice. For the vertex $w$ and its incident edge $e$, the derivative in the direction coming from $w$ is defined
$$
\biggl[\frac{dy}{d\tau}\biggr]^+_e(w):=\lim_{e\ni x\to w}\frac{dy}{d\tau}(x)
$$
The value is
$$
F_w[y]:=\sum_{e\succ w}\biggl[\frac{dy}{d\tau}\biggr]^+_e(w)
$$
is called the flow of function $y$ through vertex $w$.

The metric $\tau$ on a graph defines a (real) Hilbert space $\mathscr{ H}=L_2(\Omega)$ with scalar product:
$$
(y,u)_\mathscr{H}:=\int _\Omega y u\,d\tau = \sum_{e \in E} \int _e y u\,d\tau.
$$
Let $C(\Omega)$ be the space of continuous functions with norm $\|y\||\,{=}\sup_\Omega|y(\cdot)|$. A function $y$ is assigned to the Sobolev class $H^2(\Omega)$ if $y\in C(\Omega)$ and $dy/d\tau,d^2y/d\tau^2\in L_2(e)$ is satisfied on each edge.

Define the Kirchhoff class
$$
\mathscr{K} := \{y \in {H}^2(\Omega)\mid F_v[y]=0,\, v\in V \}
$$
The Laplace operator on a graph is given by
\begin{equation}
\label{eq2.1}
\Delta\colon \mathscr{H}\to\mathscr{H}, \quad \operatorname{Dom}\Delta=\mathscr{K},\qquad (\Delta y)|_e=\frac{d^2y}{d\tau^2}, \quad e \in E.
\end{equation}
It is densely defined and closed.

\subsection{Waves}
\label{ss2.2}
{\bad The initial boundary value problem describing wave propagation} in the graph, is of the form
\begin{alignat}{2}
\label{eq2.2}
&u_{tt}-\Delta u=0 &\quad &\text{in }\mathscr{H},\quad 0<t<T,
\\
\label{eq2.3}
&u(\,{\cdot}\,,t)\in \mathscr{K} &\quad &\text{at }0\le t\le T,
\\
\label{eq2.4}
&u|_{t=0}=u_t|_{t=0}=0 &\quad &\text{in }\Omega,
\\
\label{eq2.5}
&u=f &\quad &\text{on }\Gamma \times [0,T].
\end{alignat}
Here $T >0$ is the final moment of time, $f=f(\gamma ,t)$ -- \textit{boundary control}, $u=u^f(x,t)$ -- solution (wave). For a smooth (by $t$) control $f$ vanishing near $t=0$, the problem has a unique classical solution $u^f$.

By the definition \eqref{eq2.1}, on each edge $e$ the solution $u^f$ satisfies the homogeneous string equation $u_{tt}-u_{\tau \tau}=0$. From this, we see that the waves propagate from the boundary into $\Omega$ with unit velocity. As a consequence, if the control acts from the part of the boundary $\Sigma\subseteq\Gamma$, i.e., $\operatorname{supp}f\subset\Sigma\times[0,T]$ is satisfied, we have the relation
\begin{equation}
\label{eq2.6}
\operatorname{supp}u^f(\,{\cdot}\,,t) \subset \overline{\Omega^t_\Sigma}, \qquad t>0.
\end{equation}
In the following, we introduce (generalized) solutions of the problem \eqref{eq2.2}--\eqref{eq2.5} for controls of class $L_2(\Gamma \times [0,T])$. Their definition requires some preparation.

With $\delta_x$ we denote the Dirac measure, i.e., a functional on $C(\Omega)$ taking values according to the rule $\langle\delta_{x},y\rangle=y(x)$. Let $\delta=\delta(t)$ be the Dirac delta--function. The nearest goal is to define and describe the \textit{fundamental solution} of the problem \eqref{eq2.2}--\eqref{eq2.5} with $T=\infty$, corresponding to the control $f=\delta_\gamma\delta(t)$, which acts instantaneously from the boundary vertex $\gamma$.

To describe the solution of $u^{\delta_\gamma \delta}$ it is convenient to use the following formalism of ``impulse dynamics''.

0. The \textit{impulse} corresponds to the measure $a\delta_x$; the constant $a\ne 0$ is called its amplitude.

1. Each impulse $a\delta_{x(t)}$ moves along the edge with velocity $1$ in  one of two possible directions, so that $|\dot x(t)|=1$ holds for $x(t)\in e$.

2 (superposition principle). Impulses move independently of each other. If at   the moment $t$ there are several impulses $a_1\delta_{x(t)},\dots,a_p\delta_{x(t)}$ located at   the point $x(t)\in\Omega\setminus\Gamma$, they are combined to form the impulse $[\,a_1+\dots + a_p\,]\delta_{x(t)}$.

3 (passing through the interior vertex). Moving along the edge $e$ and passing through the inner vertex $v$, the impulse $a\delta_{x(t)}$ divides into $\mu(v)$ impulses: one reflected and $\mu(v)-1$ passed. The reflected impulse moves along $e$ in the opposite direction and has an amplitude of $(2-\mu(v)/\mu(v))a$. Each of the passed impulses moves along its (incidental to $v$) edge away from $v$ and has an amplitude of $(2/\mu(v))a$.

Therefore, the total amplitude is
$$
\frac{2-\mu(v)}{\mu(v)}\,a+[\mu(v)-1]\frac{2}{\mu(v)}\,a=a,
$$
which corresponds to the Kirchhoff law of conservation of flows $F_v[y]=0$.

4 (reflection from the boundary). As soon as the impulse $a\delta_{x(t)}$ reaches vertex $\gamma\,{\in}\,\Gamma$, it instantly inverts its direction and changes its amplitude from  $a$ to $-a$.

Adopting these rules, we can describe the solution $u^{\delta_\gamma\delta}$ as follows (recall that $\tau$ is the distance in $\Omega$):

$(*)$ at $0\le t\le\tau(\gamma,V)$ we have $u^{\delta_\gamma\delta}=\delta_{x(t)}$, where $x(t)$ --- the point of edge $e\succ\gamma$ such that $\tau(x(t),\gamma)=t$. Thus, at small times $u^{\delta_\gamma\delta}$ is a single impulse with unit amplitude entering the graph from the vertex $\gamma$ and moving along $e$ with unit velocity;

$(**)$ further evolution at times $t>\tau(\gamma,V)$ is determined by the rules 1--4.

It is easy to see that this description is quite deterministic. At each moment of time $t\ge 0$, the solution $u^{\delta_\gamma\delta}$ is a finite set of impulses moving in $\Omega$. From a physical point of view, this picture describes, for example, the propagation of sharp signals (voltage spikes) in an electrical net -- a graph made up of wires.

Thus, the fundamental solution is some space-time distribution in  $\Omega\times\{t\ge 0\}$. Its structure is such that for controls of the form $f=\delta_\gamma\varphi$ with $\varphi\in L_2[0,T]$ the time convolution
\begin{equation}
\label{eq2.7}
u^f(x,t):=\bigl[u^{\delta_\gamma\delta}(x,\,{\cdot}\,)\ast \varphi\bigr](t), \qquad x\in\Omega,\quad 0\le t\le T,
\end{equation}
is well-defined. Moreover, it can be shown that $u^f\in C([0,T];\mathscr{H})$ and, if $\varphi$ is smooth and vanishes around $t=0$, then $u^f$ provides the classical solution of the problem \eqref{eq2.2}--\eqref{eq2.5}.

Since this point, the function $u^f$ defined by the relation \eqref{eq2.7} is considered as a (generalized) solution for the control of the given class. In the more general case of controls $f\in L_2(\Gamma\times[0,T])$ of the form $f=\sum\limits_{\gamma\in\Sigma}f_\gamma$ with $f_\gamma=\delta_\gamma\varphi_\gamma$ we put
\begin{equation}
\label{eq2.8}
u^f(x,t):=\sum_{\gamma\in\Sigma}u^{f_\gamma}(x,t), \qquad x\in\Omega,\quad 0\le t\le T.
\end{equation}
It is not difficult to show that for the generalized solution the relation \eqref{eq2.6} remains valid. It shows that the metric neighborhood $\Omega^T_\Sigma$ is a part of the graph captured by waves coming from $\Sigma$, by  moment $t=T$.

\subsection{Hydra}
\label{ss2.3}

Here we introduce a space-time graph, which is used to efficiently describe waves.

Let us fix the boundary vertex $\gamma$. Considering the fundamental solution as a space-time distribution, we define the set
$$
H_\gamma:=\operatorname{supp}u^{\delta_\gamma\delta}\subset \Omega\times{\overline{\mathbb R}}_+,
$$
which we will call \textit{hydra} \cite{6}. It is essentially a space-time graph formed by impulse trajectories in the course of the evolution described by rules 1--4 and $(*)$, $(**)$ (Fig. \ref{fig1})\footnote{Illustrations borrowed from \cite{7}.}. Its edges are the characteristics of the wave equation \eqref{eq2.2}.

\begin{figure}[!htb]
\begin{center}
\includegraphics{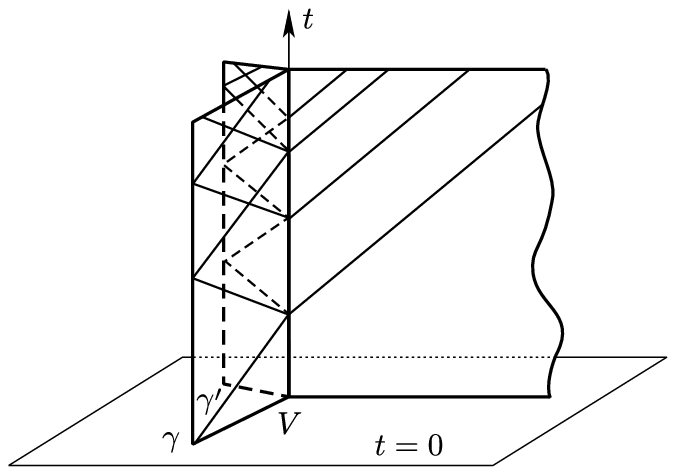}
\caption{Hydra}
\label{fig1}
\end{center}
\end{figure}

Let us define the projections
\begin{alignat*}{2}
&\pi\colon H_\gamma \ni h=(x,t)\mapsto x\in \Omega, &\qquad \pi^{-1}(x) &:=\{h\in H_\gamma\mid \pi(h)=x\};
\\
&\rho\colon H_\gamma \ni h=(x,t)\mapsto t\in \overline{\mathbb R}_+, &\qquad \rho^{-1}(t) &:=\{h\in H_\gamma\mid \rho(h)=t\}.
\end{alignat*}
On the hydra, we define the function (\textit{amplitude}) $a(\,{\cdot}\,)$ by the following rule:

-- for a point $h\in H_\gamma$ such that $\pi(h)=x\in\Omega\setminus\Gamma$ and $\rho(h)=t>0$, we have $u^{\delta_\gamma \delta}(\,{\cdot}\,,t)=a\delta_x(\,{\cdot}\,)$ and define $a(h):=a$

-- for $h\in H_\gamma$ such that $\pi(h)\in \Gamma$ and $\rho(h)>0$, put $a(h):=0$;

-- for $h\in H_\gamma$ such that $\pi(h)=\gamma$ and $\rho(h)=0$, put $a(h):=1$.

As one can see, the amplitude is a piecewise constant function defined over the entire hydra $H_\gamma$ (Fig. \ref{fig2}). Let us specify that at points of self-intersection $p$ we have $a(p)=-4/9+1/3=-1/9$ according to the rule 2 of impulse evolution. The self-intersection points like $p$ are the vertices of $H_\gamma$, which are projected to $\Omega\setminus[V\cup\Gamma]$.

\begin{figure}[!htb]
\begin{center}
\includegraphics{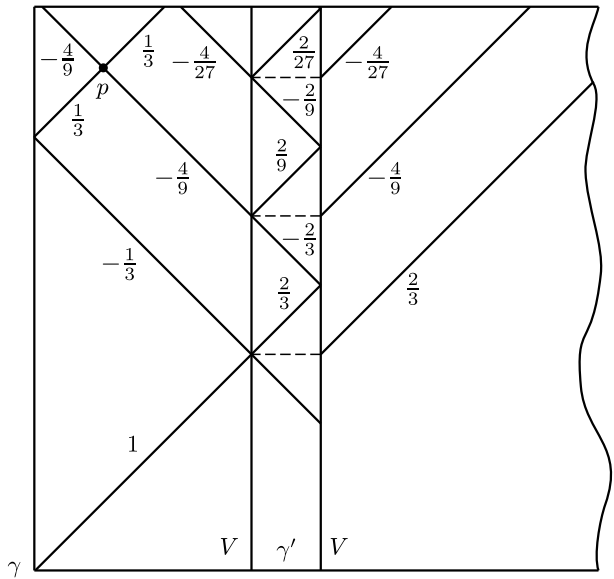}
\caption{Amplitude on the hydra}
\label{fig2}
\end{center}
\end{figure}

Let us present the representation for which the hydra was introduced (see \cite{7}). Using the notation $h=(x,t)\in H_\gamma$ and $a(h)=a(x,t)$, for the control $f=\delta_\gamma\varphi(t)$ with $\varphi\in L_2[0,T]$ according to \eqref{eq2.7} we obtain
\begin{equation}
\label{eq2.9}
u^f(x,T)=\sum_{t\in \rho(\pi^{-1}(x))}a(x,t)\,\varphi(T-t),\qquad x\in\Omega.
\end{equation}
In the general case where the control $f=\sum_{\gamma\in\Sigma}\delta_\gamma\varphi_\gamma(t)$ acts from several vertices, according to \eqref{eq2.8} we have
\begin{equation}
\label{eq2.10}
u^f(x,T)=\sum_{\gamma\in\Sigma}\,\sum_{t\in \rho(\pi^{-1}(x))}a_\gamma(x,t)\varphi_\gamma(T-t), \qquad x\in\Omega,
\end{equation}
where $a_\gamma$ are the amplitudes on the hydras $H^T_\gamma$.

These representations are quite efficient: they can be used to calculate the values of waves. However, the forthcoming analysis of the eikonal algebra will require their modification, to which we now turn. The modification uses a partitioning of the graph $\Omega$ into parts (\textit{families}) consistent with the structure of the hydra. It is described in all details and with graphic illustrations in  \cite{7}.

\subsection{Partition $\Pi$}
\label{ss2.4}

In furthermore, we are dealing with  \textit{ truncated} hydras
$$
H^T_\gamma:=H_\gamma \cap \{\Omega\times[0,T]\}.
$$

Let us introduce some general notion. Let $X$ be a set and let a reflexive symmetric (but not, generally speaking, transitive!) binary relation $\sim_0$ be given on it. The elements connected by this relation will be called \textit{neighbors}. Neighborhood determines equivalence by the following rule. We say that $x$ and $y$ are equivalent (and write $x\sim y$) if $X$ contains a finite set of elements $x_1,\dots,x_n$ such that $x\sim_0 x_1\sim_0 \dots\sim_0 x_n\sim_0 y$ holds.

By means of the relation $\sim$ introduced above, the equivalence class $[x]$ of an element $x\in X$  can be described constructively. Let us introduce an operation $\operatorname{ext}$ extending the subset $B\subset X$ by the rule
\begin{equation}
\label{eq2.11}
B\mapsto \operatorname{ext}B:=\bigcup_{b\in B}\{x\in X\mid x\sim_0 b\};
\end{equation}
denote $\operatorname{ext}^1B:=\operatorname{ext}B$ and $\operatorname{ext}^jB:=\operatorname{ext}\operatorname{ext}^{j-1}B$, $j\ge 2$. As one can easily see, the following representation holds:
\begin{equation}
\label{eq2.12}
[x]=\bigcup_{j\ge 1}\operatorname{ext}^j\{x\}.
\end{equation}
If set $X$ is finite, then the sequence $\operatorname{ext}^j$ stabilizes at some step: $\operatorname{ext}^1\{x\}\subset\dots\subset \operatorname{ext}^N\{x\}=\operatorname{ext}^{N+1}\{x\}=\dots=[x]$. This case will occur in the paper.

Consider the equivalence of this kind on the hydra. We say that the points $h,h'\in H^T_\gamma$ are neighbors ($h\overset{\gamma}\sim_0 h' $) if at least one of the conditions is satisfied: $\pi(h)=\pi(h')$ or $\rho(h)=\rho(h')$. By $\stackrel{\gamma}{\sim}$, we denote the equivalence relation, generated by such a neighborhood. The equivalence class
$$
\mathscr{L}[h]:=\{h'\in H^T_\gamma\mid h'\stackrel{\gamma}{\sim}
h\}
$$
is called a \textit{lattice}. It can be shown that this class consists of a finite number of points. For the subset $B\subset H^T_\gamma$ we define the lattice
$$
\mathscr{L}[B]:=\bigcup_{h\in B}\mathscr{L}[h].
$$
Note that the operation $B\mapsto \mathscr{L}[B]$ has the following properties:
\begin{gather*}
B\subset \mathscr{L}[B],\qquad \mathscr{L}[\mathscr{L}[B]]=\mathscr{L}[B],\qquad \mathscr{L}[B_1\cup B_2]=\mathscr{L}[B_1]\cup\mathscr{L}[B_2],
\\
\pi^{-1}(\pi(\mathscr{L}[B]))=\rho^{-1}(\rho(\mathscr{L}[B]))=\mathscr{L}[B].
\end{gather*}
The first three properties show that it is a topological closure (by Kuratowski).

To each point $x\in\overline{\Omega^T_\gamma}$ we associate the set
\begin{equation}
\label{eq2.13}
\Lambda[x]:=\pi(\mathscr{L}[\pi^{-1}(x)])\subset \overline{\Omega^T_\gamma}
\end{equation}
(the closure in the metric $\Omega$), which we will call \textit{determination set} of point $x$. This set is finite. It is easy to check that the relation
$$
x\sim x'\quad \Longleftrightarrow\quad \Lambda[x]=\Lambda[x']
$$
is equivalence, and the operation $A\mapsto\Lambda[A]:=\bigcup_{x\in A}\Lambda[x]=\pi(\mathscr L[\pi^{-1}(A)])$ is a topological closure. We say the sets $A=\Lambda[A]$ to be $\Lambda$--\textit{closed}.

On a complete hydra, a point $h\in H_\gamma$ is called \textit{corner} if $\pi(h)\in V\cup\Gamma$ or $h$ is a self-intersection point (like $p$ in fig. \ref{fig2}). The latter are the vertices of the valence hydra $4$, which are projected to $\Omega\setminus[V\cup\Gamma]$.

On the truncated hydra $H^T_\gamma$, in addition to the corner points of the complete hydra contained in it, we assign the points of the set $\rho^{-1}(T)$ to corner points. By $\operatorname{Corn}H^T_\gamma$ we denote the set of all corner points of the truncated hydra.

The lattice $\mathscr{L}[\operatorname{Corn} H^T_\gamma]$ divides the hydra into a finite number of open space-time intervals. At each interval, the amplitude $a$ takes a constant value.

The points that constitute a finite set
$$
\Theta:=\pi\bigl(\mathscr{L}[\operatorname{Corn}H^T_\gamma]\bigr)\subset \overline{\Omega^T_\gamma},
$$
are called \textit{critical}. The remaining points $x\in \overline{\Omega^T_\gamma}\setminus \Theta$ will be called \textit{regular}. The critical points divide $\overline{\Omega^T_\gamma}$ into parts. The set of regular points
$$
\Pi:=\overline{\Omega^T_\gamma}\setminus \Theta
$$
is a set of a finite number of open intervals, each of which lies on some edge $e$. Thus, $\overline{\Omega^T_\gamma}=\Pi\cup\Theta$ is a partition of the part of the graph $\Omega$ captured by waves, determined by the hydra structure $H^T_\gamma$.

Let $\omega=(c,c')\subset \Pi$ be \textit{maximal} interval consisting of regular points. Maximality means that the ends of the interval $c$ and $c'$ are critical points, so it is impossible to extend $\omega$ by keeping the interior points regular. It is easy to see that the set
\begin{equation}
\label{eq2.14}
\Phi:=\Lambda(\omega)=\pi(\mathscr{L}[\pi^{-1}(\omega)])
\end{equation}
consists of the maximum intervals$\omega_1,\dots,\omega_m$ of equal length:
$$
\Phi=\bigcup_{k=1}^m \omega_k, \qquad \operatorname{diam}\omega_1=\dots=\operatorname{diam}\omega_m=\tau(c,c')=:\epsilon_\Phi,
$$
where $\tau$ -- the distance in the graph. We say that the intervals $\omega_k$ are \textit{cells} of the family $\Phi$.

Comparing the definitions of \eqref{eq2.13} and \eqref{eq2.14}, we arrive to the  representation
\begin{equation}
\label{eq2.15}
\Phi=\bigcup_{x\in\omega}\Lambda[x],
\end{equation}
where $\omega$ is any of the cells of $\Phi$.

Let $\omega'\subset\Pi$ be the maximal interval not lying in the family $\Phi$. It defines another family $\Phi'=\Lambda[\omega']$ consisting of cells and so on. As a result, the set $\Pi$ is a finite union of non-intersecting families $\Phi^1,\dots,\Phi^J$, each family consisting of non-intersecting cells:
\begin{equation}
\label{eq2.16}
\Pi=\bigcup_{j=1}^J\Phi^j=\bigcup_{j=1}^J\bigcup_{k=1}^{m_j}\omega^j_k,
\end{equation}
where $m_j$ is the number of cells in $\Phi^j$.

\goodbreak 

In parallel to the determination set \eqref{eq2.13}, each $x\in\overline{\Omega^T_\gamma}\setminus\Gamma$ is associated with the set
$$
\Xi[x]:=\rho(\mathscr{L}[\pi^{-1}(x)])\subset [0,T].
$$
For $x\ne x'$ either $\Xi[x]=\Xi[x']$ or $\Xi[x]\cap\Xi[x']=\varnothing$ is satisfied. We also define $\Xi[B]:=\bigcup\limits_{x\in B}\Xi[x]$.

Let $\Phi=\bigcup\limits_{k=1}^{m_\Phi}\omega_k\subset\Pi$ be a family. It is easy to check that the set
\begin{equation}
\label{eq2.17}
\Psi:=\Xi[\Phi]=\bigcup_{i=1}^{n_\Phi}\psi_i\subset [0,T]
\end{equation}
consists of time intervals $\psi_i: =(t_i,\widetilde{t}_i)$ such that $0\le t_1<\widetilde{t}_1\le t_2<\widetilde{t}_2\le \dots\le t_{n_\Phi }<\widetilde{t}_{n_\Phi}\le T$ and have the same length $\widetilde{t}_i-t_{i}=\epsilon_\Phi$. The set $\Psi$ will also be called a family consisting of \textit{time} cells $\psi_i$.

In what follows, the functions $\tau^i\colon \Phi \to [0,T]$ associated with graph partitioning into families are used. They are introduced as follows \footnote{The definition given here differs from that in \cite{7}. The reason for the change will be explained later.}. For $x\in\Phi$, we set
\begin{equation}
\label{eq2.18}
\tau^i(x):=\psi_i\cap \rho(\mathscr{L}[\pi^{-1}(x)]),\qquad
i=1,\dots, n_\Phi.
\end{equation}
Since for any $x_k\in\Lambda[x]$ $\mathscr{L}[\pi^{-1}(x)]=\mathscr{L}[\pi^{-1}(x_k)]$ is satisfied, these functions are constant on determination sets: $\tau^i(x)=\tau^i(x_k)$. The definition easily implies
\begin{equation}
\label{eq2.19}
\tau^i(x)\ne \tau^{i'}(x)\quad\text{for}\quad i\ne i',\ \ x\in\Phi.
\end{equation}
Varying $x$ inside the cell $\omega=(c,c')\subset\Phi$, one varies the set $\Lambda[x]$. In this case, the value of $\tau^i(x)$ covers the cell $\psi_i=(t_i,\widetilde t_i)\subset\Psi$ and, as it is easy to see from the definition \eqref{eq2.18}, \textit{one} of two representations is valid
\begin{equation}
\label{eq2.20}
\tau^i(x)=t_i+\tau(x,c)\quad \text{or}\quad \tau^i(x)=\widetilde t_i-\tau(x,c).
\end{equation}
With this representation in mind, we can say that the functions $\tau^i$ depend linearly on $x\in\omega$.

The property \eqref{eq2.20} allows us to extend functions $\tau^i$ to critical points: if $x\in\omega=(c,c')$, $x\to c$, then $\tau^i(c)=t_i$ or $\tau^i(c)=\widetilde t_i$, depending on which representation \eqref{eq2.20} takes place.

For each family $\Phi\subset\Pi$, the set of functions $\tau^i$ is different and, if necessary, we denote it by $\tau^i_\Phi$. We add that, due to \eqref{eq2.19} and \eqref{eq2.20}, the equality $\tau^i_\Phi(x)=\tau^{i'}_{\Phi'}(x)$ for different $\Phi$ and $\Phi'$ is possible only at  critical points $x$.

The partitioning of the graph into families is motivated, in particular, by the fact that the waves $u^f$ depend on the controls $f$ locally in the following sense. As can be seen from \eqref{eq2.9}, the values of $u^f(\,{\cdot}\,,T)|_{\Phi}$ are determined by the values of $f|_{\Xi[\Phi]}$. Moreover, the conditions
\begin{equation}
\label{eq2.21}
\operatorname{supp}f\subset \Xi[\Phi]\quad\text{and}\quad\operatorname{supp}u^f(\,{\cdot}\,,T)\subset \Phi
\end{equation}
are equivalent.

The modification of the representations \eqref{eq2.9} and \eqref{eq2.10}, discussed at the end of the previous paragraph, consists of an efficient representation of waves on determination sets. It uses functions (vectors) on $\Lambda[x]$, to describe which we proceed.

\subsection{Amplitude vectors}
\label{ss2.5}

The current considerations still refer to a chosen boun\-dary vertex $\gamma\in\Gamma$, a fixed $T>0$, and the corresponding partition $\Pi$.

  The set $\{\alpha^1,\dots,\alpha^n\}$ will be called $\alpha$--\textit{set} over the set $A$. For each $x\in\in\overline{\Omega^T_\gamma}$ the determination set $\Lambda[x]$ is $\Lambda$--closed and has its $\alpha$--set over it.

It follows from \eqref{eq2.6} that the value
$$
T_\gamma:=\inf \{t>0\mid \Omega^t_\gamma=\Omega\}
$$
  is the time of filling the whole graph $\Omega$ with waves coming (with unit velocity) from the vertex $\gamma$. The lemma  \ref{l1} is valid:

\begin{lemma}
\label{l1}   
Let $x$ be a point in $\overline{\Omega^T_\gamma}$, $\Lambda[x]$ be its determination set, and let $\alpha^1,\dots,\alpha^n$ be $\alpha$--set over $\Lambda[x]$. If $T<T_\gamma$, then the vectors $\alpha^i$ are linearly independent.
\end{lemma}

\begin{proof}
  Fix some point in time $t^*=t_i\in\Xi[x]$ with $i\neq 1$. In the lattice $\mathscr{L}[\pi^{-1}(x)]$ there is necessarily a point $(x^*,t^*)$ such that that $x^*$ belongs to the boundary $\overline{\Omega^{t^*}_\gamma}\setminus{\Omega^{t^*}_\gamma}$ region of ${\Omega^{t^*}[\gamma]}$ captured by waves by the moment $t=t^*$. This boundary is \textit{nonempty} due to $t^*\le T<T_\gamma$; thus $x^*\in\Lambda[x]$, $t^*=\tau(x^*,\gamma)$ and obviously $a(x^*,t^*)=\alpha^i(x^*)\ne 0$.

  At the same time, at all points $(x^*,t)$ with $t\in\Xi[x^*]=\Xi[x]$, $t<t^*$, $a(x^*,t)=0$ holds since such $(x^*,t)$ do not lie on hydra $H^T_\gamma$. The latter corresponds to the simple fact that at the given times $t$ the waves from $\gamma$ do not have time to reach the point $x^*$ (see \eqref{eq2.6}).

  Thus, at points $(x^*,t)$ with $t=t_k\in\Xi[x]$, $t_i<t^*$\enskip ($k=1,\dots,{i-1}$) we have $a(x^*,t_k)=\alpha^k(x^*)=0$ and, at the same time, $a(x^*,t^*)=\alpha^i(x^*)\ne 0$. This eliminates the linear dependence of $\alpha^i$ on $\alpha^1,\dots,\alpha^{i-1}$. Due to the arbitrariness of $i$, we obtain linear independence of the entire set $\alpha^1,\dots,\alpha^n$. The lemma is proved.
\end{proof}

\begin{corollary}
\label{c1}
If $T<T_\gamma$, the equality $\Xi[\overline{\Omega^T_\gamma}]=[0,T]$ holds.
\end{corollary}

  In fact, the proof of the lemma \ref{l1} shows that the condition $T<T_\gamma$ ensures that the set $\rho^{-1}(t)\subset H^T_\gamma$ is nonempty at $0\le t\le T$, which is equivalent to the equality $\Xi[\overline{\Omega^T_\gamma}]=[0,T]$.

Let's return to the representation \eqref{eq2.9}. In terms of amplitude vectors, it can be written in the form of
\begin{equation}
\label{eq2.22}
u^f(x_k,T)|_{x_k\in\Lambda[x]}=\sum_{i=1}^{n[x]}\varphi(T-t_i)\alpha^i(x_k),\quad \text{for}\quad f=\delta_\gamma\varphi,\ \ n[x]:=\sharp\,\Xi[x],
\end{equation}
  representing the wave not only at  point $x$ but also on the whole determination set $\Lambda[x]$. Moreover, according to \eqref{eq2.15}, by varying the point $x$ inside the cell $\omega\subset\Phi$, we represent the wave $u^f(\,{\cdot}\,,T)$ on the whole family $\Phi$.

The final step to modify the original representation \eqref{eq2.9} is to transition in \eqref{eq2.22} to a more convenient system of amplitude vectors.

\subsection{$\beta$-representation of waves}
\label{ss2.6}

  Again, let $A=\Lambda[A]=\{x_1,\dots,x_m\}\subset\overline{\Omega^T_\gamma}$. Let us introduce the space $\mathbf{l}_2(A)$ of functions (vectors) on $A$ with scalar product
$$
\langle f,g\rangle=\sum_{x\in A}f(x)g(x)=\sum_{k=1}^{m} f(x_k)g(x_k).
$$
It contains a subspace
$$
\mathbb{A}[A]:=\operatorname{span}\{\alpha^1,\dots,\alpha^{n}\},\qquad \operatorname{dim}\mathbb{A}[A]\le n,
$$
defined by the $\alpha$-set over $A$. Using the Gram-Schmidt procedure, we pass in $\mathbb{A}[A]$ to the set
$$
\beta^i:= \begin{cases}
\dfrac{\alpha^1}{\|\alpha^1\|}, &\text{if }i=1,
\\[3mm]
\dfrac{\alpha^i-\sum _{j=1}^{i-1}\langle \alpha^i,\beta^j\rangle
\beta^j}{\|\alpha^i-\sum_{j=1}^{i-1}\langle
\alpha^i,\beta^j\rangle \beta^j\|}, &\text{if }i\ge 2\text{ and
}\alpha^i \notin \operatorname{span}\{\alpha^1, \dots,
\alpha^{i-1}\},
\\[2mm]
0, &\text{if }\alpha^i \in \operatorname{span}\{\alpha^1, \dots, \alpha^{i-1}\},
\end{cases}
$$
  where $\beta^i$ -- a vector with components $\beta^i(x_1),\dots,\beta^i(x_m)$. For nonzero elements of the set, $\langle \beta^i,\beta^j \rangle=\delta_{ij}$ holds, whereas their linear span obviously coincides with $\mathbb{A}[A]$. The set of vectors $\{\beta^1,\dots,\beta^n\}$ will be called $\beta$-set over the set $A$.

  For each $x\in\overline{\Omega^T_\gamma}$ over the determination set $\Lambda[x]$, there is a different $\beta$--set. According to the lemma, \ref{l1} for $T<T_\gamma$ the vectors $\alpha^i$ are linearly independent. As a consequence, for such $T$ all vectors $\beta^i$ are nonzero and $\operatorname{dim}\mathbb{A}[\Lambda[x]]=n\le m$ is valid.

With the use of the $\beta$-set over the set $A=\Lambda[x]$, the representation \eqref{eq2.22} takes the final form
\begin{equation}
\label{eq2.23}
u^f(x_k,T)|_{x_k\in\Lambda[x]}=\sum_{j=1}^{n[x]}c^\varphi_j\beta^j(x_k),\qquad
c^\varphi_j=\sum_{i=1}^{n[x]}\rho_{ji}\,\varphi(T-t_i),
\end{equation}
where $\rho$ is the transition matrix connecting the sets $\alpha$ and $\beta$.

\subsection{Hydra $H^T_\Sigma$}
\label{ss2.7}

  The concepts and objects introduced above correspond to a single boundary vertex $\gamma$. In what follows, we use the notation $\overset{\gamma}\sim$, $\mathscr L_\gamma$, $\Lambda_\gamma$, etc., to indicate this if necessary.

For the set of boundary vertices $\Sigma\subseteq\Gamma$ we define a space-time graph
$$
H^T_\Sigma:=\bigcup_{\gamma\in\Sigma}H^T_\gamma\subset\Omega\times[0,T].
$$
On it, the analogs of the objects introduced earlier for the single hydras $H^T_\gamma$ are defined. Let's describe them.

  The projections from $H^T_\Sigma$ to $\overline{\Omega^T_\Sigma}$ and to $[0,T]$ are $\pi((x,t)):=x$ and $\rho((x,t)):=t$; by $\pi^{-1}$ and $\rho^{-1}$ we mean the full preimages in $H^T_\Sigma$.

  By definition, the neighborhood $h\sim_0 h'$ on $H^T_\Sigma$ means that $\pi(x)=\pi(x')$ and/or $\rho(x)=\rho(x')$ holds. Neighborhood determines the equivalence of $h\overset{\Sigma}\sim h'$. By $\mathscr L_\Sigma[h]$ we denote the equivalence class (\textit{lattice}) of the point $h\in H^T_\Sigma$. The operation $H^T_\Sigma\supset B\mapsto \mathscr L_\Sigma[B]$ is a (topological) closure.

  The set $\Lambda_\Sigma[x]:=\pi(\mathscr{L}_\Sigma[\pi^{-1}(x)])$ will be called \textit{determination set} of $xin\overline{\Omega^T_\Sigma}$. Note the obvious embedding $\Lambda_\gamma[x]\subset\Lambda_\Sigma[x]$ for $\gamma\in\Sigma$. The operation $\overline{\Omega^T_\Sigma}\supset A\mapsto \Lambda_\Sigma[A]$ is a (topological) closure.

  The set of corner points $\operatorname{Corn}H^T_\Sigma$ consists of all corner points of hydras $H^T_\gamma\subset H^T_\Sigma$ plus points of (transversal) intersection of edges of different hydras $H^T_\gamma$. The critical points in $\overline{\Omega^T_\Sigma}$ are $\Theta_\Sigma:=\pi(\mathscr L_\Sigma[\operatorname{Corn}H^T_\Sigma])$, the regular points are $\Pi_\Sigma:=\overline{\Omega^T_\Sigma}\setminus\Theta_\Sigma$. There is a partitioning
$$
\Pi_\Sigma=\bigcup_{j=1}^J\Phi^j=\bigcup_{j=1}^J\bigcup_{k=1}^{m_j}\omega^j_k,\qquad \operatorname{diam}\omega^j_k=\epsilon_j:=\epsilon_\Phi,
$$
  into families and cells, which is quite similar to the partitioning \eqref{eq2.16}. For each family $\Phi\subset\Pi_\Sigma$, the set
$$
\Xi_\Sigma[\Phi]:=\rho(\mathscr L_\Sigma[\pi^{-1}(\Phi)]) =\bigcup_{i=1}^{n_\Phi}\psi_i\subset [0,T]
$$
  consists of time intervals $\psi_i=(t_i,\widetilde{t}_i)$ such that $0\le t_1<\widetilde{t}_1\le t_2<\widetilde{t}_2\le \dots \le t_{n_\Phi }<\widetilde{t}_{n_\Phi}\le T$; all intervals have the same length $\widetilde{t}_i-t_{i}=\epsilon_\Phi$. On each family the set of functions $\tau^i$ is defined by
\begin{equation}
\label{eq2.24}
\tau^i_\Phi(x):=\psi_i\cap \rho(\mathscr{L}_\Sigma[\pi^{-1}(x)]),\qquad x\in\Phi,\quad i=1,\dots,
n_\Phi.
\end{equation}
  These functions are constant on determination sets: $\tau^i_\Phi(x)=\tau^i_\Phi(x')$ for regular $x,x'\in \Lambda_\Sigma[x]$. For them \eqref{eq2.19} is satisfied and the representation \eqref{eq2.20} is valid. The latter allows us to extend the functions $\tau^i_\Phi$ to the critical points (the ends of $\omega$ cells) by continuity.

  Let $\gamma\in\Sigma$, $\Phi\subset\Pi_\Sigma$ and $x\in \Phi$. It is not difficult to show that the set $\Lambda_\Sigma[x]\cap {\Omega^T_\gamma}$ is $\Lambda_\gamma$--closed (in $\Omega^T_\gamma$). Consequently, it has $\beta$--set of vectors $\beta^1_{\gamma\Phi},\dots, \beta^{n_\Phi}_{\gamma\Phi}$. We define them on the whole $\Lambda_\Sigma[x]$, extending them from $\Lambda_\Sigma[x]\cap \Omega^T_\gamma$ to $\Lambda_\Sigma[x]$ by zero.

Repeating the construction for all $\gamma\in\Sigma$, we obtain the collection of $\beta$--sets
\begin{equation}
\label{eq2.25}
\{\beta^1_{\gamma\Phi},\dots,\beta^{n_{\Phi}}_{\gamma\Phi}\mid\gamma\in\Sigma\},
\end{equation}
  where $\beta^i_{\gamma\Phi}$ is the vector with components $(\beta^i_{\gamma\Phi})_1,\dots,(\beta^i_{\gamma\Phi})_{m_{\Phi}}$ and $m_{\Phi}=\sharp\Lambda_\Sigma[x]$. Each of the sets is orthonormalized in $\mathbf{l}_2(\Lambda_\Sigma[x])$. Accordingly, on each family $\Phi\subset\Pi_\Sigma$ there are functions $\beta^i_{\gamma\Phi}(\,{\cdot}\,)$, which take constant values $(\beta^i_{\gamma\Phi})_k$ on cells $\omega_k\subset\Phi$.

  For a given family $\Phi$ and different vertices $\gamma\in\Sigma$, the set $\{\beta^1_{\gamma\Phi},\dots,\beta^{n_{\Phi}}_{\gamma\Phi}\}$ contains the same number of vectors equal to $n_{\Phi}$, and the functions $\tau^i_\Phi$ are the same. Nevertheless, let us assume by definition
\begin{equation}
\label{eq2.26}
n_{\gamma\Phi}:=n_{\Phi},\qquad \tau_{\gamma\Phi}^i(x):=\tau_{\\Phi}^i(x),\qquad x\in\Phi, \quad \gamma\in\Sigma.
\end{equation}
This seemingly redundant notation (index $\gamma$) will prove convenient in further consideration.

\section{Eikonals}
\label{s3}

\subsection{Reachable sets and projectors}
\label{ss3.1}

  Here we consider the problem \eqref{eq2.2}--\eqref{eq2.5} as a dynamical system and equip it with control theory attributes --- spaces and operators.

The space of controls $\mathscr{F}^T :=L_2(\Gamma \times [0,T])$ with the inner product
$$
(f,g)_{\mathscr{F}^T}=\sum_{\gamma \in \Gamma} \int_{0}^T f(\gamma ,t)g(\gamma ,t)\,dt
$$
  is called an \textit{external} space of the system \eqref{eq2.2}--\eqref{eq2.5}. It contains subspaces of the controls acting from separate boundary vertices $\gamma \in \Gamma$:
$$
\mathscr{F}^T_\gamma:=\bigl\{f \in \mathscr{F}^T\bigm| \operatorname{supp}f \subset \{\gamma\} \times [0,T]\bigr\}
$$
  Each control $f \in \mathscr{F}^T_\gamma$ has the form $f=\delta_{\gamma}\varphi$ with some $\varphi\in L_2[0,T]$. The subset of boundary vertices of $\Sigma\subseteq\Gamma$ corresponds to the subspace
$$
\mathscr{F}^T_\Sigma:={\sum_{\gamma \in \Sigma}}^\oplus \mathscr{F}^T_\gamma
$$
(the summands are orthogonal in $\mathscr{F}^T$).


  The space $\mathscr{H}=L_2(\Omega)$ is called \textit{internal}; the waves $u^f(\,\cdot\,,t)$ are its time-depen\-dent elements. For the set $B\subset\Omega$ we define a subspace $\mathscr H\langle B\rangle:=\{y\in\mathscr H\mid \operatorname{supp}y\subset\overline B\}$ of functions localized in $B$.

The set of waves
$$
\mathscr{U}^s_\gamma:=\{u^f(\,{\cdot}\,,s)\mid f \in \mathscr{ F}^T_\gamma\}\subset \mathscr{H}, \qquad 0\le s\le T,
$$
  is called a \textit{reachable} (from the vertex $\gamma$ to  moment $t=s$). From the representations \eqref{eq2.9} and \eqref{eq2.10} we see that $\mathscr{U}^s_\gamma$ are (closed) subspaces in $\mathscr H$. With growth of $s$ they expand: $\mathscr{U}^s_\gamma\subset\mathscr{U}^{s'}_\gamma$ at $s<s'$.

The locality of the ``control-wave'' correspondence noted in \eqref{eq2.21} leads to a decomposition on the families
\begin{equation}
\label{eq3.1}
\mathscr{U}^T_\gamma={\sum_{\Phi\subset\Pi_\Sigma}}^\oplus \mathscr{U}^T_\gamma\langle\Phi\rangle,
\end{equation}
  where the subspace $\mathscr{U}^T_\gamma\langle\Phi\rangle\subset\mathscr H\langle\Phi\rangle$ consists of waves $u^f(\,{\cdot}\,,T)\in\mathscr U^T_\gamma$ localized in $\Phi\cap\Omega^T_\gamma$. Orthogonality of the summands is a consequence of the disjunction of families: $\Phi^j\cap\Phi^k=\varnothing$ for $j\ne k$.

  Let us fix the boundary vertex $\gamma\in\Sigma$; let $P^T_{\gamma}$ be a projector in $\mathscr{H}$ on the subspace $\mathscr{U}^T_\gamma$. Now we will discuss its properties and describe its action.

As a consequence of \eqref{eq3.1}, we have the representation
\begin{equation}
\label{eq3.2}
P^T_\gamma=\sum _{\Phi \subset \Pi_\Sigma}P^T_{\gamma}\langle\Phi\rangle,
\end{equation}
  where $P^T_{\gamma}\langle\Phi\rangle$ are projectors in $\mathscr{H}$ on $\mathscr{U}^T_\gamma\langle\Phi\rangle$. Thus, the projector $P^T_\gamma$ is reduced by the subspaces $\mathscr{U}^T_\gamma\langle\Phi\rangle$, and to characterize it one needs to describe how $P^T_{\gamma}\langle\Phi\rangle$ acts.

  As shown in \cite{7}, the projectors $P^T_{\gamma}\langle\Phi\rangle$ can be expressed through vectors \eqref{eq2.25} and their corresponding functions $\beta^i_{\gamma\Phi}(\,{\cdot}\,)$ as follows:
\begin{equation}
\label{eq3.3}
(P^T_{\gamma}\langle\Phi\rangle y)(x)= \begin{cases}
{\displaystyle\sum _{i=1}^{n_{\gamma\Phi}} \langle y|_{\Lambda_\gamma[x]}, \beta^i_{\gamma\Phi}\rangle \beta^i_{\gamma\Phi}(x)}, &x \in \Phi,
\\
0, &x \in \Omega \setminus\Phi,
\end{cases}
\end{equation}
  where $y \in \mathscr H$ is an arbitrary function on a graph\footnote{In the exact sense \eqref{eq3.3} is the representation ``almost everywhere'' in $\Omega$. For $y\in C(\Omega)$ it is true ``everywhere''.}. This representation is derived from the expression \eqref{eq2.23}, for which the vectors $\beta^i$ were introduced.

  As can be seen from \eqref{eq3.2} and \eqref{eq3.3}, if the function $y\in\mathscr{U}^T_\gamma$ projected onto $\mathscr{H}$ is continuous and such that $y|_{\Lambda_\gamma[x]}\equiv 0$, then so is $(P^T_\gamma y)|_{\Lambda_\gamma[x]}\equiv 0$. In other words, the values of the projection $P^T_\gamma y$ on the set $\Lambda_\gamma[x]$ are quite determined by the values of $y$ on $\Lambda_\gamma[x]$. This is the reason for calling $\Lambda_\gamma[x]$, and with  them also $\Lambda_\Sigma[x]$, determination sets.

  Recall that the vectors $\beta^i_{\gamma\Phi}$ are elements of the subspace $\mathbb A[\Lambda_\Sigma[x]]\subset\mathbf{l}_2(\Lambda_\Sigma[x])$. It follows from the above that the projector $P^T_\gamma$ determines in $\mathbf{l}_2(\Lambda_\Sigma[x])$ the operator $p_{\gamma\Phi}[x]$ that projects onto $\mathbb A[\Lambda_\gamma[x]]\subset\mathbb A[\Lambda_\Sigma[x]]$. Let $\chi_1,\dots, \chi_{m_\Phi}$ be the (standard) basis in $\mathbf{l}_2(\Lambda_\Sigma[x])$ consisting of indicators of points of the set $\Lambda_\Sigma[x]$; in such a basis $(\beta^i_{\gamma\Phi})_k=\langle\beta^i_{\gamma\Phi},\chi_k\rangle$ holds. According to \eqref{eq3.3} in this basis, the projector matrix $p_{\gamma\Phi}[x]$ takes the form
\begin{equation}
\label{eq3.4}
{\check p}_{\gamma\Phi}[x]=B^*_{\gamma\Phi}[x]B_{\gamma\Phi}[x],\qquad B_{\gamma\Phi}[x]:=
\begin{pmatrix}
(\beta^1_{\gamma\Phi})_1(x) & \dots &(\beta^1_{\gamma\Phi})_{m_{\Phi}}(x) \\
(\beta^{2}_{\gamma\Phi})_1(x) & \dots & (\beta^{2}_{\gamma\Phi})_{m_{\Phi}}(x) \\
\dots & \dots & \dots \\
(\beta^{n_{\gamma\Phi}}_{\gamma\Phi})_1(x) & \dots & (\beta^{n_{\gamma\Phi}}_{\gamma\Phi}(x))
\end{pmatrix} .
\end{equation}

\subsection{Eikonal}
\label{ss3.2}

  The family $\{P^s_\gamma\mid 0\le s\le T\}$ of projectors in  $\mathscr{H}$ on reachable sets $\mathscr{U}^s_\gamma$ defines the eikonal operator (briefly -- \textit{eikonal})
$$
E^T_\gamma\colon\mathscr{H}\to \mathscr{H},\qquad E^T_\gamma:=\int^T_0 s\,dP^s_\gamma.
$$
  It follows from the definition that $E^T_\gamma$ is a bounded self-adjoint positive operator. Like the projector $P^T_\gamma$, the eikonal $E^T_\gamma$ is reduced by subspaces $\mathscr{H}\langle\Phi\rangle$: for this operator the relation $E^T_\gamma\mathscr{H}\langle\Phi\rangle\subset\mathscr{H}\langle\Phi\rangle$ and the expansion
$$
E^T_\gamma=\sum_{\Phi\subset\Pi}E^T_\gamma\langle\Phi\rangle
$$
  hold, where $E^T_\gamma\langle\Phi\rangle:=E^T_\gamma|_{\Phi}$ is the part of $E^T_\gamma$ acting in $\mathscr{H}\langle\Phi\rangle$.

As shown in \cite{7}, a representation
\begin{equation}
\label{eq3.5} (E^T_\gamma\langle\Phi\rangle y)(x)= \begin{cases}
{\displaystyle\sum _{i=1}^{n_{\gamma\Phi}}
\tau^i_{\gamma\Phi}(x)\langle
y|_{\Lambda_\gamma[x]},\beta^i_{\gamma\Phi}\rangle
\beta^i_{\gamma\Phi}(x)}, &x \in \Phi,
\\
0, &x \in \Omega \setminus\Phi
\end{cases}
\end{equation}
  consistent with \eqref{eq3.3} is valid, in which $y \in \mathscr H$ is arbitrary, and the functions $\tau^i_{\gamma\Phi}$ are given by definitions\footnote{The corresponding representation in \cite{7} makes the use of incorrectly defined functions $\tau^i_\Phi$ (see note 4). The \eqref{eq3.5} representation with functions defined in \eqref{eq2.18}, \eqref{eq2.26}, eliminates this mistake.} \eqref{eq2.24} and  \eqref{eq2.26}.

  The operator $E^T_\gamma$ induces in $\mathbf{l}_2(\Lambda_\Sigma[x])$ the operator $e_{\gamma\Phi}[x]$. According to \eqref{eq3.4} and \eqref{eq3.5}, in basis $\chi_1,\dots, \chi_{m_\Phi}$ its matrix is
\begin{gather*}
{\check e}_{\gamma\Phi}[x]=B^*_{\gamma\Phi}[x]D_{\gamma\Phi}[x] B_{\gamma\Phi}[x],
\\
B_{\gamma\Phi}[x]:= \begin{pmatrix}
(\beta^1_{\gamma\Phi})_1(x) & \dots &(\beta^1_{\gamma\Phi})_{m_{\Phi}}(x) \\
(\beta^{2}_{\gamma\Phi})_1(x) & \dots & (\beta^{2}_{\gamma\Phi})_{m_{\Phi}}(x) \\
\dots & \dots & \dots \\
(\beta^{n_{\gamma\Phi}}_{\gamma\Phi})_1(x) & \dots & (\beta^{n_{\gamma\Phi}}_{\gamma\Phi})_{m_{\Phi}}(x)
\end{pmatrix},\qquad D_{\gamma\Phi}[x]=\operatorname{diag}\{\tau^i_{\gamma\Phi}(x)\}_{i=1}^{n_{\gamma\Phi}}.
\end{gather*}
  The subspace $\mathbb A_{\gamma\Phi}[x]=\operatorname{span}\{\beta^1_{\gamma\Phi}, \dots, \beta^{n_{\gamma\Phi}}_{\gamma\Phi}\}\subset \mathbf{l}_2(\Lambda_\Sigma[x])$ reduces this matrix and its non-zero block in basis $\beta^1_{\gamma\Phi}, \dots, \beta^{n_{\gamma\Phi}}_{\gamma\Phi}$ is $\operatorname{diag}\{\tau^i_{\gamma\Phi}(x)\}_{i=1}^{n_{\gamma\Phi}}$. As can be seen from \eqref{eq2.19}, all its eigenvalues $\tau^i_{\gamma\Phi}(x)$ within cells of the family are different; when $x$ changes in $\omega\subset\Phi$ they sweep the intervals $(t_i,\widetilde t_i)\subset\Xi[\Phi]$. According to the following proposition \ref{p1}, at $T<T_\gamma$ the union (over all families $\Phi\subset\Pi$) of all segments $[t_i,\widetilde t_i]$ coincides with $[0,T]$.

  From the above considerations it is easy to provide some general properties of the eikonal as an operator in $\mathscr{H}$ (see, e.g., \cite{23}).

\begin{proposition}
\label{p1}   
Let $\overline{\operatorname{Ran}E^T_\gamma}=\mathscr{U}^T_\gamma$ and $\operatorname{Ker}E^T_\gamma=\mathscr{H}\ominus\mathscr{U}^T_\gamma$ hold for the eikonal. The eikonal is reduced by the parts of the reachable set: {\vrule width0pt height10pt$E^T_\gamma\mathscr{U}^T_\gamma\langle\Phi\rangle \subset\mathscr{U}^T_\gamma \langle\Phi\rangle$}, $\Phi\subset\Pi_\Sigma$. If $T<T_\gamma$, the operator $E^T_\gamma|_{\mathscr{U}^T_\gamma}$ has a simple absolutely continuous spectrum filling the segment {\vrule width0pt height9pt$[0,T]$}.
\end{proposition}

\begin{remark}
\label{r1}
From the representation \eqref{eq3.5} it follows that at times $T>T_\gamma$ the spectrum $E^T_\gamma|_{\mathscr{U}^T_\gamma}$ is the union of segments $[0,T_0]\cup[T_1,T_2]\cup\dots \cup[T_{N-1},T_N]$, where $T_\gamma\le T_0<T_1<\dots <T_N\le T$, and each segment consists of ranges of functions\footnote{But the exact description of the spectrum is an open question. There is a hypothesis that it is always exhausted by the interval $[0,T_0]$ with  sufficiently large $T_0$. The question is related to the subtle details of the structure of the reachable sets $\mathscr{U}^T_\gamma$ of the metric graph.} $\tau^i_{\gamma\Phi}$ (closures of time cells ${\psi^i_{\gamma\Phi}}$, see \eqref{eq2.17}).
\end{remark}

\subsection{Parameterization}
\label{ss3.3}

Let us choose a family $\Phi=\bigcup_{k=1}^{m_{\Phi}}\omega_k\subset\Pi_\Sigma$; let $\omega=(c,c')\subset \Phi$ be one of its cells. Recall that all the cells are of the same length $\epsilon_\Phi=\tau(c,c')$. For $x\in\omega$ let us introduce the notation $x=x(r)$ if $\tau(c,x)=r$. Along with $x$, the determination set also turns out to be parameterized: $\Lambda_\Sigma[x(r)]=\{x_k(r)\}_{k=1}^{m_{\Phi}}$. As $r$ changes in the interval $(0,\epsilon_\Phi)$, the points of $x_k(r)$ continuously change position and cover the cells $\omega_k$. Thus, the family $\Phi$ is parametrized.

All elements of the representations \eqref{eq3.3} and \eqref{eq3.5} are also parameterized by $r$: the vectors
\begin{equation}
\begin{gathered}
\beta^i_{\gamma\Phi}=\{(\beta^i_{\gamma\Phi})_k(r)\}_{k=1}^{m_{\Phi}},
\\
(\beta^i_{\gamma\Phi})_k(r):= \beta^i_{\gamma\Phi}(x_k(r)) =(\beta^i_{\gamma\Phi})_k=\mathrm{const},\qquad 0<r<\epsilon_\Phi,
\label{eq3.6}
\end{gathered}
\end{equation}
and the functions $\tau^i_{\gamma\Phi}(r):=\tau^i_{\gamma\Phi}(x(r))$. The latter, according to \eqref{eq2.20} take the values \begin{equation} \label{eq3.7} \tau^i_{\gamma\Phi}(r)=t_{i\Phi}+r \quad\text{or}\quad \tau^i_{\gamma\Phi}(r)=\widetilde{t}_{i\Phi}-r=(t_{i\Phi}+\epsilon_\Phi)-r. \end{equation} Note that there are two parametrizations of the $\Phi$ family: the c $r=\tau(x,c)$ adopted above and the one corresponding to the parameter $r=\tau(x,c')$. They are quite equivalent. In what follows, we will assume that each family $\Phi\subset\Pi_\Sigma$ is parameterized in one of two ways.

Parametrization determines matrix representations of functions and operators on the graph.

Let $\Phi\subset\Pi_\Sigma$ be a parameterized family, $y\in\mathscr H$ a function on the graph, $x=x(r)\in\Lambda_\Sigma[x(r)]=\{x_k(r)\}_{k=1}^{m_{\Phi}}\subset\Phi$, $0<r<\epsilon_{\Phi}$. It is easy to check that the map
$$
U_\Phi\colon \mathscr H\to L_2([0, \epsilon_{\Phi};
\mathbb{R}^{m_{\Phi}}),\qquad (U_\Phi y)(r):= \begin{pmatrix}
y(x_1(r))\\ \dots \\y(x_{m_{\Phi}}(r))
\end{pmatrix}, \quad r\in(0,\epsilon_{\Phi}),
$$
is unitary. For each vertex $\gamma\in\Sigma$ we define the constant (see \ \eqref{eq3.6}) columns and the matrices formed by them
$$
\beta^i_{\gamma\Phi}= \begin{pmatrix}
(\beta^i_{\gamma\Phi})_1\\ \dots \\(\beta^i_{\gamma\Phi})_{m_{\Phi}}
\end{pmatrix} \in \mathbb{R}^{m_{\Phi}},\qquad B_{\gamma\Phi}:= \begin{pmatrix}
(\beta^1_{\gamma\Phi})_1 & \dots & (\beta^1_{\gamma\Phi})_{m_{\Phi}} \\
(\beta^{2}_{\gamma\Phi})_1 & \dots & (\beta^{2}_{\gamma\Phi})_{m_{\Phi}} \\
\dots & \dots & \dots \\
(\beta^{n_{\gamma\Phi}}_{\gamma\Phi})_1 & \dots & (\beta^{n_{\gamma\Phi}}_{\gamma\Phi})_{m_{\Phi}}
\end{pmatrix};
$$
and introduce the matrices
$$
D_{\gamma\Phi}(r):=\{\tau^i_{\gamma\Phi}(r)\,\delta_{ij}\}_{i,j=1}^{n_{\gamma\Phi}},\qquad
r\in(0,\epsilon_{\Phi}),
$$
where $\tau^i_{\gamma\Phi}(r)$ is of the form \eqref{eq3.7}. The matrix $B_{\gamma\Phi}$ does not change when the parameter $r$ is varied because its elements in the cells $\omega_1,\dots,\omega_{m_\Phi}$ of the family $\Phi$ are constant. The matrix $B_{\gamma\Phi}^*B_{\gamma\Phi}$ is also constant. The latter, since the columns of $\beta^i_{\gamma\Phi}$ constitute an orthonormalized set, is a projector in ${\mathbb R}^{m_{\Phi}}$ on the subspace
$$
\mathscr{A}_{\gamma}[\Phi]:=\operatorname{span}\{\beta^1_{\gamma\Phi},
\dots, \beta^{n_{\gamma\Phi}}_{\gamma\Phi}\}=[B^*_{\gamma\Phi}
B_{\gamma\Phi}] \mathbb{R}^{m_{\Phi}}.
$$
The matrix-projector $B_{\gamma\Phi}^*B_{\gamma\Phi}$ can be decomposed into a sum of pairwise orthogonal one-dimensional projectors
\begin{equation}
\label{eq3.8}
B_{\gamma\Phi}^*B_{\gamma\Phi}=\sum_{i=1}^{n_{\gamma\Phi}}P^i_{\gamma\Phi},\qquad
P^i_{\gamma\Phi}:=\langle\,{\cdot}\,,\beta^i_{\gamma\Phi}\rangle\,\beta^i_{\gamma\Phi},
\end{equation}
where $\langle\,{\cdot}\,,{\cdot}\,\rangle$ is the standard inner product in $\mathbb R^{m_\Phi}$. According to the \eqref{eq3.5}, we have the representation
$$
(U_\Phi\, E^{\,T}_\gamma\!\langle\Phi\rangle y)(r)=[B_{\gamma\Phi}^* D_{\gamma\Phi}(r) B_{\gamma\Phi}](U_\Phi{y})(r), \qquad r \in (0,\epsilon_{\Phi})
$$
with matrices
\begin{equation}\label{eq3.9}
B_{\gamma\Phi}^* D_{\gamma\Phi}(r)B_{\gamma\Phi}=U_\Phi
E^T_\gamma\langle\Phi\rangle U^{-1}_\Phi\stackrel{\eqref{eq3.5},
\eqref{eq3.8}}{=}
\sum_{i=1}^{n_{\gamma\Phi}}\tau_{\gamma\Phi}^i(r)
P_{\gamma\Phi}^i.
\end{equation}

Let us describe the parametrization of spaces and operators corresponding to the decomposition $\Pi_{\Sigma}$ in whole. The analogue of the decomposition \eqref{eq3.1} takes the form
$$
\mathscr{U}^T_\Sigma = \oplus \sum_{\Phi \subset \Pi_\Sigma}
\mathscr{U}^T_\Sigma\langle\Phi\rangle,\qquad
\mathscr{U}^T_\Sigma\langle\Phi\rangle:=\operatorname{span}
\{\mathscr{U}^T_\gamma\langle\Phi\rangle\mid \gamma\in\Sigma\},
$$
and each of the subspace summands reduces all eikonals simultaneously:
$$
E^T_\gamma\mathscr{U}^T_\Sigma\langle\Phi\rangle\subset\mathscr{U}^T_\Sigma\langle\Phi\rangle, \qquad \gamma\in\Sigma.
$$
Using parametrizations in families, we have
$$
U_\Phi \mathscr{H}\langle \Phi \rangle =
L_2([0,\epsilon_\Phi];\mathbb{R}^{m_\Phi}),\qquad U_\Phi
\mathscr{U}^T_\Sigma[\Phi] =
L_2([0,\epsilon_\Phi];\mathscr{A}_\Sigma[\Phi]),
$$
where $\mathscr{A}_\Sigma[\Phi]: =\operatorname{span}\{\mathscr{A}_\gamma[\Phi]\mid \gamma\in\Sigma\}$, and the parts $E^T_\gamma\langle{\Phi}\rangle$ of the eikonals multiply the elements $L_2([0,\epsilon_\Phi];{\mathbb R}^{m_\Phi})$ by matrix-function \eqref{eq3.9}.

\begin{convention}
\label{con1} 
Let us specify the notations. For spaces $\mathscr S_1,\dots,\mathscr S_n$ the sum $\mathscr{S}=\bigoplus\sum_j \mathscr{S}_j$ is the space of sets $s=\{s_1,\dots,s_n\}$, $s_j\in\mathscr S_j$ (with component operations). For operators $A_1,\dots, A_n$, $A_j\in\operatorname{End}\mathscr{S}_j$, the sum $A=\bigoplus\sum_j A_j\in \operatorname{End}\mathscr{S}$ is the operator acting on the rule $As:=\{A_1s_1,\dots,A_ns_n\}$. For matrices $M_1,\dots,M_n$, $M_j\in \mathbb M^{\varkappa_j}$, the sum of $M=\bigoplus\sum_j M_j\in\mathbb M^{\varkappa_1+\dots+\varkappa_n}$ is a blockwise=diagonal matrix with blocks $M_1,\dots, M_n$ (we also write $[M]_j=M_j$). For algebras $\{\mathfrak{A}_1,\dots, \mathfrak{A}_n\}$\enskip $\mathfrak{A}=\bigoplus_j\mathfrak{A}_j$ is a direct sum of algebras--subjects (we also write $[\mathfrak{A}]_j=\mathfrak{A}_j$ and call $\mathfrak{A}_j$ blocks).
\end{convention}

Parameterization of the whole $\Pi_\Sigma$ is executed  by the operator $U:=\bigoplus\sum\limits_{\Phi\subset\Pi_\Sigma}U_\Phi$:
\begin{align}
\notag U\mathscr{H}\langle \Omega^T_\Sigma \rangle
&=\bigoplus\sum_{\Phi\subset\Pi_\Sigma}L_2([0,\epsilon_\Phi];\mathbb{R}^{m_\Phi}),
\\
\label{eq3.10} UE^T_\gamma U^{-1}
&=\bigoplus\sum_{\Phi\subset\Pi_\Sigma}U_\Phi E^T_\gamma
\langle\Phi\rangle U^{-1}_\Phi \stackrel{\eqref{eq3.9}}{=}
\bigoplus\sum_{\Phi\subset\Pi_\Sigma}\sum_{i=1}^{n_{\gamma\Phi}}\tau^i_{\gamma\Phi}P^i_{\gamma
\Phi}, \qquad \gamma\in\Sigma.
\end{align}

\subsection{Shifted eikonals}
\label{ss3.4}

For technical reasons, while describing an algebra generated by
eikonals, it is more convenient to use operators (\textit{shifted
eikonals})
\begin{equation}
\label{eq3.11}
\dot E^T_{\gamma}:=\int_0^T(s+1)\,dP^s_{\gamma} = E^T_{\gamma}+ P^T_{\gamma}
\end{equation}
The previously established properties and representations for
$E^T_\gamma$ are obviously reformulated for $\dot E^T_\gamma$.
Thus, the analog of the representation \eqref{eq3.10} takes the
form
\begin{align}
\notag U\dot E^T_{\gamma}U^{-1}
&=\bigoplus\sum_{\Phi\subset\Pi_\Sigma}U_\Phi \dot
E^T_{\gamma}\langle\Phi\rangle U^{-1}_\Phi
=\bigoplus\sum_{\Phi\subset\Pi_\Sigma}[B_{\gamma\Phi}^* \dot
D_{\gamma\Phi}(\, {\cdot}\,) B_{\gamma\Phi}]
\\
&\!\!\!\!\!\stackrel{\eqref{eq3.10}}{=}
\bigoplus\sum_{\Phi\subset\Pi_\Sigma}\sum_{i=1}^{n_{\gamma\Phi}}\dot\tau^i_{\gamma
\Phi}P^i_{\gamma \Phi}, \qquad \gamma\in\Sigma, \label{eq3.12}
\end{align}
where $\dot
D_{\gamma\Phi}(\,{\cdot}\,):=D_{\gamma\Phi}(\,{\cdot}\,)+ I$,
$I$ is a unit matrix of appropriate dimension and
$\dot\tau_{\gamma\Phi}^i(r) :=\tau_{\gamma\Phi}^{i}(r)+1$. The
analogue of the proposition \ref{p1} is as follows.
\begin{proposition}
\label{p2} For the operator $\dot E^T_{\gamma}$,
${\operatorname{Ran}\dot E^T_{\gamma}}=\mathscr{U}^T_\gamma$,
$\operatorname{Ker}\dot
E^T_{\gamma}=\mathscr{H}\ominus\mathscr{U}^T_\gamma$ and $\dot
E^T_{\gamma}\mathscr{U}^T_\gamma\langle\Phi\rangle\subset\mathscr{U}^T_\gamma
\langle\Phi\rangle$, $\Phi\subset\Pi$. When $T<T_\gamma$ it has
eigenvalue $0$ of infinite multiplicity and a simple absolutely
continuous spectrum filling the segment $[1,T+1]$. When
$T>T_\gamma$, according to the observation \ref{r1}
\begin{equation}
\sigma(\dot
E^T_\gamma|_{\mathscr{U}^T_\gamma})\,{=}\,\sigma_{\mathrm{ac}}(\dot
E^T_\gamma) \,{=}\,[1,T_0+1]\cup[T_1+1,T_2+1]\,{\cup}\,{\cdots}\,
{\cup}\,[T_{N-1}+1,T_N+1], \label{eq3.13}
\end{equation}
where $T_\gamma\le T_0<T_1<\dots <T_N\le T$, and the segments are
the unions of the regions of values of functions
$\dot\tau^i_{\gamma\Phi}$ (shifted by  $1$ time cells
$\dot\psi^i_{\gamma\Phi}$).
\end{proposition}

\section{Algebra of eikonals}
\label{s4}

\subsection{Definitions and general facts}
\label{ss4.1}

Recall that a $C^*$-algebra {$\mathfrak{A}$ is a Banach algebra
with an involution  $a \to a^*$, which satisfies the identity
\cite{24}, \cite{25}:
$$
\|a^*a\| = \|a\|^2,\quad a\in\mathfrak{A}.
$$}
In particular, such are the algebras of bounded operators
$\mathfrak{B}(\mathscr{H})$ in  Hilbert space $\mathscr{H}$, in
which the operator conjugation plays the role of involution.
Writing ${\mathfrak A}\cong{\mathfrak B}$ would mean that
$C^*$-algebras ${\mathfrak A}$ and ${\mathfrak B}$ are related by
isometric $*$-isomorphism ({then briefly -- \textit{isomorphic}}).
For the set $S\subset \mathfrak{A}$, the minimal
$C^*$-(sub)algebra in $\mathfrak{A}$ containing $S$ is denoted by
$\vee S$.

By $\mathbb{M}^n$ is meant the algebra of real $(n\times
n)$--matrices treated as operators in $\mathbb{M}^n$ and endowed
with an appropriate (operator) norm. It is \textit{irreducible}.

Through $C([a,b],\mathbb{M}^n)$ denotes the algebra of continuous
$\mathbb{M}^n$-valued functions with norm $\|c\|=\sup_{a\le t\le
b} \|c(t)\|_{\mathbb{M}^n}$. By the same symbol, we denote the
operator (sub)algebra in $\mathfrak{B}(L_2([a,b];\mathbb{M}^n))$,
whose elements multiply square--summable $\mathbb R^n$-valued
functions by functions from $C([a,b],\mathbb{M}^n)$. The
correspondence $c\mapsto c\cdot$ establishes the isomorphism
between these algebras.

The $C^*$-subalgebra $\mathfrak A\subset \mathbb{M}^n$ is
considered irreducible if $\mathfrak A\cong \mathbb{M}^k$, where
$k\le n$ is satisfied. Such an algebra in  a suitable basis
in $\mathbb R^n$ takes a block=diagonal form and consists of two
blocks, one of which is $\mathbb{M}^k$ and the second (if
available) -- zero.

Here is a summary of known results\footnote{In  article \cite{8}
in analogous summary ``On matrix algebras'', there is a mistake:
statement 3 is incorrect. However, after appropriate corrections,
all results of the work remain valid.}

\begin{proposition}
\label{p3} Any $C^*$-subdalgebra of the algebra $\mathbb{M}^n$ is
isometric to the direct sum $\bigoplus_{k}\mathbb{M}^{n_k}$, where
$\sum_{k}{n_k}\le n$.
\end{proposition}

\begin{proposition}[(see\ \cite{24})]
\label{p4} Let $\mathfrak{P}\subset\mathbb{M}^n$ and
$C^*$-subdalgebra $\mathfrak{A}\subset C([a,b]; \mathfrak{P})$ is
such that for any $t,t'\in[a,b]$ and $p,p'\in\mathfrak{P}$ there
exists an element $u\in\mathfrak{A}$, for which $u(t)=p$,
$u(t')=p'$ is satisfied. Then $\mathfrak{A}=
C([a,b];\mathfrak{P})$.
\end{proposition}

The representation of the $C^*$-algebra $\mathfrak{A}$  is a
homomorphism $\pi\colon \mathfrak{A}\to\mathfrak{B}(H)$, where
$H$ is a Hilbert space. Equivalence of representations
$\pi\sim\pi'$ means that $\iota\,\pi(a)=\pi'(a)\iota$,
$a\in\mathfrak{A}$, where $\iota\colon H\to H'$ is an isometry of
representation spaces. A representation is \textit{irreducible} if
the operators $\pi(\mathfrak{A})$ do not have a common nonzero
invariant subspace in $H$.

The spectrum of $C^*$-algebra $\mathfrak{A}$ is the set
$\widehat{\mathfrak{A}}$ of equivalence classes of its irreducible
representations. The equivalence class (spectrum point)
corresponding to the representation $\pi$ will be denoted by
$\widehat\pi$. The spectrum is endowed with the canonical Jacobson
topology \cite{24}, \cite{25}.

Isomorphism of the algebras $\mathrm{u}\colon
\mathfrak{A}\to\mathfrak{B}$ determines the correspondence of
representations
\begin{align}
\label{eq4.1}
\widehat{\mathfrak{A}}\ni\pi\to\mathrm{u}_*\pi\in\widehat{\mathfrak{B}},
\quad(\mathrm{u}_*\pi)(b):=\pi(\mathrm{u}^{-1}(b)), \qquad
b\in\mathfrak{B},
\end{align}
which continues up to the canonical homeomorphism of the spectra:
\begin{align}
\label{eq4.2}
\widehat{\mathfrak{A}}\ni\widehat\pi\to\mathrm{u}_*\widehat\pi\in\widehat{\mathfrak{B}},\qquad
\mathrm{u}_*\widehat\pi:=\{\mathrm{u}_*\pi\mid \pi\in\widehat\pi\}.
\end{align}

\begin{proposition}
\label{p5}
Representations
\begin{equation}
\label{eq4.3}
\pi_t\colon C([a,b],\mathbb{M}^n)\to\mathbb{M}^n,\qquad \pi_t(\phi):=\phi(t),
\end{equation}
{\bad are irreducible; their equivalence classes exhaust the
spectrum of the algebra $C([a,b],\mathbb{M}^n)$}. For any
irreducible representation $\pi$ of the algebra
$C([a,b],\mathbb{M}^n)$, there exists a single point $t\in[a,b]$
such that $\pi\sim\pi_t$ holds.
\end{proposition}

Let
$$
\dot C([a,b]; \mathbb M^n):= \{\phi\in C([a,b]; \mathbb M^n)\mid
\phi(a)\in \mathbb M_a,\, \phi(b)\in \mathbb M_b\}
$$
where $\mathbb M_a,\mathbb M_b$ are $C^*$-subdalgebras $\mathbb
M^n$, which we call \textit{boundary}. It follows from the
proposition \ref{p3} that
\begin{equation}
\label{eq4.4} \mathbb M_a\cong\bigoplus_{k=1}^{n_a} \mathbb
M^{\varkappa_k},\quad \varkappa_1+\dots+\varkappa_{n_a}\le n;
\qquad \mathbb M_b\cong\bigoplus_{k=1}^{n_b}\mathbb
M^{\lambda_k},\quad \lambda_1+\dots+\lambda_{n_b}\le n
\end{equation}
holds. In the case $\mathbb M_a=\mathbb M_b=\mathbb M^n$ we have
$\dot C([a,b]; \mathbb M^n)= C([a,b]; \mathbb M^n)$. Algebras
$\dot C([a,b]; \mathbb M^n)$ we will call \textit{standard}. The
spectrum of standard algebra consists of the classes
$\widehat\pi_t$, $t\in(a,b)$ corresponding to representations
\eqref{eq4.3}, and irreducible representations, which constitute
(possibly, reducible) representations $\widehat\pi_a$,
$\widehat\pi_b$. If, for example, $n_a\ge 2$, then $\pi_a$
decomposes into irreducible representations
$$
\pi_a^k\colon \phi(a)\mapsto [\phi(a)]^k\in\mathbb M^{\varkappa_k},
$$
where $[\,{\cdots}\,]^k$ is the $k$th block matrix
in representations of \eqref{eq4.4}. In this case we say that
$\widehat\pi_a^1,\dots,\widehat\pi_a^{n_a}$ form a
\textit{cluster} in the spectrum of standard algebra. This term is
motivated by the fact that they are inseparable from each other in
Jacobson's topology \cite{27}. A similar cluster may exist at the
right end $t=b$. At the same time, all $\widehat\pi_t$ with
different $t\in(a,b)$ are separable from each other and from
clusters (see \ \cite{7}, \cite{8}). The spectrum of $C([a,b];
\mathbb M^n)$ contains no clusters.
\medskip

The central object of the paper is the \textit{eikonal algebra} of
the graph $\Omega$
$$
\mathfrak{E}^T_{\Sigma}:=\vee\{E^T_{\gamma}\mid
\gamma\in\Sigma\}\subset\mathfrak{B}(L_2(\Omega))
$$
(see \cite{7}, \cite{8}). The subalgebras
$$
\mathfrak{E}^T_{\gamma}:=\vee E^T_\gamma,\qquad \gamma\in\Sigma
$$
are said to be \textit{partial}. In considerations it is
convenient to use the ``shifted'' algebras
$$
\dot{\mathfrak{E}}^T_{\gamma}:=\vee \dot
E^T_\gamma,\qquad\dot{\mathfrak{E}}^T_{\Sigma}:=\vee\{\dot
E^T_{\gamma}\mid \gamma\in\Sigma\}=
\vee\{\dot{\mathfrak{E}}^T_{\gamma}\mid \gamma\in\Sigma\}.
$$
The transition from $\mathfrak{E}^T_{\gamma}$ to 
$\dot{\mathfrak{E}}^T_{\gamma}$ consists of adding the projector
$P^T_\gamma$, which plays the role of the unit
in $\dot{\mathfrak{E}}^T_{\gamma}$ (see \ \eqref{eq3.11}), while
the algebra $\mathfrak{E}^T_{\gamma}$ turns out to be a subalgebra
in $\dot{\mathfrak{E}}^T_{\gamma}$.

According to the functional calculus of self-adjoint operators and
due to the orthogonality of the projectors $P^i_{\gamma \Phi}$
in \eqref{eq3.12}, we have
\begin{equation}
\label{eq4.5} \varphi(\dot
E^T_\gamma)=\int_{\sigma_{\mathrm{ac}}(\dot
E^T_\gamma)}\varphi(s)\, dP^s_\gamma
=U^{-1}\biggl[\bigoplus\sum_{\Phi\subset\Pi_\Sigma}\sum_{i=1}^{n_{\gamma\Phi}}
(\varphi\circ\dot\tau^i_{\gamma\Phi})P^i_{\gamma \Phi}\biggr]U
\end{equation}
for $\varphi\in C(\sigma_{\mathrm{ac}}(\dot E^T_\gamma))$. The
correspondence $\varphi(\dot E^T_\gamma)\leftrightarrow\varphi$
given by the first equality is an isomorphism of the algebras
$\dot{\mathfrak{E}}^T_\gamma$ and $C(\sigma_{\mathrm{ac}}(\dot
E^T_\gamma))$.

\begin{convention}
\label{con2} Hereafter, unless otherwise specified, we deal only
with shifted eikonals and omit $(\,\dot{}\,)$ in  notation: $\dot
E^T_\gamma\equiv E^T_\gamma$, $\dot\tau^k_{\gamma
l}\equiv\tau^k_{\gamma l}$, $\dot\psi^k_{\gamma
l}\equiv\psi^k_{\gamma l}$,
$\dot{\mathfrak{E}}^T_{\Sigma}\equiv\mathfrak{E}^T_{\Sigma}$ and
so on.
\end{convention}

\subsection{Representations and connections between blocks}
\label{ss4.2}

From \eqref{eq3.2} and \eqref{eq3.5} we have the representation
$$
\mathfrak{E}^T_\Sigma=\vee\left\{\bigoplus\sum_{\Phi\subset\Pi_\Sigma}
E_\gamma\langle\Phi\rangle\,\,\bigg|\,\,\gamma\in\Sigma\right\}.
$$
Parameterization \eqref{eq3.10} leads to the eikonal algebra
representation
\begin{align}
\notag
\mathfrak{E}^T_\Sigma &\cong U\mathfrak{E}^T_\Sigma U^{-1}
\\
&\!\!\!\!\stackrel{\eqref{eq3.12}}{=}
\vee\biggl\{\bigoplus\sum_{\Phi\subset\Pi_\Sigma}\sum_{i=1}^{n_{\gamma\Phi}}
\tau^i_{\gamma\Phi}(\,{\cdot}\,) P^i_{\gamma\Phi}\biggm|
\gamma\in\Sigma\biggr\}\subset\bigoplus_{\Phi\subset\Pi_\Sigma}
C([0,\epsilon_\Phi],\mathbb{M}^{m_{\Phi}}) \label{eq4.6}
\end{align}
in the form of an operator algebra; its elements multiply
functions from the representation space
\begin{equation}
\label{eq4.7}
\mathscr{R}^T_\Sigma:=\bigoplus\sum_{\Phi\subset\Pi_\Sigma} L_2([0,\epsilon_\Phi];{\mathbb R}^{m_{\Phi}})
\end{equation}
by continuous matrix-functions of the corresponding form. In  more
obvious block\-matrix notation of the representation \eqref{eq4.6}
we have
\begin{align}
\notag
U\mathfrak{E}^T_\Sigma U^{-1} &= \vee\left\{\begin{pmatrix}
{\displaystyle\sum_{i=1}^{n_{\gamma \Phi^1}} \tau^i_{\gamma \Phi^1}(\cdot_1) P^i_{\gamma \Phi^1}}\\
&\ddots\\\\
&&{\displaystyle\sum_{i=1}^{n_{\gamma \Phi^J}} \tau^i_{\gamma\Phi^J}(\cdot_J) P^i_{\gamma \Phi^J}}
\end{pmatrix}\Biggm| \gamma\in\Sigma\right\}
\\
&\subset \begin{pmatrix}
C([0,\epsilon_1]; \mathbb{M}^{m_{\Phi^1}})\\
& \ddots\\\
&&C([0,\epsilon_J]; \mathbb{M}^{m_{\Phi^J}}
\end{pmatrix},
\label{eq4.8}
\end{align}
here the zero off-diagonal blocks are omitted,
$\Pi_\Sigma=\Phi^1\cup\dots \cup\Phi^J$.  The notation
$(\,\cdot_j)$ emphasizes that the arguments $r_j\in[0,\epsilon_j]$
of the functions $\tau^i_{\gamma \Phi^j}$ are different
in accordance with representation \eqref{eq4.7}.

Let us define the projector sets $\mathbb
P_{\Phi^j}:=\{P^i_{\gamma \Phi^j}\mid i=1,\dots,n_{\gamma
\Phi^j};\, \gamma\in\Sigma\}$. Then, concerning the algebras
\begin{equation}
\label{eq4.9} \mathfrak{P}_{\Phi^j}:=\vee\mathbb
P_{\Phi^j}\subseteq\mathbb{M}^{m_{\Phi^j}},
\end{equation}
the embedding in \eqref{eq4.8} is specified as follows:
\begin{equation}
\label{eq4.10}
U\mathfrak{E}^T_\Sigma{ U^{-1}}\subset\bigoplus_{j=1}^J C([0,\epsilon_j];\mathfrak{P}_{\Phi^j}).
\end{equation}
As can be foreseen from \eqref{eq4.8} and \eqref{eq4.10}, the description of the structure of the eikonal algebra is reduced to the specification of the connections between its blocks $[U\mathfrak{E}^T_\Sigma U^{-1}]_j(\,\cdot_j)$ corresponding to different families $\Phi^j\subset\Pi_\Sigma$.  It is these
connections that distinguish $U\mathfrak{E}^T_\Sigma U^{-1}$ from
the algebra in right side of \eqref{eq4.10}, which has quite
independent blocks. The following lemma is  a step in studying the
relations between the blocks of the algebra
$U\mathfrak{E}^T_\Sigma{ U^{-1}}$.

Let us introduce the projectors
$$
\mathcal P^i_{\gamma \Phi^j}{:=} \begin{pmatrix}
O_1 \\
& \ddots\\
&& P^i_{\gamma \Phi^j}\\
&& \ddots\\
&&&& O_{\mathcal J}
\end{pmatrix}
{\in} \begin{pmatrix}
\mathfrak{P}_{\Phi^1} \\
& \ddots\\\
&& \mathfrak{P}_{\Phi^j}\\
&& \ddots\\
&&&& \mathfrak{P}_{\Phi^J}
\end{pmatrix}
{=}\bigoplus_{j=1}^J \mathfrak{P}_{\Phi^j},
$$
where $O_k$ are the zero matrices of appropriate dimensions. Like
$P^i_{\gamma \Phi^j}$, these projectors are pairwise orthogonal.
Let us also define the ``points'' $\mathbf{r}:=\{r_1,\dots ,r_J\}$
with coordinates $r_j\in[0,\epsilon_j]$, the matrices
\begin{equation}
\label{eq4.11} (U{ E}^T_\gamma{ U^{-1}})(\mathbf{r}) :
=\bigoplus\sum_{j=1}^J\sum_{i=1}^{n_{\gamma\Phi^j}}\tau^i_{\gamma
\Phi^j}(r_j)
P^i_{\gamma\Phi^j}=\sum_{j=1}^J\sum_{i=1}^{n_{\gamma\Phi^j}}\tau^i_{\gamma\Phi^j}(r_j)
\mathcal P^i_{\gamma\Phi^j} \in\bigoplus_{j=1}^J
\mathfrak{P}_{\Phi^j}
\end{equation}
(see \ \eqref{eq3.12}) and the matrix algebras
$$
(U\mathfrak{E}^T_\Sigma U^{-1})(\mathbf{r}):=\vee \{(UE^T_\gamma
U^{-1})(\mathbf{r})\mid \gamma\in\Sigma\}.
$$

\begin{lemma}
\label{l2} Fix $\gamma$, $i$, $j$. Let the points $\mathbf{r}$ and
$\mathbf{r}'$ be such that their coordinates satisfy
$r_j\in(0,\epsilon_j)$ and ${r_j}\ne {r_j'}$. Then there exists an
element $e\in U\mathfrak{E}^T_\Sigma U^{-1}$, for which
$e(\mathbf{r})=\mathcal{P}^i_{\gamma j}$ and $e(\mathbf{r}')=O$
holds.
\end{lemma}

\begin{proof}
Due to the pairwise orthogonality of the projectors
in \eqref{eq4.11}, for $s\in\mathbb N$ we have
$$
\bigl((U E^T_\gamma U^{-1})^s\bigr)(\mathbf{r})
=\sum_{k=1}^J\sum_{l=1}^{n_{\gamma\Phi^k}}\bigl((\tau^l_{\gamma
\Phi^k}(r_k)\bigr)^s\, \mathcal P^l_{\gamma \Phi^k}
$$
As a consequence, for the polynomial $q=q(t)=a_\nu t^\nu+\dots+a_1t$, the
$$
\bigl(q(U E^T_\gamma U^{-1})\bigr)(\mathbf{r}) =\sum_{k=1}^J\sum_{l=1}^{n_{\gamma\Phi^k}}q\bigl(\tau^l_{\gamma \Phi^k}(r_k)\bigr) \mathcal P^l_{\gamma \Phi^k}
$$

By the condition on the coordinates of the points $\mathbf{r}$ and
$\mathbf{r}'$ and the property \eqref{eq2.19}, in the set of
numbers
$$
\{\tau^l_{\gamma\Phi^k}(\eta)\mid k=1,\dots,J;\,\,
l=1,\dots,n_{\Phi^J};\,\, \eta=r_k,r'_k\}
$$
the number $\tau^i_{\gamma j}(r_j)$ occurs once. Choose a
polynomial so that it equals $1$ at $t=\tau^i_{\gamma j}(r_j)$ and
$0$ at all other points in the set. For $e:=q(UE^T_\gamma
U^{-1})$, we obviously have $e(\mathbf{r})= \mathcal P^i_{\gamma
j}$, $e(\mathbf{r}')=O$. The lemma is proved.
\end{proof}

\begin{corollary}
\label{c2} If the coordinates of the point $\mathbf{r}$ do not
take extreme values, i.e. $r_j\notin\{0,\epsilon_j\}$ for all $j$,
then the
    \begin{equation}
        \label{eq4.12}
        (U\mathfrak{E}^T_\Sigma U^{-1})(\mathbf{r}) =\bigoplus_{j=1}^J \mathfrak{P}_{\Phi^j}
    \end{equation}
If the coordinates are such that $0<a_j\le r_j\le b_j<\epsilon_j$
is satisfied for a fixed $j$, and the other coordinates are
arbitrary, then the following relation holds
    \begin{equation}
        \label{eq4.13}
        {{[U\mathfrak{E}^T_\Sigma U^{-1}]}_j}\big|_{a_j\le r_j\le b_j}=C([a_j,b_j]; \mathfrak{P}_{\Phi^j}).
    \end{equation}
\end{corollary}

The first relation is a straightforward consequence of the lemma
statement; the second is easily derived from the first with the
use of proposition \ref{p4}. The equality \eqref{eq4.13} indicates
that there are no connections between the blocks
$[U\mathfrak{E}^T_\Sigma U^{-1}]_j$ under the assumed coordinate
constraints.

If the coordinates of the point $\mathbf{r}$ take extreme values,
then the equality \eqref{eq4.12} is violated: in matrix algebra
$[U\mathfrak{E}^T_\Sigma U^{-1}](\mathbf{r})$ some links between
its blocks may appear. It is these links that distinguish the left
and right-hand sides in the embedding \eqref{eq4.10}: the right
one consists of standard algebras and has no such links. Let us
explain this by examples.

From the definition and properties of the functions
$\tau^i_{\gamma \Phi^j}$ (see \ \eqref{eq2.19}) it follows that
the equality $\tau_{\gamma
\Phi^j}^i(r_j)=\tau_{\gamma\Phi^{j'}}^{i'}(r_{j'})$ is possible
only at extreme values of coordinates, i.e. \,\ at
$r_j\in\{0,\epsilon_j\}$ and $r_{j'}\in\{0,\epsilon_{j'}\}$. Let
the functions $\tau^i_{\gamma\Phi^j}$ be such that
$\tau^i_{\gamma\Phi^j}(\epsilon_j)=\tau^{i+1}_{\gamma\Phi^j}(\epsilon_j)=\tau$
holds. In this case, the $j$block of the eikonal in \eqref{eq4.8}
looks like
$$
\bigl[[UE^T_\gamma U^{-
1}](\mathbf{r})\bigr]_j=\tau^1_{\gamma\Phi^j}(\epsilon_j)
P^1_{\gamma \Phi^j}+\dots+\tau( P^i_{\gamma\Phi^j}+
P^{i+1}_{\gamma \Phi^j}) +\dots+\tau^{n_{\gamma\Phi^j}}_{\gamma
\Phi^j}(\epsilon_j)P^{n_{\gamma\Phi^j}}_{\gamma \Phi^j}
$$
Therefore, the projectors $P^i_{\gamma \Phi^j}$ and
$P^{i+1}_{\gamma \Phi^j}$ do not enter into the family of algebra
generators $[U\mathfrak{E}^T_\Sigma U^{-1}](\mathbf{r}$
separately, but as a sum $P^i_{\gamma \Phi^j}+ P^{i+1}_{\gamma
\Phi^j}$ and thus the number of generators decreases by  $1$. As a
consequence, a violation of \eqref{eq4.12} is possible (and
in substantial cases happens): instead of equality, only
$[U\mathfrak{E}^T_\Sigma
U^{-1}](\mathbf{r})\subset\bigoplus_{j=1}^J \mathfrak{P}_{\Phi^j}$
is guaranteed.

Similarly, projectors belonging to different blocks can be linked.
Let $\tau^i_{\gamma \Phi^j}$ be such that
$\tau^i_{\gamma\Phi^j}(\epsilon_j)=\tau^{i'}_{\gamma\Phi^{j'}}(\epsilon_{j'})=\tau$
for some different $j$ and $j'$. In this case, in \eqref{eq4.8}
the blocks with numbers $j$ and$j'$ look like
\begin{align*}
\bigl[[UE^T_\gamma U^{- 1}](\mathbf{r})\bigr]_j
&=\tau^1_{\gamma\Phi^j}(\epsilon_j) P^1_{\gamma \Phi^j}+\dots+\tau
P^i_{\gamma\Phi^j} +\dots+\tau^{n_{\gamma\Phi^j}}_{\gamma
\Phi^j}(\epsilon_j) P^{n_{\gamma\Phi^j}}_{\gamma \Phi^j},
\\
\bigl[[U\, E^T_\gamma U^{-1}](\mathbf{r})\bigr]_{j'}
&=\tau^1_{\gamma \Phi^{j'}}(\epsilon_{j'})P^1_{\gamma
\Phi^{j'}}+\dots+\tau P^{i'}_{\gamma \Phi^{j'}}
+\dots+\tau^{n_{\gamma\Phi^{j'}}}_{\gamma\Phi^{j'}}(\epsilon_{j'})
P^{n_{\gamma\Phi^{j'}}}_{\gamma\Phi^{j'}}
\end{align*}
and turn out to be connected (via the projectors $P^i_{\gamma
\Phi^j}$ and $P^{i'}_{\gamma \Phi^{j'}}$), which can also lead
to a decrease in the number of generators of algebra
$[U\mathfrak{E}^T_\Sigma U^{-1}](\mathbf{r})$.

\subsection{Reducibility}
\label{ss4.3}

Algebras $\mathfrak{P}_{\Phi^j}$ defined by projector sets
$\mathbb P_{\Phi^j}$ in \eqref{eq4.9} are, generally speaking,
reducible: by the proposition \ref{p3} we have
$$
\mathfrak{P}_{\Phi^j} = \bigoplus_{k=1}^{q_j}\mathfrak{P}_{\Phi^j}^k,\qquad
\mathfrak{P}_{\Phi^j}^k \cong \mathbb M^{\varkappa_{j,k}}, \quad
\varkappa_{j,1}+\dots+ \varkappa_{j,q_k}\le m_{\Phi^j}
$$
and the embedding \eqref{eq4.10} takes the form
\begin{equation}
\label{eq4.14}
U\mathfrak{E}^T_\Sigma{ U^{-1}}\subset\bigoplus_{j=1}^J C\biggl([0,\epsilon_j];\bigoplus_{k=1}^{q_j}\mathfrak{P}_{\Phi^j}^k\biggr)=
\bigoplus_{j=1}^J\bigoplus_{k=1}^{q_j}C([0,\epsilon_j];\mathfrak{P}_{\Phi^j}^k)
\end{equation}
with \textit{unreducible} $\mathfrak{P}_{\Phi^j}^k$. In
the algebraic theorem \ref{t1} below, among other results, the
following equality will be established
\begin{equation}
\label{eq4.15}
\mathfrak{P}_{\Phi^j}^k = \vee \mathbb P_{\Phi^j}^k,
\end{equation}
where $\mathbb P_{\Phi^j}^k$ is some subset of $\mathbb P_{\Phi^j}$ such that
\begin{equation}
\label{eq4.16} \bigcup_{k=1}^{q_j} \mathbb P_{\Phi^j}^k=\mathbb
P_{\Phi^j},\qquad \mathbb P_{\Phi^j}^k \cap \mathbb
P_{\\Phi^j}^{k'}=\varnothing\quad\text{text}\quad k\neq k'.
\end{equation}
{\bad Thus, reducing the algebra $U\mathfrak{E}^T_\Sigma{ U^{-1}}$
to the sum of irreducible blocks in \eqref{eq4.14} will turn
into an adequate grouping of projectors within each of} $\mathbb
P_{\Phi^j}$ sets. When grouping, it is convenient to go to new
numbering.
\medskip

The first step is to move to the continuous numbering of the
algebras entering in \eqref{eq4.14}:
\begin{align}
\notag &\mathfrak{P}_{\Phi^1}^1,\, \dots,\,
\mathfrak{P}_{\Phi^1}^{q_1};\quad \mathfrak{P}_{\Phi^2}^1,\,
\dots,\, \mathfrak{P}_{\Phi^2}^{q_2} \quad\dots;\quad
\mathfrak{P}_{\Phi^J}^1,\, \dots,\, \mathfrak{P}_{\Phi^J}^{q_J}
\\
&\qquad\Longrightarrow \mathfrak{P}_1,\, \dots,\,
\mathfrak{P}_{q_1};\quad \mathfrak{P}_{q_1+1},\, \dots,\,
\mathfrak{P}_{q_1+q_2} \quad \dots;\quad
\mathfrak{P}_{q_1+s+q_{J-1}},\, \dots,\, \mathfrak{P}_L,
\label{eq4.17}
\end{align}
where $L:={q_1+\dots+q_{J}}$. Similarly {(\,\,\ by the formal
replacement of $\mathfrak{P}$ by $\mathbb P$ in \eqref{eq4.17})},
let us go from sets $\mathbb P_{\Phi^j}^k$ to sets $\mathbb
P_{l}$, where $l=1,\dots, L$. Each set $\mathbb P_{l}$, in turn,
can be broken down into subsets corresponding to individual
vertices of $\gamma\in\Sigma$:
$$
\mathbb P_l=\bigcup_{\gamma\in\Sigma}\mathbb P_{l}^{\gamma},\qquad
\mathbb P_{l}^{\gamma} =\{P^i_{\gamma' \Phi^j}\in\mathbb P_{l}\mid
\gamma'=\gamma\}=\{P^{i_1}_{\gamma\Phi^j},\dots,P^{i_{n_{\gamma
l}}}_{\gamma \Phi^j}\}
$$
where $n_{\gamma l}:=\#\mathbb P_{l}^{\gamma}$. Finally, we
renumber the projectors inside each $\mathbb P_{l}^{\gamma}$:
$$
P^{i_1}_{\gamma \Phi^j},\ \dots,\ P^{i_k}_{\gamma \Phi^j},\ \dots,
\ P^{i_{n_{\gamma l}}}_{\gamma \Phi^j}\quad \Longrightarrow\quad
P^1_{\gamma l},\ \dots,\ P^{k}_{\gamma l},\dots,\ P^{n_{\gamma
l}}_{\gamma l}
$$
and the functions $\tau_{\gamma \Phi^j}^i$ and projectors
$\mathcal P_{\gamma \Phi^j}^i$:
\begin{align*}
&\tau^{i_1}_{\gamma \Phi^j},\, \dots,\, \tau^{i_k}_{\gamma
\Phi^j},\,\dots,\, \tau^{i_{n_{\gamma l}}}_{\gamma \Phi^j} \quad
\mathcal P^{i_1}_{\gamma \Phi^j},\, \dots,\, \mathcal
P^{i_k}_{\gamma \Phi^j},\, \dots,\, \mathcal P^{i_{n_{\gamma
l}}}_{\gamma \Phi^j}
    \\
&\qquad\Longrightarrow \tau^1_{\gamma l},\, \dots,\,
\tau^{k}_{\gamma l},\, \dots,\, \tau^{{n_{\gamma l}}}_{\gamma l}
\quad \mathcal P^1_{\gamma l},\, \dots,\, \mathcal P^{k}_{\gamma
l},\, \dots,\, \mathcal P^{n_{{\gamma l}}}_{\gamma l}.
\end{align*}
In new numbering, the embedding \eqref{eq4.14} is written as
\begin{equation}
    \label{eq4.18}
U\mathfrak{E}^T_\Sigma{ U^{-1}}\subset\bigoplus_{l=1}^{L}
C([0,\varepsilon_l];\mathfrak{P}_{l}),
\end{equation}
where $\varepsilon_1=\dots =\varepsilon_{q_1}=\epsilon_1$;
$\varepsilon_{q_1+1}=\dots =\varepsilon_{q_1+q_2}=\epsilon_2$;
$\dots$, and the eikonal representation \eqref{eq3.10} and
relation \eqref{eq4.12} take the form consistent
with \eqref{eq4.18}:
\begin{equation}
    \label{eq4.19}
UE^T_\gamma U^{- 1}= \bigoplus\sum_{l=1}^L
\biggl[\sum_{k=1}^{n_{\gamma l}}\tau^k_{\gamma l} P^k_{\gamma
l}\biggr] =\sum_{l=1}^L \biggl[\sum_{k=1}^{n_{\gamma
l}}\tau^k_{\gamma l} \mathcal P^k_{\gamma l}\biggr]
\end{equation}
and
$$
(U\mathfrak{E}^T_\Sigma U^{-1})(\mathbf{r})=\bigoplus_{l=1}^L \mathfrak{P}_l,
$$
where the coordinates of the point $\mathbf{r} =
\{r_1,\dots,r_{L}\}$ do not take extreme values.

Note that it is at the renumbering that the $\gamma$ index for the
functions $\tau^i_\Phi$ becomes necessary (see \ \eqref{eq2.26}).
This is because the originally equal functions
$\tau_{\gamma \Phi^j}^i=\tau_{\gamma'\Phi^j}^i$ for a fixed family
$\Phi^j$ can correspond to the projectors $P_{\gamma \Phi^j}^i$
and $P_{\gamma'\Phi^j}^i$, which, after renumbering, end up in 
different blocks $\mathfrak{P}_l$ and $\mathfrak{P}_{l'}$.

As the corollary \ref{c2} and the following comments show,
connections between blocks of the algebra $U\mathfrak{E}^T_\Sigma{
U^{-1}}$ are possible only on the boundaries of intervals
$[0,\varepsilon_j]$. For a detailed description of these
connections, it is convenient to use the following formalism.
\medskip

Fix the vertex $\gamma\in\Sigma$ and consider the set $(k,l,r_l)$
corresponding to the value $\tau_{\gamma l}^k(r_l)$ of the
function $\tau_{\gamma l}^k$. We say that the sets $(k,l,r_l)$ and
$(k',l',r_{l'})$ are related and write $(k,l,r_l)\leftrightarrow
(k',l',r_{l'})$ if $\tau_{\gamma l}^k(r_l)=\tau_{\gamma
l'}^{k'}(r_{l'})$ is satisfied. By the properties \eqref{eq2.19}
and \eqref{eq2.20}, such equalities are possible only for extreme
values of the parameters $r_l\in\{0,\varepsilon_l\}$ and
$r_{l'}\in\{0,\varepsilon_{l'}\}$; only such values occur in the
proposition below. These same properties lead to the following
properties of the $\leftrightarrow$-- relation.
\begin{proposition}
            \label{p6}
For a given set $(k,l,r_l)$, one and only one of the following
conditions is satisfied:
\begin{itemize}
\item[1)] there is no set $(k,l,r_l)$ other than $(k,l',r_{l'})$
such that $(k,l,r_l)\leftrightarrow(k',l',r_{l'})$;

\item[2)] {\bad there is a single set $(k',l',r_{l'})$ such that
$(k,l,r_l)\leftrightarrow(k',l',r_{l'})$}, with $l'=l$,
$r_{l'}=r_{l}$, $k'\neq k$;

\item[3)] {\bad there exists a single set $(k',l',r_{l'})$ such
that $(k,l,r_l)\leftrightarrow(k',l',r_{l'})$}, with $l'\neq l$.
\end{itemize}
 \end{proposition}

In accordance with proposition \ref{p6}, let us segregate the sets
$(k,l,r_l)$ into their corresponding types $\mathbf 1$, $\mathbf
2$, and $\mathbf 3$.

\begin{lemma}
\label{l3} Let $\mathbf{r}=\{r_1,\dots, r_l,\dots,r_\mathcal L\}$,
$r_l\in\{0,\varepsilon_l\}$.

1. If $(k,l,r_l)\in \mathbf 1$, then for any
$\widetilde{\mathbf{r}}$ with  coordinate $\widetilde r_{l}\neq
r_l$ there exists an element $e\in U\mathfrak{E}^T_\Sigma U^{-1}$
such that
$$
e(\mathbf{r})={\mathcal P}^k_{\gamma l},\qquad
e(\widetilde{\mathbf{r}})=O.
$$

2.  If $(k,l,r_l)\in \mathbf{2}$ and
$(k,l,r_l)\leftrightarrow(k',l,r_l)$, then for any
$\widetilde{\mathbf{r}}$ with $\widetilde r_{l}\neq r_l$ there
exists an element $e\in U\mathfrak{E}^T_\Sigma U^{-1}$ such that
$$
e(\mathbf{r})={\mathcal P}^k_{\gamma l}+{\mathcal P}^{k'}_{\gamma
l},\qquad e(\widetilde{\mathbf{r}})=O.
$$

3. 
{\bad If $(k,l,r_l)\,{\in}\, \mathbf 3$ 
and $(k,l,r_l)\,{\leftrightarrow}\,(k',l',r_{l'})$,
then for $\mathbf{r}\,{=}\,\{r_1,\dots, r_l,\dots,r_{l'},\dots,r_\mathcal L\}$} 
and any $\widetilde{\mathbf{r}}$ with coordinates $\widetilde
r_{j}\neq r_j$, $\widetilde r_{j'}\neq r_{j'}$ there exists an
element $e\in U\mathfrak{E}^T_\Sigma U^{-1}$ such that
$$
e(\mathbf{r})={\mathcal P}^k_{\gamma l}+{\mathcal P}^{k'}_{\gamma
l'},\qquad e(\widetilde{\mathbf{r}})=O.
$$
\end{lemma}

\begin{proof}[{\rm is quite similar to the proof of the lemma \ref{l2}: like the latter,
it reduces to choosing a suitable polynomial $q$. The reason for
the ``gluing'' of the projectors in sum is the same as in the
examples after the corollary \ref{c2}}].
\end{proof}

As noted, the difference between the algebra
$U\mathfrak{E}^T_\Sigma U^{-1}$ and $\bigoplus_{l=1}^L
C([0,\varepsilon_l]; \mathfrak{P}_l)$ consists of possible
connections between matrix algebra blocks $(U\mathfrak{E}^T_\Sigma
U^{-1})(\mathbf{r})$ that may appear when the coordinates of $r_l$
take extreme values. The \textit{boundary algebra} that will be
used to study these connections, is
\begin{equation}
\label{eq4.20}
\partial (U\mathfrak{E}^T_\Sigma U^{-1}):=\bigl\{ e(\mathbf{0})\oplus e(\boldsymbol{\varepsilon})\bigm|
e\in U\mathfrak{E}^T_\Sigma U^{-1}\bigr\} \subset \biggl[ \bigoplus_{l=1}^L \mathfrak{P}_l\biggr]\oplus\biggl[ \bigoplus_{l=1}^L \mathfrak{P}_l\biggr],
\end{equation}
where $\mathbf{0}=\{0,\dots,0\}$ and
$\boldsymbol{\varepsilon}=\{\varepsilon_1,\dots,\varepsilon_L\}$.
To describe them, let us first consider some general results
concerning the structure of matrix algebras of
type $\bigoplus_{l=1}^L \mathfrak{P}_l$.

\subsection({Algebras generated by one-dimensional projectors}){{\bad Algebras generated by one-dimensional projectors}} 
\label{ss4.4}

{\bad Let} 
$\mathscr{G}$ be a Hilbert space with the inner product
$\langle\,{\cdot}\,,{\cdot}\,\rangle$, and let a set of
one-dimensional projectors $P^1,\dots,P^n:
P^i=\langle\,{\cdot}\,,\beta^i\rangle\beta^i$, $\|\beta^i\|=1$, be
given; denote
$$
\mathscr A:=\operatorname{span}\{\beta^1,\dots,\beta^n\},\qquad \mathfrak P:=\vee\{P^1,\dots,P^n\}\subset\mathfrak{B}(\mathscr{G}).
$$
Let us endow the set $B:=\{\beta^1,\dots,\beta^n\}$ with a
reflexive and symmetric relation by putting
$\beta^i\sim_0\beta^{i'}$ if $\langle\beta^i,\beta^{i'}\rangle\neq
0$. It determines the equivalence: we put $\beta^i\sim\beta^{i'}$
if there are vectors $\beta^{i_1},\dots,\beta^{i_k}$ such that
$\beta^i \sim_0\beta^{i_1}\sim_0 \cdots\sim_0\beta^{i_k}\sim_0
\beta^{i'}$. Then we transfer this equivalence to the projectors
by taking $P^i\sim P^{i'}$ if $\beta^i\sim\beta^{i'}$.

Represent $B=B_1\,{\cup}\,{\cdots}\,{\cup}\, B_q$ 
as the partition on the equivalence classes and denote $\mathscr
A_k:=\operatorname{span}B_k$. To obtain such a representation one
can use the procedure \eqref{eq2.11}, \eqref{eq2.12}.

The definition of $\sim$ implies $\mathscr A_k\perp \mathscr A_l$
for $k\neq l$. Thus we have a decomposition $\mathscr A=\mathscr
A_1\oplus\dots\oplus \mathscr A_q$, which obviously reduces all
projectors $P^i$.

\begin{proposition}
\label{p7}
The algebra $\mathfrak P$ is reduced by subspaces $\mathscr A_k$, and the expansion
\begin{equation}
\label{eq4.21} \mathfrak P=\bigoplus_{k=1}^q\mathfrak
P|_{\mathscr{A}_k}
\end{equation}
holds, where $\mathfrak
P|_{\mathscr{A}_k}\cong\mathbb{M}^{\varkappa_k}$,
$\varkappa_k=\operatorname{dim}\mathscr A_k$.
\end{proposition}

Reducibility is obvious; the equality for the dimensions follows
from the fact that any $\beta^i\in\mathscr A_k$ is cyclic
in $\mathscr A_k$ for the part $\mathfrak{P}|_{\mathscr A_k}$. 
\medskip

The following considerations model the situation that will be
encountered in the study of eikonal algebra. Namely, the possible
connections between its blocks in the representation \eqref{eq4.8}
are discussed in abstract form.

Let $\mathscr{G}_k$, $k=1,2,3$ be three Hilbert spaces, each
containing a set of one-dimen\-sional projectors $P^1_k,
\dots,P^{n_k}_k\colon
P^i_k=\langle\,{\cdot}\,,\beta^i_k\rangle\beta^i_k$,
$\|\beta^i_k\|=1$, where $\beta^i_k$ are vectors from the sets
$$
B_k:=\{\beta^i_k\mid i=1,\dots,{n}_k\}\subset \mathscr{G}_k.
$$
The projectors generate the algebras
$$
\mathfrak{P}_1=\vee\{P^1_1,\dots,P_1^{n_1}\},\qquad \mathfrak{P}_2=\vee\{P_2^1,\dots,P_2^{n_2}\},\qquad \mathfrak{P}_3=\vee\{P_3^1,\dots,P_3^{n_3}\}.
$$
Let us construct the algebra
\begin{equation}
\label{eq4.22} \mathfrak{P} :=
\mathfrak{P}_1\oplus\mathfrak{P}_2\oplus\mathfrak{P}_3
\subset\mathfrak{B}(\mathscr{G}_1\oplus
\mathscr{G}_2\oplus\mathscr{G}_3),
\end{equation}
{\bad which is a subalgebra of the algebra of bounded operators acting in the space $\mathscr{G}_1\oplus\mathscr{G}_2\oplus\mathscr{G}_3$,} 
with the  generators
\begin{equation}
\begin{gathered}
\mathcal{P}_1^i:=P_1^i \oplus O_2 \oplus O_3,\qquad i=1,\dots,n_1,
\\
\mathcal{P}_2^i:= O_1\oplus P^i_2\oplus O_3,\qquad i=1,\dots,n_2,
\\
\mathcal{P}_3^i:= O_1\oplus O_2 \oplus P^i_3,\qquad i=1,\dots,n_3,
\end{gathered}
\label{eq4.23}
\end{equation}
where $O_k$ is the null operator acting in the $k$th component of
the space $\mathscr{G}_1\oplus\mathscr{G}_2\oplus\mathscr{G}_3$.
The algebras $\mathfrak{P}_k$ in \eqref{eq4.22} will be referred
to as \textit{blocks} of the algebra $\mathfrak{P}$. Note that
in the current considerations, roughly speaking, the algebras
$\mathfrak{P}_1$ and $\mathfrak{P}_2$ correspond to some pair of
selected blocks in \eqref{eq4.8}, whereas $\mathfrak{P}_3$ is
``all the rest''.

Let us say that an algebra $\mathfrak{Q}\subset\mathfrak{P}$
\textit{separates} (does not connect) the blocks $\mathfrak{P}_1$
and $\mathfrak{P}_2$ in \eqref{eq4.22}, if along with any element
$q_1\oplus q_2\oplus q_3\in\mathfrak{Q}$ it contains elements
$q_1\oplus O_2\oplus q_3'$ and $O_1\oplus q_2\oplus q_3''$, where
$q_3'$, $q_3''$ are some elements $\mathfrak{P}_3$. Otherwise, we
say that these blocks are \textit{connected}. Similarly, a
connection (or lack thereof) is defined for any pair of blocks
in \eqref{eq4.22}.

Note the obvious fact: if algebra $\mathfrak{Q}$ admits a system
of generators, each of which is of the form either $q_1\oplus
O_2\oplus q_3'$ or $O_1\oplus q_2\oplus q_3''$, then it does
separate the blocks $\mathfrak{P}_1$ and $\mathfrak{P}_2$.
\medskip

Let
\begin{equation}
\label{eq4.24}
\mathbb P:=\{\mathcal{P}_k^i\mid i=1,\dots,{n}_k;\, k=1,2,3\}
\end{equation}
be the complete set of generators of algebra $\mathfrak{P}$. Let
$\mathcal{T}\colon \mathbb P\to\mathbb P$ be a map (involution) on
it such that the following holds: if
$\mathcal{T}(\mathcal{P}_k^i)=\mathcal{P}_{k'}^{i'}$ then
$\mathcal{T}(\mathcal{P}_{k'}^{i'})=\mathcal{P}_{k}^i$ and one
(and only one) of the conditions is satisfied:
$$
\mathcal{T}(\mathcal{P}_k^i)=\mathcal{P}_k^i\quad\text{or}\quad
\mathcal{T}(\mathcal{P}_k^i)\mathcal{P}_k^i=\mathcal{P}_k^i\mathcal{T}(\mathcal{P}_k^i)=0.
$$
This map induces a partition of the set $\mathbb{P}$ into pairs
$\{\mathcal P,\mathcal T(\mathcal P)\}$ with the components
in each pair either identical or orthogonal to each other. It can
be seen that such $\mathcal{T}$ do exist, and in the plural.
However, in the eikonal algebra the map $\mathcal{T}$ will be
quite concrete and determined by the values of the functions
$\tau^k_{\gamma l}(r)$ at $r=0$ and $r=\varepsilon$.

Using this map, let us determine the projectors $\mathcal{Q}_k^i$:
$$
\mathcal{Q}_k^i := \begin{cases}
\mathcal{P}_k^i, &\text{if }\mathcal{T}(\mathcal{P}_k^i)=\mathcal{P}_k^i,
\\
\mathcal{P}_k^i+\mathcal{T}(\mathcal{P}_k^i), &\text{if } \mathcal{T}(\mathcal{P}_k^i)\mathcal{P}_k^i=0,
\end{cases}
$$
among which there may be identical ones. If $\mathcal{T}$ is not
identical, then due to coincidences their total number will
obviously be less than $n_1+n_2+n_3$. Let us form the algebra
\begin{equation}
\label{eq4.25}
\mathfrak{Q}:=\vee\{\mathcal{Q}_k^i\mid i=1,\dots,{n}_k;\,k=1,2,3\}\subset \mathfrak{P}.
\end{equation}
It is determined by the map $\mathcal T$. However, different
admissible $\mathcal T$ and $\mathcal T'$ can correspond to the
same algebra $\mathfrak{Q}$. This fact is used below in the proof
of theorem \ref{t1}.
\medskip

Next, we discuss under what conditions the algebra $\mathfrak{Q}$
introduced in this way, separates or connects the blocks
$\mathfrak{P}_1$ and $\mathfrak{P}_2$.

Let us divide the set of projectors $\mathbb{P}$ defined
in \eqref{eq4.24} into parts $\mathbb{P}_k:=\{\mathcal{P}_k^i\mid
i=1,\dots,{n}_k\}$, $k=1,2,3$. To each part, we correspond the
matrices
$$
G(\mathbb P_k):=\{\|\mathcal{P}^i_k\mathcal{P}^j_k\|\}_{i,j=1}^{n_k}=
\begin{pmatrix}
\|\mathcal{P}^1_k\mathcal{P}^1_k\| &\cdots &\|\mathcal{P}^1_k\mathcal{P}^{n_k}_k\| \\
\vdots & \ddots &\vdots\\
\|\mathcal{P}^{n_k}_k\mathcal{P}^1_k\| & \cdots &
\|\mathcal{P}^{n_k}_k\mathcal{P}^{n_k}_k\|
\end{pmatrix}
$$
and ${G}(\mathcal{T}(\mathbb{P}_k))=\{\|\mathcal
T(\mathcal{P}^i_k)\mathcal T(\mathcal{P}^j_k)\|\}_{i,j=1}^{n_k}$.
Let
$$
B:=\{\beta^i_k\mid i=1,\dots,{n}_k;\,k=1,2,3\} = B_1\cup B_2\cup
B_3,\qquad B_k=\{\beta^i_k\mid i=1,\dots,{n}_k\}.
$$
By virtue of the obvious equalities
$\|\mathcal{P}^i_k\,\mathcal{P}^j_{k}\|=|\langle\beta^i_k,\beta^j_{k}\rangle|$
we have
$$
G(\mathbb{P}_k)= \begin{pmatrix}
|\langle\beta_k^1,\beta_k^1\rangle| &\cdots &||\langle\beta_k^1,\beta_k^{{n}_k}\rangle| \\
\vdots & \ddots &\vdots\\
|\langle\beta_k^{{n}_k},\beta_k^1\rangle| & \cdots &
|\langle\beta_k^{{n}_k},\beta_k^{{n}_k}\rangle|
\end{pmatrix}.
$$
The following result about the connection of blocks is used later in
the study of the structure of eikonal algebra.

\begin{theorem}
\label{t1} Let the blocks $\mathfrak{P}_1$ and $\mathfrak{P}_2$ of
the algebra \eqref{eq4.22} be such that each of their
corresponding sets $B_1$ and $B_2$ is an equivalence class with
respect to relation $\sim$. Then the algebra $\mathfrak{Q}$ can
connect $\mathfrak{P}_1$ and $\mathfrak{P}_2$ only if
$$
\mathcal{T}(\mathbb{P}_1)=\mathbb{P}_2\quad\text{and}\quad
G(\mathbb{P}_1)=G(\mathcal{T}(\mathbb{P}_1))
$$
holds. If these conditions are satisfied, the following relations
are valid
$$
\operatorname{dim}\operatorname{span}
B_1=\operatorname{dim}\operatorname{span}B_2=:l, \qquad
\mathfrak{P}_1\cong\mathfrak{P}_2\cong\mathbb M^{l},
$$
and the algebra $\mathfrak{Q}$ is represented as
$$
\mathfrak{Q}=\mathfrak{Q}_{12}\oplus\mathfrak{Q}_3,
$$
where
$$
\mathfrak{Q}_{12}=\vee\{\mathcal{P}\oplus
\mathcal{T}(\mathcal{P})\mid
\mathcal{P}\in\mathbb{P}_1\}\subset\mathfrak{P}_1\oplus\mathfrak{P}_2,\qquad
\mathfrak{Q}_3\subseteq\mathfrak{P}_3.
$$
Moreover, the algebra $\mathfrak{Q}$ separates $\mathfrak{Q}_{12}$
and $\mathfrak{Q}_3$.
\end{theorem}


\begin{proof}
1. For $k,k'=1,2,3$, $k'\neq k$, put
$$
\mathbb{Q}_{k k'}:= \{\mathcal{Q}_k^i=\mathcal{P}_k^i+\mathcal{T}(\mathcal{P}_k^i)\mid \mathcal{T}(\mathcal{P}_k^i)\in\mathbb P_{k'}\}
$$
and note the equality following from the involution property of
$\mathcal T$:
$$
\mathbb{Q}_{k' k} =
\{\mathcal{Q}_{k'}^{i'}=\mathcal{P}_{k'}^{i'}+\mathcal{T}(\mathcal{P}_{k'}^{i'})\mid
\mathcal{T}(\mathcal{P}_{k'}^{i'})\in\mathbb P_{k}\} =
\mathbb{Q}_{k k'}.
$$
When $k=k'=1,2,3$, we put
$$
\mathbb{Q}_{k k}:= \Biggl\{\mathcal{Q}_k^i=
\begin{cases}
\mathcal{P}_k^i+\mathcal{T}(\mathcal{P}_k^i), &\text{if }\mathcal{P}_k^i\ne \mathcal{T}(\mathcal{P}_k^i)\in\mathbb P_{k},
\\
\mathcal{P}_k^i, &\text{if }\mathcal{P}_k^i=\mathcal{T}(\mathcal{P}_k^i),
\end{cases}\Biggm| i=1,\dots,n_k\Biggr\}\Biggr\}.
$$
Thus, all $\mathbb{Q}_{k k'}$ consist of one- and two-dimensional
projectors, and the algebra $\mathfrak{Q}$ is represented as
\begin{equation}
\label{eq4.26}
\mathfrak{Q}=\vee[\mathbb{Q}_{11}\cup\mathbb{Q}_{22}\cup\mathbb{Q}_{33}\cup\mathbb{Q}_{12} \cup\mathbb{Q}_{13}\cup\mathbb{Q}_{23}].
\end{equation}
In view of the form of the projectors \eqref{eq4.23}, it is easy
to see that the algebra $\mathfrak{Q}$ can connect  the blocks
$\mathfrak{P}_1$ and $\mathfrak{P}_2$ only if $\mathbb{Q}_{12}\neq
\varnothing$ holds; otherwise, it separates them.

2. Suppose that algebra $\mathfrak{Q}$ admits the representation
\begin{equation}
\label{eq4.27}
\mathfrak{Q}=\vee [\mathbb{Q}_{11}'\cup\mathbb{Q}_{22}'\cup\mathbb{Q}_{33}\cup\mathbb{Q}_{13}\cup \mathbb{Q}_{23}],
\end{equation}
where
$$
\mathbb{Q}_{11}':=\mathbb{Q}_{11}\cup\,\{\mathcal{P}_1^i\mid
\mathcal{P}_1^i+\mathcal
T(\mathcal{P}_1^i)\in\mathbb{Q}_{12}\},\qquad \mathbb{Q}_{22}':
=\mathbb{Q}_{22}\cup\,\{\mathcal{P}_2^i\mid
\mathcal{P}_2^i+\mathcal T(\mathcal{P}_2^i)\in\mathbb{Q}_{12}\}.
$$
Then it separates the blocks $\mathfrak{P}_1$ and
$\mathfrak{P}_2$. Indeed, in this case,
in the algebra $\mathfrak{Q}$ along with elements $\mathcal
Q^i_1=\mathcal{P}_1^i+\mathcal T(\mathcal{P}_1^i)$ and $\mathcal
Q^i_2=\mathcal{P}_2^i+\mathcal T(\mathcal{P}_2^i)$ all the
projectors $\mathcal{P}_1^i={P}^i_1\oplus O_2\oplus O_3$ and
$\mathcal{P}_2^i=O_1\oplus{P}^i_2\oplus O_3$ will enter
separately. By introducing them into the set of generators of the
algebra ${\mathfrak{Q}}$ instead of the elements of $\mathbb
Q_{12}$, it is easy to see that the separation takes place.

Here we explain the following. The definition of map $\mathcal T$
excludes the simultaneous presence of elements $\mathcal
Q^i_k=\mathcal P^i_k+\mathcal T(\mathcal P^i_k)$ and $\mathcal
P^i_k$ in the list of generators $\{\mathcal{Q}_k^i\mid
i=1,\dots,{n}_k; \,k=1,2,3\}$ of the algebra $\mathfrak Q$ (cf. \
\eqref{eq4.25}). Nevertheless, they may well be contained in the
algebra itself. This corresponds to the possibility to
replace $\mathcal T$ with another map $\mathcal T'\colon
\mathbb{P}\to\mathbb P$ so that its corresponding algebra
$\mathfrak Q'$ is the same as the original $\mathfrak Q$.

3. Let's show that if
$$
\mathbb{Q}_{11}\cup\mathbb{Q}_{13}\cup\mathbb{Q}_{22}\cup\mathbb{Q}_{23}\neq
\varnothing,
$$
then for the algebra $\mathfrak{Q}$ the representation
\eqref{eq4.27} is valid and hence it separates the blocks
$\mathfrak{P}_1$ and $\mathfrak{P}_2$. The following
considerations are quite analogous for the parts
$\mathbb{Q}_{11}\cup\mathbb{Q}_{13}$ and
$\mathbb{Q}_{22}\cup\mathbb{Q}_{23}$ due to the analogy of
$\mathfrak{P}_1$ and $\mathfrak{P}_2$. Let us perform them for the
case $\mathbb{Q}_{11}\cup\mathbb{Q}_{13}\neq\varnothing$. In this
case we have $\mathbb{Q}_{12}\ne \varnothing$, because in the
opposite case the blocks $\mathfrak{P}_1$ and $\mathfrak{P}_2$ are
definitely separated.

Each element $\mathcal{Q}_1^i\in\mathbb{Q}_{11}\cup\mathbb{Q}_{12}\cup\mathbb{Q}_{13}$ is
\begin{equation}
\label{eq4.28}
\mathcal{Q}_1^i:=\begin{cases}
\mathcal{P}_1^i+\mathcal{T}(\mathcal{P}_1^i), &\text{if } \mathcal{T}(\mathcal{P}_1^i)\neq\mathcal{P}_1^i,
\\
\mathcal{P}_1^i, &\text{if } \mathcal{T}(\mathcal{P}_1^i)= \mathcal{P}_1^i.
\end{cases}
\end{equation}
Each vector $\beta_1^i$ from the set $B_1$ is associated with a
projector $\mathcal{P}_1^i=P_1^i \oplus
O_2\oplus{O}_3=\langle\,{\cdot}\,, \beta^i_1\rangle\beta^i_1\oplus
O_2 \oplus O_3$, which, on its turn, determines the corresponding
projector $\mathcal{Q}_1^i$ of the form \eqref{eq4.28}. This
enables us to define the map $\mathbf{b}\colon B_1\to
\mathbb{Q}_{11}\cup\mathbb{Q}_{12}\cup\mathbb{Q}_{13}$ by the rule
$$
\mathbf{b}(\beta_1^i) := \mathcal{Q}_1^i.
$$
Note that if $\mathcal{Q}_1^i=\mathcal P^i_1+\mathcal T(\mathcal
P^i_1)=\mathcal P^i_1+\mathcal P^{i'}_1\in\mathbb{Q}_{11}$ then
there are two vectors $\beta_1^i$, $\beta_1^{i'}$ such that
$\mathbf{b}(\beta_1^i)=\mathbf{b}(\beta_1^{i'})=\mathcal{Q}_1^i=\mathcal{Q}_1^{i'}$.
In what follows, $\mathbf{b}^{-1}(\,\cdot\,)$ denotes the complete
preimage.

By the conditions of the lemma, $B_1$ is an equivalence class with
respect to $\sim$. The condition
$\mathbb{Q}_{11}\cup\mathbb{Q}_{13}\neq \varnothing$ follows, that
there is a pair of vectors
$\beta_1^i\in\mathbf{b}^{-1}(\mathbb{Q}_{12})$ and
$\beta_1^{i'}\in\mathbf{b}^{-1}(\mathbb{Q}_{11}\cup\mathbb{Q}_{13})$
such that $\langle\beta_1^i,\beta_1^{i'}\rangle\neq 0$. In fact,
the absence of such a pair would mean that
$$
\operatorname{span}B_1=\operatorname{span}\mathbf{b}^{-1}(\mathbb Q_{12})\oplus\,\operatorname{span}\mathbf{b}^{-1}(\mathbb Q_{11}\cup\mathbb Q_{13}),
$$
which is impossible by the definition of the  equivalence
${\sim}$.

\goodbreak 

The chosen pair of vectors determines the projectors
$\mathcal{Q}_1^i= \mathcal P^i_1+\mathcal T(\mathcal
P^i_1)=\mathbf{b}(\beta^i_1)\in\mathbb Q_{12}$ and
$\mathcal{Q}_1^{i'}=
\mathbf{b}(\beta^{i'}_1)\in\mathbb{Q}_{11}\cup\mathbb{Q}_{13}$,
and, along with   them, the element
$$
\widetilde{\mathcal{Q}}_1^i:= \mathcal{Q}_1^i\mathcal{Q}_1^{i'}\mathcal{Q}_1^i\in \mathfrak{Q}.
$$
According to \eqref{eq4.22} we have the representations
\begin{gather*}
\mathcal{Q}_1^i = (\mathcal{Q}_1^i)_1\oplus(\mathcal{Q}_1^i)_2 \oplus (\mathcal{Q}_1^i)_3,\qquad \mathcal{Q}_1^{i'} = (\mathcal{Q}_1^{i'})_1\oplus(\mathcal{Q}_1^{i'})_2 \oplus (\mathcal{Q}_1^{i'})_3,
\\
\widetilde{\mathcal{Q}}_1^i = (\widetilde{\mathcal{Q}}_1^i)_1 \oplus(\widetilde{\mathcal{Q}}_1^i )_2 \oplus (\widetilde{\mathcal{Q}}_1^i )_3,
\end{gather*}
where
$(\mathcal{Q}_1^i)_k,(\mathcal{Q}_1^{i'})_k,(\widetilde{\mathcal{Q}}_1^i)_k\in\mathfrak{P}_k$, $k=1,2,3$; this holds
$$
(\widetilde{\mathcal{Q}}_1^i)_k = (\mathcal{Q}_1^i)_k(\mathcal{Q}_1^{i'})_k(\mathcal{Q}_1^i)_k, \qquad k=1,2,3.
$$
By the choice of vector $\beta^i_1$ we have $\mathcal{Q}_1^i=\mathbf{b}(\beta^i_1)\in\mathbb{Q}_{12}$. Therefore, $(\mathcal{Q}_1^i)_3=O_3$ and hence $(\widetilde{\mathcal{Q}}_1^i)_3=O_3$. Similarly, if $\mathcal{Q}_1^{i'}\in\mathbb{Q}_{11}\cup\mathbb{Q}_{13}$, then $(\mathcal{Q}_1^{i'})_2=O_2$, hence $(\widetilde{\mathcal{Q}}_1^i)_2=O_2$. Thus, only the component $(\widetilde{\mathcal{Q}}_1^i)_1$ can be nonzero. In this situation, there are two possibilities.

1)   Let $\mathcal{Q}_1^{i'}\in\mathbb{Q}_{13}$. Then $(\mathcal{Q}_1^{i'})_1= P_1^{i'}$ and $(\mathcal{Q}_1^i)_1= P_1^i$, and $(\widetilde{\mathcal{Q}}_1^i)_1$ is
\begin{equation}
\label{eq4.29}
(\widetilde{\mathcal{Q}}_1^i)_1 = P_1^i P_1^{i'}P_1^i = \langle\beta_1^i,\beta_1^{i'}\rangle^2 P_1^i.
\end{equation}

2) If now $\mathcal{Q}_1^{i'}\in\mathbb{Q}_{11}$. If
$\mathcal{Q}_1^{i'}=\mathcal{P}_1^{i'} $, then considerations
quite analogous to those above lead to the same equality
\eqref{eq4.29}. Consider the case, where
$\mathcal{Q}_1^{i'}=\mathcal{P}_1^{i'} +
\mathcal{T}(\mathcal{P}_1^{i'})$, where
$\mathcal{T}(\mathcal{P}_1^{i'})= \mathcal{P}_1^j\in\mathbb{P}_1$
and $\mathcal{P}_1^j\mathcal{P}_1^{i'}=O_1$. The projector
$\mathcal{P}_1^j$ corresponds to the vector $\beta_1^j\in B_1$.
Then $(\widetilde{\mathcal{Q}}_1^i)_1$ takes the form
\begin{equation}
\label{eq4.30}
(\widetilde{\mathcal{Q}}_1^i)_1 = P_1^i (P_1^{i'}+P_1^j)P_1^i =
[\langle\beta_1^i,\beta_1^{i'}\rangle^2+\langle\beta_1^i,\beta_1^j\rangle^2]P_1^i.
\end{equation}
Comparing \eqref{eq4.29} with \eqref{eq4.30}, we arrive at the
following relation
$$
\widetilde{\mathcal{Q}}_1^i = c \mathcal{P}_1^i,\qquad
c\ge\langle\beta_1^i,\beta_1^{i'}\rangle^2>0.
$$
This means that the algebra $\mathfrak{Q}$ contains separately the
projectors $\mathcal{P}_1^i\in\mathbb{P}_1$ and
$\mathcal{T}(\mathcal{P}_1^i)\in\mathbb{P}_2$. About this result
let us say that the projector $\mathcal
Q_1^i=\mathcal{P}_1^i+\mathcal{T}(\mathcal{P}_1^i)$ is decomposed
into independent (in algebra $\mathfrak{Q}$) one-dimensional parts
$\mathcal{P}_1^i$ and $\mathcal{T}(\mathcal{P}_1^i)$.

Next, consider the map $\mathcal T'\colon
{\mathcal{P}}\to\mathcal{T}$, which differs from $\mathcal{T}$ by
values on only two projectors $\mathcal{P}_1^i$ and
$\mathcal{T}(\mathcal{P}_1^i)$, and put
$$
\mathcal{T}'(\mathcal{P}_1^i):=\mathcal{P}_1^i,\qquad \mathcal{T}'(\mathcal{P}_1^i):=\mathcal{T}(\mathcal{P}_1^i).
$$
Thus the algebras $\mathfrak{Q}$ and $\mathfrak{Q}'$ defined by
the maps $\mathcal{T}$ and $\mathcal{T}'$ obviously coincide, and
$\mathfrak{Q}'$ has its own representation of the form
\eqref{eq4.26}:
$$
\mathfrak{Q}'=\vee[\mathbb{Q}_{11}'\cup\mathbb{Q}_{22}'\cup\mathbb{Q}_{33}'\cup\mathbb{Q}_{12}' \cup\mathbb{Q}_{13}'\cup\mathbb{Q}_{23}']=\mathfrak{Q},
$$
and its relation to \eqref{eq4.26} is as follows:
\begin{gather*}
\mathbb{Q}_{11}'=\mathbb{Q}_{11}\cup\{\mathcal{P}_1^i\}, \qquad \mathbb{Q}_{22}'=\mathbb{Q}_{22}\cup\{\mathcal{T}(\mathcal{P}_1^i)\}, \qquad \mathbb{Q}_{12}'=\mathbb{Q}_{12}\setminus\{\mathcal{Q}_1^i\},
\\
\mathbb{Q}_{13}'=\mathbb{Q}_{13},\qquad\mathbb{Q}_{23}'=\mathbb{Q}_{23}, \qquad\mathbb{Q}_{33}'=\mathbb{Q}_{33}.
\end{gather*}
Thus, the decomposition of $\mathcal Q^i_1$ has resulted in  that
the part $\mathbb Q_{12}$ responsible for linking the blocks
$\mathfrak{P}_1$ and $\mathfrak{P}_2$ has decreased by one
projector.

Repeating the considerations for the part
$\mathbb{Q}_{12}'\subset\mathbb{Q}_{12}$, we see that another
projector can be removed from it too without changing the algebra
$\mathfrak{Q}$. Continuing the procedure for a finite number of
steps will lead to the decay of all the projectors contained
in $\mathbb{Q}_{12}$ and, as a consequence, to representation of
\eqref{eq4.27}.

So, it is shown that the condition
$\mathbb{Q}_{11}\cup\mathbb{Q}_{13}\cup\mathbb{Q}_{22}\cup\mathbb{Q}_{23}=
\varnothing$ is necessary for the algebra to link the blocks
$\mathfrak{P}_1$ and $\mathfrak{P}_2$. Note that the condition
$\mathbb{Q}_{11}\cup\mathbb{Q}_{13}=\varnothing$ is equivalent to
that $\mathcal{T}(\mathbb{P}_1)\subset\mathbb{P}_2$, and the
condition $\mathbb{Q}_{22}\cup\mathbb{Q}_{23}=\varnothing$ --that
$\mathcal{T}(\mathbb{P}_2)\subset\mathbb{P}_1$. Since the map
$\mathcal T$ is involutive, it follows that
$\mathbb{Q}_{11}\cup\mathbb{Q}_{13}\cup\mathbb{Q}_{22}\cup\mathbb{Q}_{23}=
\varnothing$ is satisfied if and only if
$\mathcal{T}(\mathbb{P}_1)=\mathbb{P}_2$.

4. Since then we will assume that the condition
$\mathcal{T}(\mathbb{P}_1)=\mathbb{P}_2$ is satisfied, which means
that for the algebra $\mathfrak{Q}$ connecting blocks
$\mathfrak{P}_1$ and $\mathfrak{P}_2$, the representation
$$
\mathfrak{Q}=\vee[\mathbb{Q}_{12}\cup\mathbb{Q}_{33}]
$$
holds. From this we see that
$\mathfrak{Q}=\mathfrak{Q}_{12}\oplus\mathfrak{Q}_{3} $, where
$\mathfrak{Q}_{12}:=\vee \mathbb{Q}_{12}$, $\mathfrak{Q}_{3} =\vee
\mathbb{Q}_{33}$, with algebra $\mathfrak{Q}$ separating blocks
$\mathfrak{Q}_{12}$ and $\mathfrak{Q}_3$.

Next, we specify the structure of the algebra $\mathfrak{Q}_{12}$.
It is convenient to use a matrix notation for this purpose:
$$
\mathfrak{Q}_{12} = \vee\biggl\{\mathcal{Q}^i:= \begin{pmatrix} P_1^i &\\ & P_2^{i'} \end{pmatrix}\biggm| \mathcal{T}(\mathcal{P}_1^i)=\mathcal{P}_2^{i'};\, i=1,\dots,n_1 \biggr\}
$$
(zero elements omitted) and the representation
$$
\mathcal{Q}^i= \begin{pmatrix}
\langle\,{\cdot}\,,\beta_1^i\rangle\beta_1^i &
\\
& \langle\,{\cdot}\,,\beta_2^{i'}\rangle\beta_2^{i'}
\end{pmatrix}
$$
through the vectors $\beta_k^i$ corresponding to the projectors
$P_1^i$ and $P_2^{i'}$. Now suppose that
${G}(\mathbb{P}_1)\neq{G}(\mathcal{T}(\mathbb{P}_1))$. Note that
this is possible only if each of $B_1$ and $B_2$ has more than one
element. Under this assumption, the algebra $\mathfrak{Q}_{12}$
contains $\mathcal{Q}^i$ and $\mathcal{Q}^j$ such that
$|\langle\beta_1^i,\beta_1^j\rangle\neq\langle\beta_2^{i'},\beta_2^{j'}\rangle|$.
For the product
$\mathcal{Q}^i\mathcal{Q}^j\mathcal{Q}^i\in\mathfrak{Q}_{12}$ we
get the representation
$$
\mathcal{Q}^i\mathcal{Q}^j\mathcal{Q}^i= \begin{pmatrix}
\langle\beta_1^i,\beta_1^j\rangle^2 \langle\,{\cdot}\,,\beta_1^i\rangle \beta_1^i&\\
& \langle\beta_2^{i'},\beta_2^{j'}\rangle^2\langle\,{\cdot}\,, \beta_2^{i'}\rangle\beta_2^{i'}
\end{pmatrix}
$$
and arrive at the relations
\begin{gather*}
\mathcal{Q}^i = \begin{pmatrix} P_1^i &\\ & O_2\end{pmatrix}+
\begin{pmatrix} O_1 &\\ & P_2^{i'}\end{pmatrix}
\\
\mathcal{Q}^i\mathcal{Q}^j\mathcal{Q}^i =
\langle\beta_1^i,\beta_1^j\rangle^2\begin{pmatrix} P_1^i &\\ & O_2
\end{pmatrix}+
\langle\beta_2^{i'},\beta_2^{j'}\rangle^2\begin{pmatrix} O_1 &\\ &
P_2^{i'}\end{pmatrix}
\end{gather*}
From these, considering the inequality
$|\langle\beta_1^i,\beta_1^j\rangle|\neq|\langle\beta_2^{i'},\beta_2^{j'}\rangle|$,
we conclude that the algebra $\mathfrak{Q}_{12}$ contains each of
the projectors
$$
\begin{pmatrix}P_1^i &\\& O_2 \end{pmatrix} \quad\text{and}\quad \begin{pmatrix} O_1 &\\ & P_2^{i'}\end{pmatrix}
$$
separately. Hence, the algebra $\mathfrak{Q}$ separates the blocks
$\mathfrak{P}_1$ and $\mathfrak{P}_2$, and the following
representation holds
$$
\mathfrak{Q}=\mathfrak{P}_1\oplus\mathfrak{P}_2\oplus \mathfrak{Q}_3.
$$

5. So the condition
${G}(\mathbb{P}_1)={G}(\mathcal{T}(\mathbb{P}_1))$ is also
necessary for the algebra $\mathfrak{Q}$ to connect blocks
$\mathfrak{P}_1$ and $\mathfrak{P}_2$. When it is satisfied, we
have the representation
$$
\mathfrak{Q}_{12}=\vee\biggl\{ \begin{pmatrix} \mathcal P &\\\ &
\mathcal{T}(\mathcal P) \end{pmatrix} \biggm| \mathcal
P\in\mathbb{P}_1 \biggr\}=\vee\{\mathcal P\oplus\mathcal
T(\mathcal P)\mid \mathcal P\in\mathbb P_1\}.
$$
This concludes the proof of the theorem \ref{t1}.
\end{proof}

Let us point out an important circumstance. If the conditions of
the \ref{t1} theorem are satisfied, then the existence of a
connection between blocks $\mathfrak{P}_1$ and $\mathfrak{P}_2$
excludes the connection of any of these blocks with the block
$\mathfrak{P}_3$. This follows from the fact about the separation
of blocks $\mathfrak{Q}_{12}$ and $\mathfrak{Q}_3$, noted in
italics at the beginning of section  4 of the proof.
\medskip

The following result of D.\,V. Korikov \cite{26}, enables one to
clarify the structure of the algebra $\mathfrak{Q}_{12}$ from
theorem \ref{t1}. Let the equality
$\mathcal{T}(\mathbb{P}_1)=\mathbb{P}_2$ be fulfilled and
$n_1=n_2=:n$ accordingly. Let us choose a consistent numbering in
these sets:
$$
\mathbb{P}_1=\{\mathcal{P}_1^i\mid i=1,\dots,{n}\},\qquad
\mathbb{P}_2=\{\mathcal{P}_2^i\mid
\mathcal{P}_2^i=\mathcal{T}(\mathcal{P}_1^i);\,i=1,\dots,n\}.
$$
With every projector $\mathcal{P}_k^i$ we match a one-dimensional subspace $L_k^i$:
$$
L_k^i:= P_k^i \mathscr{G}_k = \operatorname{span}\{\beta_k^i\}\subset \mathscr{G}_k.
$$
An angle between subspaces $L$ and $M$ is defined by the relation
$$
\phi(L,M):=\arccos \|P_{L}P_{M}\|\in\biggl[0,\frac{\pi}2\biggr],
$$
where $P_{L}$, $P_{M}$ are the corresponding orthogonal
projectors. With the families of subspaces
$$
\mathfrak{L}_k:=\{L^1_k,\dots,L^n_k\},\qquad k=1,2,
$$
we associate the sets of angles
\begin{alignat*}{2}
\varphi^i_k &:=\phi(L^i_k,L^1_k+\dots+L^{i-1}_k), &\qquad i &=1,\dots,n,
\\
\varphi^{ij}_k &:=\phi(L^i_k,L^j_k), &\qquad i,j &=1,\dots,n,\, i<j,
\\
\varphi^{ij,l}_k &:=\phi(L^i_k+L^j_k,L^l_k), &\qquad i,j,l &=1,\dots,n,\, i<j<l.
\end{alignat*}
A straightforward application of the unitary equivalence criterion
for families of subspaces from \cite{26} leads to the following
result.

\begin{lemma}
\label{l4} Let the sets $\mathbb{P}_1$ and $\mathbb{P}_2$ satisfy
the conditions of theorem \ref{t1}. Then the map $\mathcal{T}$
extends from the generators $\mathcal{P}_k^i$ to the isomorphism
of the algebras $\mathcal{I}\colon
\mathfrak{P}_1\to\mathfrak{P}_2$, and the algebra
$\mathfrak{Q}_{12}$ takes the form
$$
\mathfrak{Q}_{12} = \{A\oplus {\mathcal{I}} A\mid A\in\mathfrak{P}_1\}
$$
if and only if the equalities
$$
\varphi^i_1=\varphi^i_2,\qquad \varphi^{ij}_1=\varphi^{ij}_2,\qquad \varphi^{ij,l}_1=\varphi^{ij,l}_2
$$
for all $i$, $j$, $l$ take place.
\end{lemma}

It is also easy to show that violating \textit{at least one} of
the equalities in the lemma \ref{l4} causes the algebra
$\mathfrak{Q}_{12}$ to \textit{not connect} blocks
$\mathfrak{P}_1$ and $\mathfrak{P}_2$. Thus, the simultaneous
fulfillment of the conditions of theorem \ref{t1} and lemma
\ref{l4} guarantees that the blocks are connected and that there exists
an isomorphism $\mathcal{I}$, while the lack of at least one of
the conditions leads to the separation of the corresponding blocks by
the algebra $\mathfrak{Q}$.

Consider a more general situation. Let there be $N$ Hilbert spaces
$\mathscr{G}_k$, $k=1,\dots,N$. Each $\mathscr{G}_k$ has its own
set of one-dimensional projectors
$$
\mathbb P_k:=\{P^1_k,\dots,P^{n_k}_k\},\qquad P^i_k=\langle\,{\cdot}\,,\beta^i_k\rangle\,\beta^i_k,\quad \|\beta^i_k\|=1,
$$
defined by the set of vectors $B_k:=\{\beta^i_k\mid
i=1,\dots,{n}_k\}\subset \mathscr{G}_k$, and each $B_k$ is an
equivalence class with respect to $\sim$. The projectors generate
the algebras
$$
\mathfrak{P}_k:=\vee\mathbb P_k=\vee\{P^1_k,\dots,P_k^{n_k}\}
\cong \mathbb{M}^{l_k},
$$
where $l_k:=\operatorname{dim}\operatorname{span}B_k$. Let us
consider the algebra
\begin{equation}
\label{eq4.31}
\mathfrak{P} := \bigoplus_{k=1}^N\mathfrak{P}_k
\end{equation}
with the generators
$$
\mathcal{P}_k^i:= O_1\oplus \dots \oplus {O}_{k-1} \oplus P_k^i
\oplus{O}_{k+1}\oplus\dots \oplus O_N,
$$
where ${O}_k$ is the null operator acting in the $k$th component
of the space $\bigoplus\sum_{k=1}^N\mathscr{G}_k$. We will speak
of the components $\mathfrak{P}_k$ in \eqref{eq4.31}
as \textit{blocks} of the algebra $\mathfrak{P}_*$.

On the complete set of generators
$$
\mathbb{P}:=\mathbb P_1\cup\dots\cup\mathbb
P_N=\{\mathcal{P}_k^i\mid i=1,\dots,{n}_k;\,k=1,\dots, N\}
$$
of the algebra $\mathfrak P$ we define a map
$\mathcal{T}\colon\mathbb{P}\to\mathbb{P}$ such that if
$\mathcal{T}(\mathcal{P}_k^i)=\mathcal{P}_{k'}^{i'}$ then
$\mathcal{T}(\mathcal{P}_{k'}^{i'})=\mathcal{P}_{k}^i$ and one
(and only one) of the conditions is satisfied:
$$
\mathcal{T}(\mathcal{P}_k^i)=\mathcal{P}_k^i\quad\text{or}\quad
\mathcal{T}(\mathcal{P}_k^i)\mathcal{P}_k^i=\mathcal{P}_k^i\mathcal{T}(\mathcal{P}_k^i)=0.
$$
Using this map, we define the projectors $\mathcal{Q}_k^i$:
$$
\mathcal{Q}_k^i := \begin{cases}
\mathcal{P}_k^i, &\text{if }\mathcal{T}(\mathcal{P}_k^i)=\mathcal{P}_k^i,
\\
\mathcal{P}_k^i+\mathcal{T}(\mathcal{P}_k^i), &\text{if } \mathcal{T}(\mathcal{P}_k^i)\mathcal{P}_k^i=0,
\end{cases}
$$
among which there may be identical ones. We construct the algebra
\begin{equation}
\label{eq4.32}
\mathfrak{Q}:=\vee\{\mathcal{Q}_k^i\mid i=1,\dots,{n}_k;\,k=1,\dots, N\}\subset \mathfrak{P}
\end{equation}
and describe its structure using the results of theorem \ref{t1} and lemma \ref{l4}.

\goodbreak 

From the blocks constituting the algebra $\mathfrak{P}$
in \eqref{eq4.31}, let us form all possible pairs of
$\{\mathfrak{P}_k, \mathfrak{P}_{k'}\}$ with $k\ne k'$ and select
those of them, in   which the components are connected through the
algebra $\mathfrak{Q}$ (in  the same sense as $\mathfrak{P}_1$,
$\mathfrak{P}_2$ in theorem \ref{t1} and lemma \ref{l4}). This
selection is unique because, as noted after the proof of the
theorem, each of $\mathfrak{P}_k$ can be connected with only one
$\mathfrak{P}_{k'}$. Let us renumber the blocks in \eqref{eq4.31},
distinguishing pairs of connected blocks and independent blocks:
\begin{equation}
\label{eq4.33}
\underbrace{\mathfrak{P}_1, \mathfrak{P}_{2}};\,\dots;\,\underbrace{\mathfrak{P}_{2k-1},
\mathfrak{P}_{2k}};\,\dots;\,\underbrace{\mathfrak{P}_{2N_1-1},
\mathfrak{P}_{2N_1}};\,\mathfrak{P}_{2N_1+1};\,\dots;\, \mathfrak{P}_{2N_1+j};\,\dots;\,\mathfrak{P}_{N},
\end{equation}
and grouping the components of the set $\mathbb P$ accordingly:
$$
\underbrace{\mathbb{P}_1,
\mathbb{P}_{2}};\,\dots;\,\underbrace{\mathbb{P}_{2k-1},
\mathbb{P}_{2k}};\,\dots; \,\underbrace{\mathbb{P}_{2N_1-1},
\mathbb{P}_{2N_1}};\,\mathbb{P}_{2N_1+1};\, \dots;\,
\mathbb{P}_{2N_1+j};\,\dots;\, \mathbb{P}_{N}.
$$
This numbering will be used hereafter. It can be seen that this
grouping results in the reduction of the map $\mathcal T$ in the
following sense:
\begin{equation}
\begin{gathered}
\mathcal T(\{\mathbb P_{2k-1},\mathbb P_{2k}\})=\{\mathbb
P_{2k-1},\mathbb P_{2k}\},\quad \mathcal T(\mathbb
P_{2k-1})=\mathbb P_{2k},\qquad k=1,\dots,N_1
\\
\mathcal{T}(\mathbb P')=\mathbb P',\quad \text{where}\quad\mathbb P':=\mathbb P_{2N_1+1}\cup\dots\cup\mathbb P_{N}
\end{gathered}
\label{eq4.34}
\end{equation}
The blocks $\mathfrak{P}_{2N_1+j}$ are distinguished by the fact
that they are pairwise separated (not connected) by the
algebra $\mathfrak{Q}$. As it follows from the separation, if the
projector
$$
\mathcal{Q}_{2N_1+j}^i=\mathcal{P}_{2N_1+j}^i+\mathcal{T}(\mathcal{P}_{2N_1+j}^i)
$$
is such that
$\mathcal{T}(\mathcal{P}_{2N_1+j}^i)\notin\mathbb{P}_{2N_1+j}$,
then the algebra $\mathfrak{Q}$ separately includes the projectors
$\mathcal{P}_{2N_1+j}^i$ and
{\vrule width0pt height9pt$\mathcal{T}(\mathcal{P}_{2N_1+j}^i)$}. 
This enables us to replace $\mathcal T|_{\mathbb P'}$ with a new
map $\mathcal T'\colon {\mathbb P'}\to{\mathbb P'}$ that is
defined by $\mathcal{P}_k^i\in\mathbb{P}'$ using the rules
$$
\mathcal T'(\mathcal{P}_k^i):= \begin{cases}
\mathcal{P}_k^i, &\text{if }\mathcal{T}(\mathcal{P}_k^i)\notin \mathbb{P}_k,
\\
\mathcal{P}_k^i, &\text{if }\mathcal{T}(\mathcal{P}_k^i)= \mathcal{P}_k^i,
\\
\mathcal{T}(\mathcal{P}_k^i), &\text{if
}\mathcal{T}(\mathcal{P}_k^i) \in \mathbb{P}_k \text{ and }
\mathcal{T}(\mathcal{P}_k^i) \neq \mathcal{P}_k^i.
\end{cases}
$$
It is easy to see that the map
$$
\widetilde{\mathcal T}\colon \mathbb P\to\mathbb P,\qquad\widetilde{\mathcal T}:=\begin{cases}
\mathcal T &\text{on }\mathbb P\setminus\mathbb P',
\\
\mathcal T' &\text{on }\mathbb P',
\end{cases}
$$
defines the same algebra $\mathfrak{Q}$ as $\mathcal T$, and
in addition to \eqref{eq4.34} is given by independent blocks
$\widetilde{\mathcal T}(\mathbb P_{2N+j})=\mathbb P_{2N+j}$.
Within these blocks, it either acts identically or maps the
projector to orthogonal to it.

The map $\mathcal T$ determines the algebra ${\mathfrak{Q}}$
through the generating projectors $Q^i_k$ according to
\eqref{eq4.32}. Quite similarly, the map $\widetilde{\mathcal T}$
determines corresponding projectors $\widetilde Q^i_k$ generating
the same algebra. The special feature of the latter is in their form:
by construction we have
$$
\widetilde{\mathcal{Q}}^i_k=\begin{cases} {\mathcal
P}_k^i+\widetilde{\mathcal T}({\mathcal P}_k^i), &\text{if }
\widetilde{\mathcal T}({\mathcal P}_k^i)\in\mathbb{P}_k,\ \
\widetilde{\mathcal T}({\mathcal P}_k^i)\neq{\mathcal P}_k^i,
\\
{\mathcal P}_k^i, &\text{in other cases},
\end{cases}
$$
i.e., all two-dimensional $\widetilde{\mathcal{Q}}^i_k$ are sums
of projectors belonging to the same block $\mathfrak{P}_k$.

Let $\mathbb{Q}_k:=\{\widetilde{\mathcal{Q}}_k^i\mid
i=1,\dots,n_k\}$. The following result summarizes the above
considerations. Recall that the blocks are numbered according to
\eqref{eq4.33}.

\begin{proposition}
\label{p8}
Algebra $\mathfrak{Q}$ has the following form
\begin{equation}
\label{eq4.35} \mathfrak{Q}= \biggl[\bigoplus_{k=1}^{N_1}
\mathfrak{Q}_k^{\mathrm{I}}\biggr]\oplus
\biggl[\bigoplus_{k=2N_1+1}^{N}
\mathfrak{Q}_k^{\mathrm{II}}\biggr],
\end{equation}
where $\mathfrak{Q}_k^{\mathrm{I}}:=\{A\oplus \mathcal TA\mid
A\in\mathfrak{P}_{2k-1}\}\subset\mathfrak{P}_{2k-1}\oplus\mathfrak{P}_{2k}$
and $\mathfrak Q^{\mathrm{II}}_k =\vee \mathbb Q_k\subset
\mathfrak P_k$.
\end{proposition}

\subsection{Canonical form}
\label{ss4.5}

Let us use the results established above to describe the structure
of the eikonal algebra. The description will be reduced to some
canonical form (representation).

The grouping of projectors in \eqref{eq4.16}, carried
out according to decomposition \eqref{eq4.15}, is the partitioning
of the sets $\mathbb P_{\Phi^j}$ into equivalence classes $\mathbb
P_{\Phi^j}^k$ (they are also the classes $\mathbb P_{l}$) with
respect to $\sim$. It is motivated by the proposition \ref{p7} and
prepares the application of the theorem \ref{t1}.

The role of the algebra $\mathfrak P$ of the theorem \ref{t1} (see
\ \eqref{eq4.31}) is played by the algebra
$$
\mathfrak P^\partial:=\biggl[ \bigoplus_{l=1}^{L} \mathfrak{P}_l\biggr]\oplus\biggl[ \bigoplus_{l=1}^{L}\mathfrak{P}_l\biggr].
$$
It consists of $2L$ of irreducible blocks and is represented as
$\mathfrak P^\partial=\vee \mathbb P^\partial$, where
\begin{align*}
{\mathbb P}^\partial &:= \{\mathcal{P}_{\gamma l}^k\oplus O\mid k=1,\dots,n_{\gamma l};\, l=1,\dots,L;\, \gamma\in\Sigma\}
\\
&\qquad\cup\{ O\oplus\mathcal{P}_{\gamma l}^k\mid k=1,\dots,n_{\gamma l};\, l=1,\dots,L;\, \gamma\in\Sigma\}
\end{align*}
where $\mathcal{P}_{\gamma l}^k$ are the projectors from
\eqref{eq4.19}, $O$  is the zero element of the algebra
$\bigoplus\limits_{l=1}^{L} \mathfrak{P}_l$. By denoting
$$
\mathcal P_{\gamma l}^{k r}:= \begin{cases}
\mathcal{P}_{\gamma l}^k\oplus O, &r= 0,
\\
O\oplus\mathcal{P}_{\gamma l}^k, &r =\varepsilon_l,
\end{cases}
$$
we have
$$
{\mathbb P}^\partial= \{\mathcal{P}_{\gamma l}^{k r}\mid
k=1,\dots,n_{\gamma l},\, l=1,\dots,L,\, r=0,\varepsilon_l; \,
\gamma\in\Sigma\}
$$
Set the map $\mathcal{T}\colon {\mathbb P}^\partial\to{\mathbb
P}^\partial$ using the formalism introduced earlier, which defines
the relations between the sets $(k,l,r_l)$:
\begin{equation}
\label{eq4.36}
\mathcal{T}(\mathcal P_{\gamma l}^{k r_{l}}):= \begin{cases}
\mathcal{P}_{\gamma l'}^{k' r_{l'}}, &\text{if }(k',l',r_{l'})\leftrightarrow(k,l,r_l), \\
\mathcal{P}_{\gamma l}^{k r_l}, &\text{if } (k,l,r_l)\text{ not related to any } (k',l',r_{l'}).
\end{cases}
\end{equation}
It is not difficult to check that the set $\mathbb P^\partial$ and
the map $\mathcal T$ satisfy all the conditions of theorem
\ref{t1}, and that the algebra $\mathfrak Q$ they define according
to \eqref{eq4.25} coincides with the boundary algebra
\eqref{eq4.20}:
$$
\mathfrak Q=\partial ({U}{\mathfrak{E}}^T_\Sigma U^{-1}).
$$

In bringing the eikonal algebra to canonical form, a
\textit{junction} of its blocks connected through the boundary
algebra will be used. Let us describe this construction.

Let there be two standard algebras $\mathfrak{A}=\dot
C([0,\varepsilon];\mathfrak{P})$ and $\mathfrak{B}=\dot
C([0,\varepsilon'];\mathfrak{P}')$ such that
$$
\mathfrak{A}(\varepsilon):=\{a(\varepsilon)\mid a\in\mathfrak{A}\}=\mathfrak{P},\qquad
\mathfrak{B}(0):=\{b(0)\mid b\in\mathfrak{B}\}=\mathfrak{P}',
$$
and let $\mathfrak{P}\cong\mathfrak{P}'$ holds through isomorphism
$\mathcal I\colon \mathfrak{P}'\to\mathfrak{P}$. Let us define the
algebra
$$
\mathfrak{A}^{\oplus}\mathfrak{B}:= \{a\oplus b\mid
a\in\mathfrak{A},\, b\in\mathfrak{B}, \, a(\varepsilon)=\mathcal I
b(0)\}
$$
Algebras of this kind will appear in course of transforming
to canonical form in the situation, when in the boundary algebra
$\partial (U\mathfrak{E}^T_\Sigma U^{-1})\,{=}\,\mathfrak Q$ 
there are blocks of type $\mathfrak{Q}_{k}^I$ in representation
\eqref{eq4.35} connecting the boundary values of
$\mathfrak{A}(\varepsilon)$ and $\mathfrak{B}(0)$ of some pair of
blocks $\mathfrak{A}$ and $\mathfrak{B}$ of eikonal algebra.

For $a\in\mathfrak{A}$ and $b\in\mathfrak{B}$ such that
$a(\varepsilon)=\mathcal I b(0)$, we define an element $a\sqcup
b\in C([0,\varepsilon+\varepsilon'],\mathfrak{P})$ by the rule
\begin{equation}
\label{eq4.37}
(a\sqcup b)(r):=\begin{cases}
a(r), &r\in[0,\varepsilon],
\\
\mathcal I b(r-\varepsilon), &r\in[\varepsilon,\varepsilon+\varepsilon'],
\end{cases}
\end{equation}
which we call a \textit{junction} of $a$ and $b$. Then we define
the \textit{junction of algebras} $\mathfrak{A}$ and
$\mathfrak{B}$:
$$
\mathfrak{A}\sqcup\mathfrak{B}:= \{a\sqcup b\in
C([0,\varepsilon+\varepsilon'];\mathfrak{P})\mid
a\in\mathfrak{A},\, b\in\mathfrak{B},\, a(\varepsilon)=\mathcal I
b(0) \}.
$$
We see that $\mathfrak{A}\sqcup\mathfrak{B}$ is a subalgebra
in $C([0,\varepsilon+\varepsilon'];\mathfrak{P})$, which is a
standard algebra, for which the representation
$$
\mathfrak{A}\sqcup\mathfrak{B}=\{c\in C([0,\varepsilon+\varepsilon'];\mathfrak{P})\mid
c(0)\in\mathfrak{A}(0),\, c(\varepsilon+\varepsilon')\in\mathcal I
[\mathfrak{B}(\varepsilon')]\}
$$
and the equalities
$$
(\mathfrak{A}\sqcup\mathfrak{B})(0) = \mathfrak{A}(0),\qquad (\mathfrak{A}\sqcup\mathfrak{B})(\varepsilon+\varepsilon') = \mathfrak{B}(\varepsilon')
$$
hold.

Note that the algebras $\mathfrak{A}^{\oplus}\mathfrak{B}$ and
$\mathfrak{A}\sqcup\mathfrak{B}$ are isomorphic. To summarize the
considerations, let us say that the algebras $\mathfrak{A}$ and
$\mathfrak{B}$ admit a junction \textit{through the ends}
$r=\varepsilon$ and $r'=0$. Note that the algebras
$\mathfrak{A}(0)$ and $\mathfrak{B}(\varepsilon')$ have no
influence on the possibility of the junction and its result.

Obviously, by changing the definition of \eqref{eq4.37}, we can
introduce the junction $\mathfrak{A}\sqcup\mathfrak{B}\subset
C([0,\varepsilon+\varepsilon'];\mathfrak{P})$
\begin{equation}
\begin{aligned}
&\text{through the ends of }r=\varepsilon \text{ and } r'=\varepsilon', \text{ if } \mathfrak{A}(\varepsilon)=\mathcal I[\mathfrak{B}(\varepsilon')]
\\
&\text{through the ends of }r=0 \text{ and } r'=\varepsilon', \text{ if } \mathfrak{A}(0)=\mathcal I[\mathfrak{B}(\varepsilon')];
\\
&\text{through the ends } r=0 \text{ and } r'=0, \text{ if } \mathfrak{A}(0)=\mathcal I[\mathfrak{B}(0)].
\end{aligned}
\label{eq4.38}
\end{equation}

On the algebra elements $\mathfrak{A}=\dot
C([0,\varepsilon],\mathfrak{P})$ we define \textit{transposition}
$t\colon a\mapsto a^t$, $a^t(r): =a(\varepsilon-r)$,
$r\in[0,\varepsilon]$, and take $\mathfrak{A}^t:=\{a^t\mid
a\in\mathfrak{A}\}$. The isomorphism ${\mathcal M}\colon
\mathfrak{P}\to\mathfrak{P}$ defines the transformation
$\check{\mathcal M}\colon C([0, \varepsilon],\mathfrak{P})\to
C([0,\varepsilon],\mathfrak{P})$, $(\check{\mathcal M}a)(r):
=\mathcal M [a(r)]$, $r\in[0,\varepsilon]$. For the standard
algebra $\mathfrak{A}$ we take $\check{\mathcal
M}\mathfrak{A}:=\{\check{\mathcal M} a\mid
a\in\mathfrak{A}\}\subset C([0,\varepsilon],\mathfrak{P})$.

The algebras $\mathfrak{A}^t$ and $\mathcal M\mathfrak{A}$ are
also standard; it is easy to see that
$\check{\mathfrak{A}}\cong{\mathfrak{A}}^t\cong{\mathfrak{A}}$ and

$$
\mathfrak{A}^t(0)=\mathfrak{A}(\varepsilon),\quad
\mathfrak{A}^t(\varepsilon)=\mathfrak{A}(0); \qquad
(\check{\mathcal M}\mathfrak{A})(0)= \mathcal
M[\mathfrak{A}(0)],\quad (\check{\mathcal
M}\mathfrak{A})(\varepsilon)=\mathcal M[\mathfrak{A}(\varepsilon)]
$$
holds.

The transformation to the canonical form of the eikonal algebra is
quite prepared and can be described by the following procedure.

\textsl{Step}~1.~Define the isomorphism $\mathbf{U}_0\colon
\mathfrak{E}^T_{\Sigma}\to U\mathfrak{E}^T_\Sigma
U^{-1}\stackrel{\eqref{eq4.18}}{\subset} \bigoplus\limits_{l=1}^L
C([0,\varepsilon_l],\mathfrak{P}_l)$ by defining it on the
generators:
$$
\mathbf{U}_0 E^T_{\gamma} := U E^T_{\gamma}U^{-1}\stackrel{\eqref{eq4.19}}{=} \bigoplus\sum_{l=1}^L \biggl[\sum_{k=1}^{n_{\gamma l}}\tau^k_{\gamma l} P^k_{\gamma l}\biggr].
$$
The algebra blocks $\mathbf{U}_0 \mathfrak{E}^T_{\Sigma}$ are
$$
[\mathbf{U}_0 \mathfrak{E}^T_{\Sigma}]_l:=\vee\biggl\{\sum_{k=1}^{n_{\gamma l}}\tau^k_{\gamma l}\biggm| \gamma\in\Sigma\biggr\}\subset C([0,\varepsilon_l];\mathfrak{P}_l).
$$
Let the blocks $[\mathbf{U}_0 \mathfrak{E}^T_{\Sigma}]_l$ and
$[\mathbf{U}_0 \mathfrak{E}^T_{\Sigma}]_{l'}$ be such, that their
boundary values $[\mathbf{U}_0
\mathfrak{E}^T_{\Sigma}]_l(\varepsilon_l)$ and $[\mathbf{U}_0
\mathfrak{E}^T_{\Sigma}]_{l'}(0)$ form one block of type
${\mathfrak{Q}}_k^I$ in~boundary algebra
{\vrule width0pt height10pt$\partial(\mathbf{U}_0 \mathfrak{E}^T_{\Sigma}):=\partial (U\mathfrak{E}^T_\Sigma U^{-1}) = \mathfrak Q$} 
(see proposition \ref{p8}). In~this case, the map $\mathcal T$
(see ~\eqref{eq4.36}) determines an isomorphism $\mathcal I\colon
\mathfrak{P}_{l'}\to\mathfrak{P}_l$:
$$
\mathcal I (P^{k'}_{\gamma l'}) = P^{k}_{\gamma l },\quad \text{if}\quad \mathcal T(\mathcal P^{k' 0}_{\gamma l'})=\mathcal P^{k \varepsilon_l}_{\gamma l }
$$
As a consequence, the algebra $[\mathbf{U}_0
\mathfrak{E}^T_{\Sigma}]_l\sqcup[\mathbf{U}_0
\mathfrak{E}^T_{\Sigma}]_{l'}$ that is the junction of these
blocks through the ends $r_l=\varepsilon_l$ and $r_{l'}=0$, is
well defined. It is not difficult to check that the elements
$$
\bigl( [\mathbf{U}_0 E^T_{\gamma}]_l\sqcup [\mathbf{U}_0 E^T_{\gamma}]_{l'}\bigr)(r) :=\begin{cases}
{\displaystyle\sum_{k=1}^{n_{\gamma l}}\tau^k_{\gamma l}(r) P^k_{\gamma l}}, &r\in[0,\varepsilon_l],
\\[4mm]
{\displaystyle\mathcal I \biggl[\sum_{k'=1}^{n_{\gamma l'}}\tau^{k'}_{\gamma l'} (r-\varepsilon_l)P^{k'}_{\gamma l'}\biggr]}, &r\in[\varepsilon_l,\varepsilon_l+\varepsilon_{l'}],
\end{cases}
$$
constitute the system of its generators.


By the properties of the map $\mathcal T$, the existence of a
connection between the blocks implies $n_{\gamma l}=n_{\gamma l'}$
for any $\gamma\in\Sigma$, and the junctions $[\mathbf{U}_0
E^T_{\gamma}]_l\sqcup [\mathbf{U}_0 E^T_{\gamma}]_{l'}$ can be
represented in~ the following way:
$$
\bigl( [\mathbf{U}_0 E^T_{\gamma}]_l\sqcup [\mathbf{U}_0
E^T_{\gamma}]_{l'}\bigr)(r) =\sum_{k=1}^{n_{\gamma
l}}[\tau^k_{\gamma l }\sqcup\tau^{k'}_{\gamma l}](r)P^k_{\gamma l}
\qquad r\in[0,\varepsilon_l+\varepsilon_{l'}],
$$
where
$$
[\tau^k_{\gamma l}\sqcup\tau^{k'}_{\gamma l' }](r):= \begin{cases}
\tau^k_{\gamma l}(r), &r\in[0,\varepsilon_l],
\\
\tau^{k'}_{\gamma l'}(r-\varepsilon_l), &r\in[\varepsilon_l,\varepsilon_l+\varepsilon_{l'}].
\end{cases}
$$
The functions $[\tau^k_{\gamma l}\sqcup\tau^{k'}_{\gamma l'}]$ are
continuous due to the equality $\tau^k_{\gamma
l}(\varepsilon_{l})=\tau^{k'}_{\gamma l'}(0)$ determining the
relation between the sets $(k,l,\varepsilon_{l})$ and $(k',l',0)$,
which in turn determines the map $\mathcal T$ and isomorphism
$\mathcal I$. Recall that each $\tau^k_{\gamma l}$ is a linear
function of one of two kinds:
$$
\text{or }\tau^k_{\gamma l}(r)= t^k_{\gamma l}+r, \text{or } \tau^k_{\gamma l}(r) =t^k_{\gamma l }-r,
$$
where $t^k_{\gamma l }=\mathrm{const}\ge 0$.

Let $\tau^k_{\gamma l}(r)= t^k_{\gamma l}+r$ for certainty. Then
the condition $\tau^k_{\gamma
l}(\varepsilon_{l})=\tau^{k'}_{\gamma l'}(0)$ and the fact that
the equality $\tau^k_{\gamma l}(r_{l})=\tau^{k'}_{\gamma
l'}(r_{l'})$ is possible only in~the case of boundary values of
parameters $r_l$ and $r_{l'}$ (in our case $r_l=\varepsilon_l$,
$r_{l'}=0$), follow to $\tau^{k'}_{\gamma l'}(r)= t^{k'}_{\gamma
l'}+r$, with $t^{k'}_{\gamma l'}= t^{k}_{\gamma l}+\varepsilon_l$.
This leads to the following equalities:
\begin{align*}
[\tau^k_{\gamma l}\sqcup\tau^{k'}_{\gamma l'}](r) &= \begin{cases}
t^k_{\gamma l }+r, &r\in[0,\varepsilon_l],
\\
(t^{k}_{\gamma l}+\varepsilon_{l})+(r-\varepsilon_l), &r\in[\varepsilon_l,\varepsilon_l+\varepsilon_{l'}],
\end{cases}
\\
&=\begin{cases}
t^k_{\gamma l}+r, &r\in[0,\varepsilon_l],
\\
t^{k}_{\gamma l}+r, &r\in[\varepsilon_l,\varepsilon_l+\varepsilon_{l'}],
\end{cases} = t^k_{\gamma l }+r,\qquad r\in[0,\varepsilon_l+\varepsilon_{l'}].
\end{align*}
Analogous reasoning is true in the~case $\tau^k_{\gamma
l}(r)=t^k_{\gamma l}-r$. Thus, the junction of the functions
$\tau^k_{\gamma l}$ and $\tau^{k'}_{\gamma l'}$ is a linear
function of the same kind as $\tau^k_{\gamma l}$ and
$\tau^{k'}_{\gamma l'}$ themselves. It implies that the junction
of the blocks $[\mathbf{U}_0
\mathfrak{E}^T_{\Sigma}]_l\sqcup[\mathbf{U}_0
\mathfrak{E}^T_{\Sigma}]_{l'}$ is a standard algebra with
~generators of the same kind as the original blocks.

In~the case considered, the boundary values $[\mathbf{U}_0
\mathfrak{E}^T_{\Sigma}]_l(\varepsilon_l)$ and $[\mathbf{U}_0
\mathfrak{E}^T_{\Sigma}]_{l'}(0)$ are connected. The cases of
other possible connections between the boundary values, admitting
the junction of blocks (see \ \eqref{eq4.38}), are treated quite
analogously.

\textsl{Step}~2.~Let the blocks $[\mathbf{U}_0
\mathfrak{E}^T_{\Sigma}]_l$ and $[\mathbf{U}_0
\mathfrak{E}^T_{\Sigma}]_{l'}$ be admissible for junction. Let us
define the map
$$
\mathbf{U}_1\colon\mathfrak{E}^T_{\Sigma}\to\mathbf{U}_1
\mathfrak{E}^T_{\Sigma}\subset\biggl[\bigoplus_{\lambda=1\
(\lambda\neq l,l')}^{L} C([0,\varepsilon_\lambda];
\mathfrak{P}_\lambda)\biggr]\oplus
C([0,\varepsilon_{l}+\varepsilon_{l'}];\mathfrak{P}_{l}),
$$
by defining it on the generators
$$
\mathbf{U}_1 E^T_{\gamma} := \biggl[\bigoplus\sum_{\lambda=1\
(\lambda\neq l,l')}^{L} [\mathbf{U}_0 E^T_{\gamma}]_\lambda\biggr]
\oplus \bigl([\mathbf{U}_0 E^T_{\gamma}]_{l}\sqcup [\mathbf{U}_0
E^T_{\gamma}]_{l'}\bigr),\qquad \gamma\in\Sigma .
$$
It is not difficult to see that the map $\mathbf{U}_1$ is
isomorphic: the isomorphism of the algebras
$\mathbf{U}_0\mathfrak{E}^T_{\Sigma}$ and
$\mathbf{U}_1\mathfrak{E}^T_{\Sigma}$ is quite similar to
isomorphism $\mathfrak{A}^{\oplus}\mathfrak{B}$ and
$\mathfrak{A}\sqcup\mathfrak{B}$. The isomorphism $\mathbf{U}_1$
transfers the eikonal algebra into an algebra of the same kind as
$\mathbf{U}_0\mathfrak{E}^T_{\Sigma}$, but with ~$1$ fewer blocks;
in~the boundary algebra
$\partial(\mathbf{U}_1\mathfrak{E}^T_{\Sigma})$ the number of
blocks of type~${\mathfrak{Q}}_{k}^I$ also becomes by ~$1$ less.

Consistently replacing all pairs of connected blocks by their
junctions in a similar way, we arrive at some isomorphism
$\mathbf{U}_N\colon
\mathfrak{E}^T_{\Sigma}\to\mathbf{U}_N\mathfrak{E}^T_{\Sigma}$
such that that the~boundary algebra
$\partial(\mathbf{U}_N\mathfrak{E}^T_{\Sigma})$ no longer contains
any block of type $\mathfrak{Q}_{k}^I$. Thus, the map
$\mathbf{U}_N$ transforms the eikonal algebra into an~algebra of
the same structure as $\mathbf{U}_0\mathfrak{E}^T_{\Sigma}$, but
with~independent blocks. Note also that in the process of
transformation to~the representation
$\mathbf{U}_N\mathfrak{E}^T_{\Sigma}$ several (more than two)
blocks of the original representation
$\mathbf{U}_0\mathfrak{E}^T_{\Sigma}$ can be joined into~one new
block.

\textsl{Step}~3.~
The final step once again uses the irreducibility of the
algebras $\mathfrak{P}_l$. For each of them, let us find a
transformation realizing the isomorphism $\mathfrak{P}_l\,{\cong}\,\mathbb M^{\varkappa_l}$. 
As a consequence, an isomorphism $\mathbf{U}\colon \mathfrak{E}^T_{\Sigma}\,{\to}\, \bigoplus_{l=1}^\mathcal L \dot C([0,\zeta_l];\mathbb{M}^{\varkappa_l})$, 
defined on the generators by the equalities
$$
\mathbf{U} E^T_{\gamma}:=\bigoplus\sum_{l=1}^\mathcal L\biggl[\sum_{k=1}^{s_{\gamma l }}\widetilde\tau_{\gamma l}^k(\cdot_l) \widetilde P_{\gamma l}^k\biggr],\qquad \gamma\in\Sigma,
$$
where $\widetilde P_{\gamma l}^k\in\mathbb{M}^{\varkappa_l}$ are
one-dimensional (matrix) projectors, pairwise orthogonal for a
fixed vertex $\gamma\in\Sigma$, and $\widetilde\tau_{\gamma
l}^k$~are the linear functions of one of two kinds: either
$\widetilde\tau^k_{\gamma l}(r)=\widetilde t^{\,k}_{\gamma l}+r$,
or $\widetilde\tau^k_{\gamma l}(r)=\widetilde t^{\,k}_{\gamma
l}-r$, $r\in[0,\zeta_l]$, where $\widetilde t^{\,k}_{\gamma l
}\ge0$ are constants, and each $\zeta_l$ is the sum of some set of
the lengths $\varepsilon_k$.

Thus, by successive transformations of the original parametric
form \eqref{eq4.6}, we obtain a form of the same structure, but
with independent blocks which are standard algebras. Let us
formulate the final result.

\begin{theorem}\label{t2}
There is an isomorphism $\mathbf U$, which transforms the algebra
$\mathfrak{E}_{\Sigma}^T$ and its generating elements-eikonals to
the representation
\begin{equation}
\label{eq4.39} \mathbf U \mathfrak{E}_{\Sigma}^T
=\bigoplus_{l=1}^\mathcal L \dot C([0,\zeta_l];
\mathbb{M}^{\varkappa_l}), \qquad \mathbf{U} E_{\gamma}^T =
\bigoplus\sum_{l=1}^\mathcal L\biggl[\sum_{k=1}^{s_{\gamma
l}}\widetilde\tau_{\gamma l}^k \widetilde P_{\gamma l}^k\biggr],
\quad \gamma\in\Sigma.
\end{equation}
Here, $\widetilde\tau_{\gamma l}^k$~are linear functions of
$r_l\in[0,\zeta_l]$ such that $|d\widetilde \tau_{\gamma
l}^k/dr_l|=1$, and their ranges are the segments of length
$\zeta_l$, which can only have common ends. $\widetilde P_{\gamma
l}^k\in \mathbb{M}^{\varkappa_l}$ are projectors, pairwise
orthogonal for each~$\gamma$ and such that $\vee\{\widetilde
P_{\gamma l}^k\mid k=1,\dots ,s_{\gamma l} \,
\gamma\in\Sigma\}=\mathbb{M}^{\varkappa_l}$.
\end{theorem}

We call a representation (form) of this kind \textit{canonical}.
It is not unique, but it can be shown that any two such
representations differ from each other by block numbering
$[\mathbf{U}\mathfrak{E}^T_{\Sigma}]_l$, by their transposition
$[\mathbf{U}\mathfrak{E}^T_{\Sigma}]_l\to[\mathbf{U}\mathfrak{E}^T_{\Sigma}]_l^t$
and isomorphisms
$[\mathbf{U}\mathfrak{E}^T_{\Sigma}]_l\to\check{\mathcal
M}[\mathbf{U}\mathfrak{E}^T_{\Sigma}]_l$. The ambiguity associated
with~transposition obviously corresponds to two directions of
change in the argument of each function $\widetilde\tau_{\gamma
l}^k$ on the interval $[0,\delta_l]$ (two possible
parameterizations of the $l$th block).

As a comment on the \ref{t2} theorem, we note the following. The
original parametric form of the eikonal algebra \eqref{eq4.8}
corresponded to a partition of the graph $\Omega^T_\Sigma$ into
families $\Phi^1,\dots,\Phi^J$. Our conjecture is that the
transition to the canonical form corresponds to some new
partitioning. This is not proven, but it is supported by
well-known examples \cite{7}, \cite{8}.

\section{Transformation to canonical form}
\label{s5}

\subsection{Canonical representation and spectra}
\label{ss5.1}

The main advantage of the ca\-no\-ni\-cal form is that complete
information about the spectrum of the algebra
$\mathfrak{E}^T_\Sigma$ is easily extracted from it and a number
of its invariants are revealed. Let us proceed to~their
description.

Let us compare the content of the theorem \ref{t2} and the proposition \ref{p2}.

Given the ambiguity of the canonical form noted after the Theorem
\ref{t2}, we assume that the numbering and parametrization of its
blocks are fixed. Let us rewrite the form itself in~new convenient
notations (replacing $\zeta_l$ by $\varepsilon_l$ and $s_{{\gamma
l}}$ by $n_{\gamma l}$ , and removing $(\,\widetilde{\ }\,)$):
\begin{equation}
\label{eq5.1} \mathbf U \mathfrak{E}_{\Sigma}^T
=\bigoplus_{l=1}^\mathcal L \dot C([0,\varepsilon_l];
\mathbb{M}^{\varkappa_l}),\qquad \mathbf{U} E_{\gamma}^T =
\bigoplus\sum_{l=1}^\mathcal L\biggl[\sum_{k=1}^{n_{\gamma
l}}\tau_{\gamma l}^k P_{\gamma l}^k\biggr], \quad \gamma\in\Sigma.
\end{equation}

Let
$$
\psi_{\gamma l}^k :=\{\tau_{\gamma l}^k(r_l)\mid
r_l\in(0,\varepsilon_{l})\}=(\tau^k_{\gamma l}(0),\tau_{\gamma
l}^k(\varepsilon_l))
$$
be the time cells corresponding to the canonical representation.
The right-hand side of the representation \eqref{eq5.1} for the
operator $\mathbf{U} E_{\gamma}^T$ defines its spectrum and, since
$\mathbf{U}$ is an isomorphism, the eikonal spectrum:
\begin{equation}
\label{eq5.2} \sigma_{\mathrm{ac}}({E}^T_{\gamma})
=[1,T_{1}^{\gamma}]\cup[T_2^{\gamma},T_3^{\gamma}]
\cup\dots\cup[T_{N_{\gamma}-1}^{\gamma}, T_{N_{\gamma}}^{\gamma}]
=\bigcup_{l=1}^{\mathcal{L}}\bigcup_{k=1}^{n_{\gamma
l}}\overline{\psi_{\gamma l}^k},
\end{equation}
where each of the segments $[T^\gamma_{i-1},T^\gamma_i]$, in ~ its
turn, is covered by cells $\overline{\psi_{\gamma l}^k}$ that
either do not overlap or have common ends. By the proposition
\ref{p2}, the same is true for the cells
$\overline{\psi^i_{\gamma\Phi}}$ associated with~parametric
representation~\eqref{eq3.12}. Comparing the decompositions of
\eqref{eq5.2} and \eqref{eq3.13}, we conclude that the canonical
representation corresponds to a new (canonical) cutting of the
eikonal spectrum into time cells. It can be shown that each of the
new cells consists of the old ones, i.e., the transition to the
canonical representation corresponds to the enlargement of the
time cells.

As noted in ~p. ~\ref{ss4.1}, the spectrum of the standard algebra
$\dot C([a,b]; \mathbb M^n)$ consists of (equivalence classes)
irreducible representations corresponding to interior points of
the segment $[a,b]$ and clusters (if any)
$\{\widehat\pi^1_a,\dots,\widehat\pi^{n_a}_a\}$ and
$\{\widehat\pi^1_b,\dots,\widehat\pi^{n_b}_b\}$ adjacent to~its
ends. Such a spectrum equipped with the Jacobson topology is
homeomorphic to a space, which is naturally called a
\textit{segment with~split ends}. It is described by the following
construction (see, e.g., \cite{27}).

Consider $n_a$ semisegments $[a,b)$ and $n_b$ semisegments $(a,b]$
with~topology from $\mathbb R$. Identify the interior points of
all semisegments having the same coordinates. The resulting
factor-space $\mathscr{S}_{[a,b]}$ consists of the part
$S_{(a,b)}$ homeomorphic to $(a,b)$ and two sets of pairwise
inseparable points (two clusters) $K_a$ and $K_b$, which consist
of $n_a$ and $n_b$ points and correspond to ends of $a$ and $b$
respectively. Each of the clusters is inseparable from
$S_{(a,b)}$. The part
$S_{(a,b)}:=\operatorname{int}\mathscr{S}_{[a,b]}$ is the set of
interior points, each of which has a neighborhood homeomorphic to
an (open) interval of the real axis.

According to the first of the representations \eqref{eq5.1}, the
spectrum of the algebra $\mathbf U \mathfrak{E}_{\Sigma}^T$ is
homeomorphic to the disjunct union of segments
$$
\mathscr{S}^T_\Sigma=\mathscr{S}_{[0,\varepsilon_1]}\cup\dots
\cup\mathscr{S}_{[0,\varepsilon_\mathcal L]},
\qquad\mathscr{S}_{[0,\varepsilon_l]}: =K^l_{0}\cup
S_{(0,\varepsilon_l)}\cup K^l_{\varepsilon_l}.
$$
Each segment $\mathscr{S}_{[0,\varepsilon_l]}$ is characterized as
a \textit{linearly--connected component} of $\mathscr{S}^T_\Sigma$
space, and the part
$S_{(0,\varepsilon_l)}=\operatorname{int}\mathscr{S}_{[0,\varepsilon_l]}$~--
as the set of its interior points.

The spectra of the algebras $\mathfrak{E}_{\Sigma}^T$ and $\mathbf
U \mathfrak{E}_{\Sigma}^T$ are connected by the homeomorphism
$\mathbf U_*\colon
\widehat{\mathfrak{E}_{\Sigma}^T}\to\widehat{\mathfrak{E}_{\Sigma}^T}$
(see ~\eqref{eq4.2}). Hence, $\widehat{\mathfrak{E}_{\Sigma}^T}$
is homeomorphic to the space $\mathscr{S}^T_\Sigma$ and admits the
representation
\begin{equation}
\label{eq5.3}
\widehat{\mathfrak{E}_{\Sigma}^T}=\mathscr{S}_1\cup\dots
\cup\mathscr{S}_\mathcal L, \qquad
\mathscr{S}_l=\eta(\mathscr{S}_{[0,\varepsilon_l]})=\mathscr{K}^l_0\cup
S_l\cup\mathscr{K}^l_{\varepsilon_l},
\end{equation}
where
$\eta\colon\mathscr{S}^T_\Sigma\to\widehat{\mathfrak{E}_{\Sigma}^T}$~is
a homeomorphism, $S_l=\eta(S_{(0, \varepsilon_l)})$,
$\mathscr{K}^l_0=\eta(K^l_0)$,
$\mathscr{K}^l_{\varepsilon_l}=\eta(K^l_{\varepsilon_l})$. It
determines the partitioning of the spectrum into
linearly--connected components and, as such, has an invariant
topological meaning. The latter also applies to the ~partitioning
of segments $\mathscr{S}_l$ into interior points $S_l$ and
clusters $\mathscr{K}^l_0$, $\mathscr{K}^l_{\varepsilon_l}$. The
set of interior points of the spectrum is
$\operatorname{int}\widehat{\mathfrak{E}_{\Sigma}^T}=S_1\cup\dots
\cup S_\mathcal L$.

\subsection{Coordinates}
\label{ss5.2}

Recall that the algebras
$\mathfrak{E}_{\gamma}^T=\vee\{{{E}^T_{\gamma}}\}\subset\mathfrak{E}_{\Sigma}^T$
corresponding to individual eikonals, are called partial. We
denote by $\pi$ and $\widehat\pi$ an irreducible representation
and its equivalence class.

The correspondence $\varphi({E}^T_{\gamma})\leftrightarrow\varphi$
is an isomorphism of the algebras $\mathfrak{E}_{\gamma}^T$ and
$C(\sigma_{\mathrm{ac}}( E^T_\gamma))$, defined by the first
equality in~\eqref{eq4.5}\enskip (see \cite{23}). Each
$\mathfrak{E}_{\gamma}^T$ is a commutative subalgebra
in~$\mathfrak{E}_{\Sigma}^T$. Its spectrum (set of characters)
$\widehat{\mathfrak{E}_{\gamma}^T}$ is exhausted by Dirac
measures:
$$
\widehat{\mathfrak{E}_{\gamma}^T}=\{\widehat{\delta}_t \mid t\in
\sigma_{\mathrm{ac}}(E^T_{\gamma})\},\qquad
\mathfrak{E}_{\gamma}^T\ni
\varphi(E^T_\gamma)\stackrel{\delta_t}{\mapsto}
\varphi(t)\in\mathbb R,
$$

\vskip-1mm 

\noindent in this case $\widehat\delta_t=\{\delta_t\}$\enskip (see
\cite{25}, \cite{28}) holds. Thus, each character corresponds to a
point (number) $t$ from the union of segments \eqref{eq5.2}, which
we will regard as its $\gamma$-coordinate. Note also that any
reducible matrix representation of the algebra
$\mathfrak{E}_{\gamma}^T$ is of the form
\begin{equation}
\label{eq5.4} \rho
\sim\delta_{t_1}\oplus\dots\oplus\delta_{t_p},\qquad
\rho(\varphi(E^T_\gamma))\sim\operatorname{diag}\{\varphi(t_1),\dots,\varphi(t_p)\}
\end{equation}

\vskip-1mm 

\noindent The numbers $t_1,\dots,t_p$ that uniquely define the
class $\widehat\rho$, are said to be its $\gamma$--coor\-di\-na\-tes.
The next step provides a coordinatization of the spectrum of the
algebra $\mathfrak{E}_{\Sigma}^T$.

If a commutative algebra has a finite number of generators, they
constitute the coordinates on its spectrum (see \ \cite[Ch.~III,
Theorem~6]{28}). Here we use an adequate analog of this fact for
the noncommutative $C^*$--algebra whose generators are
self-adjoint operators with simple spectrum. Namely, to each
$\widehat\pi\in\widehat{\mathfrak{E}_{\Sigma}^T}$ we correspond
the set $\{\widehat\pi|_{\mathfrak{E}_{\gamma}^T}\mid
\gamma\in\Sigma\}$ of contractions to partial algebras. Each
element of the set is a representation of the form \eqref{eq5.4},
already provided with $\gamma$--coordinates
$t^1_{\gamma},t^2_{\gamma},\dots$\,. The correspondence
\begin{equation}
\label{eq5.5}
\widehat{\mathfrak{E}_{\Sigma}^T}\ni\widehat\pi \to \bigl\{\{t^1_{\gamma},t^2_{\gamma},\dots\}\bigm| \gamma\in\Sigma\bigr\},\qquad t^k_{\gamma}=t^k_{\gamma}(\widehat\pi),
\end{equation}

\vskip-1mm 

\noindent provides the required spectrum coordinates. Let's
clarify some details.

\goodbreak 

Correspondence \eqref{eq5.5}, in general, is not injective, but it
is when contracted onto
$\operatorname{int}\widehat{\mathfrak{E}_{\Sigma}^T}$. Since the
set
$\widehat{\mathfrak{E}_{\Sigma}^T}\setminus\operatorname{int}\widehat{\mathfrak{E}_{\Sigma}^T}$
consisting of points forming clusters is finite (so, almost all
points of the spectrum are interior), the use of the term
coordinates seem motivated.

If $\widehat\pi\in S_l\subset\mathscr{S}_l$ (see \ \eqref{eq5.3}),
then $\mathbf U_*\widehat\pi$ is an interior point of the spectrum
$\widehat{\mathbf U_*\mathfrak{E}_{\Sigma}^T}$ and, by the first
of the representations in \eqref{eq5.1},  it corresponds to a
certain value of the parameter $r\in (0,\varepsilon_l)$. In that
case, we denote the point $\widehat\pi$ by $\widehat\pi_r$ and
find out what its $\gamma$-coordinates are.

Let $\varphi\in C(\sigma_{\mathrm{ac}}( E^T_\gamma))$. From the
equivalent representations constituting $\widehat\pi_r$, we choose
$\pi_r$ by the condition
$$
({\mathbf U_*}\pi_r)({\mathbf U}E^T_\gamma)=\sum_{k=1}^{n_{\gamma l}}\tau_{\gamma l}^k(r) P_{\gamma l}^k
$$
(see \eqref{eq5.1}). For such $\pi_r$ we have the relations
\begin{align}
\notag &\pi_r\bigl(\varphi( E^T_\gamma)\bigr)
\stackrel{\eqref{eq4.1}}{=}(\mathbf U_*\pi_r)\bigl(\mathbf U
\varphi( E^T_\gamma)\bigr)=(\mathbf U_*\pi_r)\bigl(\varphi(\mathbf
UE^T_\gamma)\bigr)
\\
\notag &\stackrel{\eqref{eq5.1}}{=}(\mathbf
U_*\pi_r)\biggl(\varphi\biggl(\bigoplus\sum_{l=1}^\mathcal
L\sum_{k=1}^{n_{\gamma l}}\tau_{\gamma l}^k P_{\gamma
l}^k\biggr)\biggr)
\\
\notag &\,\,= (\mathbf
U_*\pi_r)\biggl(\bigoplus\sum_{l=1}^\mathcal
L\sum_{k=1}^{n_{\gamma l}}(\varphi\circ \tau_{\gamma l}^k)
P_{\gamma l}^k\biggr)
\\
&\,\, =\sum_{k=1}^{n_{\gamma l}}(\varphi\circ\tau_{\gamma l}^k)(r)
P_{\gamma l}^k= \sum_{k=1}^{n_{\gamma l}}\varphi(t_{\gamma l}^k)
P_{\gamma l}^k =\sum_{k=1}^{n_{\gamma l}}\delta_{t^k_{\gamma
l}}(\varphi) P_{\gamma l}^k, \label{eq5.6}
\end{align}
where
\begin{equation}
\label{eq5.7} t_{\gamma l}^k = t_{\gamma l}^k(\widehat\pi_r)
:=\tau_{\gamma l}^k(r), \qquad r\in(0,\varepsilon_l),\quad
k=1,\dots,n_{\gamma l}
\end{equation}
Comparing the beginning and end of this calculation, we conclude
that the representation~$\pi_r|_{\mathfrak{E}^T_\gamma}$ is
equivalent to the given representation
$\bigoplus\sum_{k=1}^{n_{\gamma l}}\delta_{t^k_{\gamma l}}$ of
algebra~$\mathfrak{E}^T_\gamma$ and the numbers $\tau_{\gamma
l}^k(r)$ are $\gamma$-coordinates $\widehat\pi_r\in
S_l\subset\operatorname{int}\widehat{\mathfrak{E}^T_\Sigma}$. They
are well defined because choosing another $\pi\ne \pi_r$,
$\pi\in\widehat\pi_r$ will only result in~ replacing the
projectors $P^k_{\gamma l}$ with unitary equivalent projectors
in~\eqref{eq5.6}.

When $\varphi(t)=t$, the relation \eqref{eq5.6} takes the form
\begin{equation}
    \label{eq5.8}
    \pi_r( E^T_\gamma)=\sum_{k=1}^{n_{\gamma l}}\tau^k_{\gamma l}(r)\,P_{\gamma l}^k,\qquad r\in (0,\varepsilon_l),\quad\gamma\in\Sigma,
\end{equation}
which we use below.

As can be seen from \eqref{eq5.7}, with variations of the point
$\widehat\pi$ in~$S_l$ its coordinates $t_{\gamma
l}^k(\widehat\pi)$ covers the time cells $\psi_{\gamma
l}^k\subset\sigma_{\mathrm{ac}}( E^T_\gamma)$, and the following
is fulfilled
\begin{equation}
    \label{eq5.9}
    \psi_{\gamma l}^k=\{t_{\gamma l}^k(\widehat\pi)\mid \widehat\pi\in S_l\},\qquad k=1,\dots,n_{\gamma l}, \quad l=1,\dots,\mathcal L,\quad \gamma\in\Sigma,
\end{equation}
which indicates the invariant nature of these cells, and hence the
partition~\eqref{eq5.2}. We mean that they are uniquely determined
by the eikonal algebra \textit{itself}; more precisely --- by the
structure of \eqref{eq5.3} of its spectrum. The cell lengths
$\varepsilon_l$, $l=1,\dots,\mathcal L$, are numerical invariants
of the algebra ${\mathfrak{E}^T_\Sigma}$.

When $k$ is fixed, the correspondence
$\operatorname{int}\psi^k_{\gamma l}\ni t^k_{\gamma
l}(\widehat\pi)\leftarrow\widehat\pi\in S_l$ is bijective and
determines the natural parameterization of the segment\footnote{Or
rather, one of two possible parameterizations.} $\mathscr{S}_l$.
For its interior points, let us put
\begin{equation}
    \label{eq5.10}
    r(\widehat\pi):=|t^k_{\gamma l}(\widehat\pi)-\tau^k_{\gamma l}(0)|\in(0,\varepsilon_l),\qquad\widehat\pi\in S_l,
\end{equation}
where $\tau^k_{\gamma l}(0)$~-- the left end of cell
$\psi^k_{\gamma l}$. Then we extend the parametrization to the
``ends'' $\mathscr{S}_l\setminus S_l$ by continuity, taking for
them $r=0$ and $r=\varepsilon_l$ respectively. It follows from
\eqref{eq5.7} that $r(\widehat{\pi}_r)=r$ is satisfied for the
interior points.

Having chosen the parameterization in the way described above, we
obviously parameterize the remaining cells $\psi^{k'}_{\gamma l}$,
$k'\ne k$, and define functions $\tau^k_{\gamma}(r)$,
$r\in(0,\varepsilon_l)$, $k=1,\dots,n_{\gamma l}$, according to
\eqref{eq5.7}. They are also invariants of the eikonal algebra.

\subsection{Transformation}
\label{ss5.3}

The canonical form of \eqref{eq5.1} was obtained by
``reformatting'' the parametric representation of \eqref{eq4.8}.
Let us show how to arrive at ~it
by starting from the algebra $\mathfrak{E}^T_\Sigma$ itself and using its invariants. 

So we have the eikonals $E^T_\gamma$, $\gamma\in\Sigma$, and the
algebra $\mathfrak{E}^T_\Sigma$ generated by them.

\textsl{Step}~1.~Find the spectrum $\widehat{
\mathfrak{E}^T_\Sigma}$, equip it with the topology (Jacobson) and
separate in ~ it the segments and their components according to
\eqref{eq5.3}.

\textsl{Step}~2.~Find the spectra $\sigma_{\mathrm{ac}}(
E^T_\gamma)$ and introduce $\operatorname{int}\widehat{
\mathfrak{E}^T_\Sigma}$\enskip $\gamma$--coordinates. Determine
the cells $\psi^k_{\gamma l}$ by \eqref{eq5.9}. Parameterize the
segments by \eqref{eq5.10}.

\textsl{Step}~3.~For each $l$, choose representations
$\pi_r\colon \widehat{\mathfrak{E}^T_\Sigma}\to\mathbb M^{\varkappa_l}$, $\pi_r\,{\in}\,\widehat\pi_r\,{\in}\, S_l$, 
$r\in(0,\varepsilon_l)$, such that for all $e\in
\mathfrak{E}^T_\Sigma$\enskip the $\mathbb
M^{\varkappa_l}$--valued functions $\pi_r(e)$ are continuous at
$r\in(0,\varepsilon_l)$ and extend by continuity\footnote{This
choice is possible already by the existence of a canonical
representation.} on $[0,\varepsilon_l]$.

\textsl{Step}~4. Determine the eigenvalues $\tau^k_{\gamma l}(r)$
and eigenprojectors $P^k_{\gamma l}$ of matrices
$\pi_r(E^T_\gamma)$ (see \ \eqref{eq5.8}). Using them, define the
map
$$
\mathbf{U}\colon E^T_\gamma \mapsto \sum_{k=1}^{n_{\gamma
l}}\tau^k_{\gamma l}(\,{\cdot}\,)P_{\gamma l}^k\in \dot
C([0,\varepsilon_l];\mathbb M^{\varkappa_l}),\quad
l=1,\dots,\mathcal L,\quad \gamma\in\Sigma,
$$
and extend it from the ~generators--eikonals to the isomorphism of
the algebras ${\mathbf U}\colon
\mathfrak{E}^T_\Sigma\to\bigoplus\limits_{l=1}^\mathcal L \dot
C([0,\varepsilon_l];\mathbb M^{\varkappa_l})$.

The map of $\mathbf U$ drives $\mathfrak{E}^T_\Sigma$ up to
a~canonical form (or more exactly to~one of its versions: see the
comment after the theorem \ref{t2}).

All the results obtained relate to~shifted eikonals $\dot
E^T_\gamma$\enskip (see\ the \ref{con2} agreement). Their
reformulation for the original $E^T_\gamma=\dot
E^T_\gamma-P^T_\gamma$ is obvious and is actually reduced
to~replacing in~\eqref{eq5.1} the $\tau^k_{\gamma l}$ functions
with $\tau^k_{\gamma l}-1$.

\goodbreak 

\subsubsection*{Comments}

The map ${\mathbf V}\colon
\mathfrak{E}^T_\Sigma\to\bigoplus_{l=1}^\mathcal L \dot
C(\mathscr{S}_l;\mathbb M^{\varkappa_l})$, defined on the
generators by the relation
$$
\bigl(\mathbf{V}(E^T_\gamma)\bigr)(\widehat\pi):=
\sum_{k=1}^{n_{\gamma l}} \tau^k_{\gamma l}(r(\widehat\pi))
P_{\gamma l}^k,\qquad \widehat\pi\in\mathscr{S}_l,\quad
l=1,\dots,\mathcal L,\quad\gamma\in\Sigma,
$$
realizes the elements of the algebra $\mathfrak{E}^T_\Sigma$ as
matrix-valued functions on its spectrum and thus provides an
\textit{invariant functional model} of the eikonal algebra.
To~this model, with all its attributes --- cells $\psi^k_{\gamma
l}$, functions~$\tau^k_{\gamma l}$, lengths $\varepsilon_l$,
dimensions $\varkappa_l$ --- one can proceed by starting from any
isomorphic copy of the algebra $\mathfrak{E}^T_\Sigma$ and using
the procedure \textsl{Step}~1--4. This is important for the
inverse problem, since its data defines one such copy. The
question is~to what extent do these attributes determine the
structure of the graph $\Omega$.

The following observation may be useful in the inverse problem.
At points on the spectrum $\widehat{\mathfrak{E}^T_\Sigma}$ we
introduce the relation $\widehat\pi\stackrel{\gamma}{\sim_0}\widehat\pi'$
if among the $\gamma$-coordinates of these points there are matching ones,
i.e. \,\ $\tau^k_{\gamma l}(\widehat\pi)=\tau^{k'}_{\gamma l'}(\widehat\pi')$
are satisfied. Then let's define 
equivalence $\widehat\pi\sim\widehat\pi'$ if one can find points
$\widehat\pi_1,\dots,\widehat\pi_p\in\widehat{\mathfrak{E}^T_\Sigma}$
and vertices $\gamma_1,\dots, \gamma_{p+1}\in\Sigma$ such that
$\widehat\pi\stackrel{\gamma_1}\sim_0\widehat\pi_1\stackrel{\gamma_2}\sim_0\cdots
\stackrel{\gamma_p}\sim_0\widehat\pi_p\stackrel{\gamma_{p+1}}{\sim_0}\widehat\pi'$.
It can be shown that the factorization of the spectrum with
respect to the relation $\widehat\pi\sim\widehat\pi'$ identifies
only the points included in~clusters. Factorization transforms
$\widehat{\mathfrak{E}^T_\Sigma}$ into a ~space homeomorphic to
some graph. Examples show that the resulting graph is homeomorphic
to a wave-filled region $\Omega^T_\Sigma$ factorized by some
relation which has a simple geometric meaning.

In known examples \cite{6}, \cite{8}, the clusters appear when an
interior vertex of the graph $\Omega$ is overlapped by waves
coming from at least two boundary vertices. We assume that this is
a general fact. It is also an interesting question: can we
characterize the presence of cycles in~$\Omega^T_\Sigma$ in~terms
of the algebra $\mathfrak{E}^T_\Sigma$ (see \cite{6})? The
question is open.

\end{fulltext}

\end{document}